\documentclass[twoside]{article}

\usepackage[accepted]{aistats2025}
\usepackage{hyperref,natbib} 

\usepackage{times,tabularx,subfigure,booktabs,makecell,booktabs}
\usepackage{amssymb,amsfonts,amsmath,graphicx,cite,caption,wrapfig}
\usepackage{enumerate}
\usepackage[shortlabels]{enumitem}
\usepackage{multicol}
\usepackage{empheq}
\usepackage{color}
\usepackage{algorithm,algorithmic}
\usepackage{cases}
\usepackage[normalem]{ulem}

\allowdisplaybreaks
\let\underbrace\LaTeXunderbrace

\usepackage{tikz}

\newcommand{\Acal}{{\cal A}}

\newcommand{\Ecal}{{\cal E}}
\newcommand{\Fcal}{{\cal F}}

\newcommand{\Lcal}{{\cal L}}
\newcommand{\Mcal}{{\cal M}}

\newcommand{\Ocal}{{\cal O}}
\newcommand{\Pcal}{{\cal P}}

\newcommand{\Scal}{{\cal S}}

\newcommand{\Xcal}{{\cal X}}

\newcommand{\1}{{\mathbf{1}}}

\newtheorem{prop}{Proposition}
\newtheorem{lem}{Lemma}
\newtheorem{thm}{Theorem}
\newtheorem{cor}{Corollary}
\newtheorem{definition}{Definition}

\newtheorem{assump}{Assumption}
\newtheorem{remark}{Remark}

\newtheorem{example}{Example}

\newcommand{\qed}{\hfill $\blacksquare$}

\begin{document}

\runningtitle{Learning in Herding Mean Field Games: Single-Loop Algorithm with Finite-Time Convergence Analysis}

\twocolumn[

\aistatstitle{Learning in Herding Mean Field Games: \\ Single-Loop Algorithm with Finite-Time Convergence Analysis}

\aistatsauthor{Sihan Zeng \And Sujay Bhatt \And  Alec Koppel \And Sumitra Ganesh}

\aistatsaddress{JPMorgan AI Research, United States} 
]


\begin{abstract}
We consider discrete-time stationary mean field games (MFG) with unknown dynamics and design algorithms for finding the equilibrium with finite-time complexity guarantees. Prior solutions to the problem assume either the contraction of a mean field optimality-consistency operator or strict weak monotonicity, which may be overly restrictive. In this work, we introduce a new class of solvable MFGs, named the ``fully herding class'', which expands the known solvable class of MFGs and for the first time includes problems with multiple equilibria. We propose a direct policy optimization method, Accelerated Single-loop Actor Critic Algorithm for Mean Field Games (ASAC-MFG), that provably finds a global equilibrium for MFGs within this class, under suitable access to a single trajectory of Markovian samples. Different from the prior methods, ASAC-MFG is single-loop and single-sample-path. We establish the finite-time and finite-sample convergence of ASAC-MFG to a mean field equilibrium via new techniques that we develop for multi-time-scale stochastic approximation. We support the theoretical results with illustrative numerical simulations.

When the mean field does not affect the transition and reward, an MFG reduces to a Markov decision process (MDP) and \texttt{ASAC-MFG} becomes an actor-critic algorithm for finding the optimal policy in average-reward MDPs, with a sample complexity matching the state-of-the-art. Previous works derive the complexity assuming a contraction on the Bellman operator, which is invalid for average-reward MDPs. We match the rate while removing the untenable assumption through an improved Lyapunov function.
\end{abstract}

\section{INTRODUCTION}
The mean field game (MFG) framework, introduced in \citet{huang2006large,lasry2007mean}, provides an infinite-population approximation to the $N$-agent Markov game with a large number of homogeneous agents. 
It addresses the increasing difficulty in solving Markov games as $N$ scales up and finds practical applications in many domains, including resource allocation \citep{li2020resource,mao2022mean}, 
wireless communication \citep{narasimha2019mean,jiang2019mean}, and the management of power grids \citep{alasseur2020extended,zhang2021vehicle}.

A mean field equilibrium (MFE) describes the notion of solution in an MFG, and is a pair consisting of a policy and a mean field: the policy performs optimally in a Markov decision process (MDP) determined by the mean field, whereas the mean field is the induced stationary distribution of the states when every agent in the infinite population adopts the policy. In the discrete-time setting without explicit knowledge of the environment dynamics, reinforcement learning (RL) provides an important tool for finding an MFE using samples of state transitions and rewards. 
Currently two classes of MFGs are known to be solvable by RL with finite-time convergence guarantees, which we now review.

\noindent
\textbf{Existing classes of solvable MFGs}\vspace{0.1cm}\\
$\bullet\,\,$ \textit{Contractive MFG}: 
In the context of MFGs, the optimality operator maps a mean field distribution to the optimal policy in the MDP determined by the mean field, whereas the consistency operator returns the induced mean field for a given policy. A series of recent works \citep{guo2019learning,xie2021learning,anahtarci2023q,mao2022mean,zaman2023oracle,yardim2023policy} make the assumption that the composition of the optimality and consistency operator is a contractive mapping, or enforce the contraction through an entropy regularization. Fixed-point iteration methods (such as the extension of Q-learning to MFGs) are developed and analyzed by exploiting the contraction. A related concurrent work \citep{zhang2024stochastic} relies on a ``sufficiently Lipschitz MDP'' assumption which plays a similar role of ensuring that the MFE is the unique fixed point of a contractive mapping.
However, as pointed out in \citet{yardim2024mean}, the contraction assumption usually only holds when an impractically large regularization is added. Since the policy at the regularized equilibrium quickly approaches a uniform distribution as the regularization weight increases, solving such a regularized problem is usually \textit{uninformative} about the original game. It is also worth noting that contractive MFGs can only have a unique equilibrium, while multiple equilibria usually exist in MFGs used to model practical problems \citep{nutz2020convergence,dianetti2024multiple}.\\
$\bullet\,\,$ \textit{Strictly Weak Monotone MFG}: 
Another line of work relies on the assumption that the MFG is strictly weak monotone \citep{perrin2020fictitious,perolat2021scaling,geist2021concave,zhang2024learning}. 
This condition intuitively means that the agent is ``\textit{discouraged from taking similar state-action pairs as the rest of the population}'' and can be interpreted as an ``aversion to crowded areas'' \citep{perolat2021scaling}.
Strictly weakly monotone MFGs are again limited in that they need to have a unique equilibrium.

\begin{figure}
	\begin{center}
		\includegraphics[width=.4\textwidth]{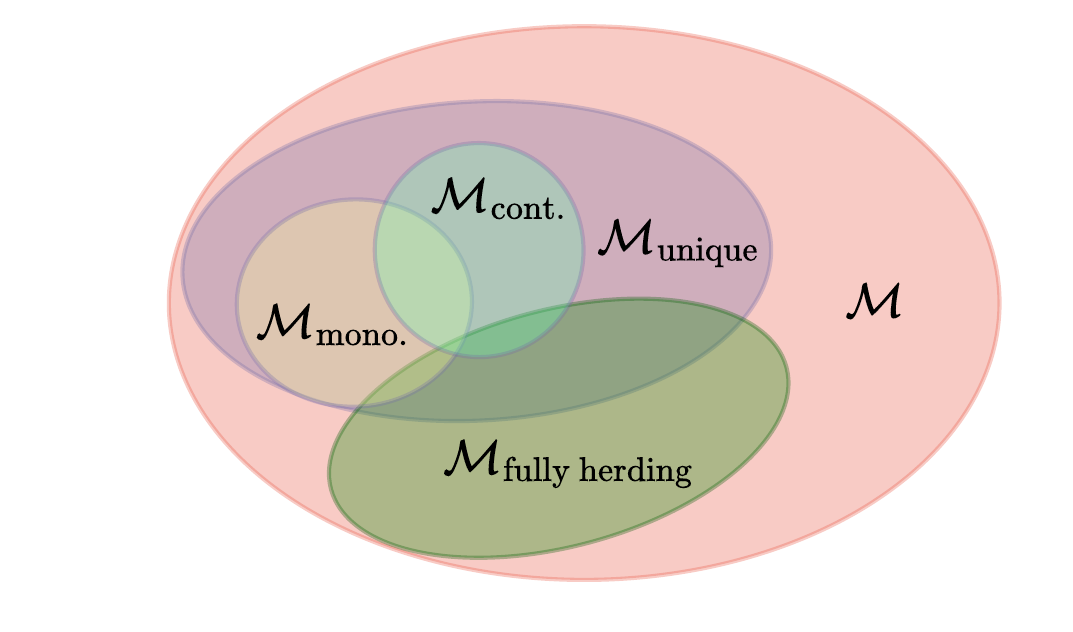}
	\end{center}
    \vspace{-15pt}
	\caption{{ \small{$\mathcal{M}$ denotes the class of all MFGs. $\mathcal{M}_{\text{cont.}}$ and $\mathcal{M}_{\text{mono.}}$ are the MFG classes satisfying contraction and strict weak monotonicity, and they are subsets of $\mathcal{M}_{\text{unique}}$ which is the class of MFGs having a unique equilibrium. The proposed algorithm,~\texttt{ASAC-MFG}, solves MFGs in 
    $\mathcal{M}_{\text{fully herding}}$~[cf.~Def.~\ref{def:h-mfg}].}}}
    \vspace{-10pt}
	\label{fig:Venn}
\end{figure}

Figure.~\ref{fig:Venn} gives an overview of the class of solvable MFGs. The aim of this paper is to expand the class of  solvable MFGs and to design an easily implementable, provably convergent, and efficient algorithm for solving MFGs. We summarize our main contributions below.

\subsection{Main Contributions}
\begin{enumerate}
%
%
    \item We introduce the $\kappa$-herding class of MFGs and the fully herding class which is a special case with $\kappa=0$, and show that the fully herding MFGs are perfectly solvable. It is known from \citet{yardim2024mean} that solving general MFGs (even with Lipschitz transition kernel and reward function) is a PPAD-complete problem, conjectured to be computationally intractable \citep{daskalakis2009complexity}. Notably, we show that the fully herding class contains MFG instances admitting more than one equilibrium. Such MFGs do not satisfy either contraction or strict weak monotonicity and are not previously known to be solvable. In this sense, our work complements and expands on the finding of \citet{yardim2024mean} and enlarges the class of solvable MFGs. As pointed out in \citet{guo2024mf,cui2021approximately}, MFGs with multiple equilibria are very common but significantly more challenging to solve than those with a unique equilibrium. 

\item We propose a single-loop, single-sample-path policy optimization algorithm \texttt{ASAC-MFG} for finding the equilibrium for MFGs in the herding class, and explicitly characterize its finite-time and finite-sample complexity. 
For MFGs in the fully herding class ($\kappa=0$), \texttt{ASAC-MFG} converges to a global MFE with a rate of $\widetilde{\Ocal}(k^{-1/4})$; for general herding MFGs with $\kappa>0$, it converges to a $\sqrt{\kappa}-$approximate equilibrium at the same rate. As our algorithm draws exactly one sample in each iteration, the finite-time complexity translates to a finite-sample complexity of the same order.
We note that any Lipschitz MFG can be shown to belong to the $LL_V$-herding class due to a simple Lipschitz continuity bound, where $L$ and $L_V$ are the Lipschitz constants introduced later. As a result, for MFGs not in the fully herding class, \texttt{ASAC-MFG} is still stable and can provably find a (despite non-ideal) $\Ocal(\sqrt{LL_V})$ approximate solution. In comparison, the existing fixed-point iteration algorithms based on the contractive MFG assumption may in theory exhibit arbitrarily unstable behavior when the assumption fails.

\item
The single-loop and single-sample-path structure makes our algorithm easily implementable. While single-loop and single-sample-path algorithms are widely used to solve RL and games in practice due to simplicity, their theoretical understanding is not as complete as their nested-loop counterparts. Specifically for MFGs, there does not currently exist a finite-time convergent single-loop and single-sample-path algorithm (see Table 1).
Our work fills in the important gap.
Note that our ability to make the algorithm single-loop and single-sample-path is not due to the herding condition. The analysis of \texttt{ASAC-MFG} is a technical innovation made by adapting the advances in accelerated two-time-scale stochastic approximation \citep{zeng2024fast} and generalizing them to a three-time-scale setting. We discuss the innovation in detail in Section~\ref{sec:algorithm}.

\item We can regard a Markov decision process (MDP) as a degenerate MFG in which the transition kernel and reward are independent of the mean field. Recognizing this connection, we note that a simplified version of the proposed method becomes an actor-critic algorithm in an average-reward MDP and is guaranteed to converge to a stationary point of the policy optimization objective with rate $\widetilde{\Ocal}(1/\sqrt{k})$. This matches the state-of-the-art complexity of the actor-critic algorithm \citep{chen2024finite}. 
Existing works derive the complexity assuming a contraction  on the Bellman operator, which is invalid for average-reward MDPs. We maintain the rate but remove this assumption through introducing an alternative Lyapunov function.

\end{enumerate}

\begin{table*}[!ht]
\centering
\setlength{\tabcolsep}{5pt}
\small\sl
\begin{tabular}{ccccc}
        \toprule
        & \makecell{Assumption} & \makecell{Single \\Sample Path} & \makecell{Single\\ Loop} & \makecell{Sample Complexity}\\
        \midrule
        \citet{guo2019learning} & Contraction & No & No & Regularization Dependent\\
        \citet{xie2021learning} & Contraction & Yes* & Yes* & $\widetilde{\Ocal}(\epsilon^{-5})$, regularized solution\\
        \citet{mao2022mean} & Contraction & No & No & $\widetilde{\Ocal}(\epsilon^{-5})$, regularized solution\\
        \citet{zaman2023oracle} & Contraction & Yes & No & $\widetilde{\Ocal}(\epsilon^{-4})$, regularized solution\\
        \citet{yardim2023policy} & Contraction & Yes & No & $\widetilde{\Ocal}(\epsilon^{-2})$, regularized solution\\
        \citet{zhang2024learning} & Strict Weak Monotonicity & No & No & $\widetilde{\Ocal}(\epsilon^{-4})$, original solution\\
        \textbf{This Work} & \textbf{Fully Herding Class ($\kappa=0$)} & \textbf{Yes} & \textbf{Yes} & \textbf{\makecell{$\widetilde{\Ocal}(\epsilon^{-4})$, original solution}}\\
        \textbf{This Work} & \textbf{Partially Herding Class($\kappa>0$)} & \textbf{Yes} & \textbf{Yes} & \textbf{\makecell{$\widetilde{\Ocal}(\epsilon^{-4})$, $\sqrt{\kappa}$-optimal solution}}\\
        \bottomrule
        \bottomrule
        \end{tabular}
\vspace{-3pt}
\caption{Assumption, structure, and complexity of existing algorithms with finite-sample analysis.
* The algorithm in \citet{xie2021learning} is single-loop and single-sample-path under an oracle that returns the stationary distribution of states for any $\pi,\mu$. \citet{mao2022mean} also relies on such an oracle. Our work, in comparison, is oracle-free. 
}
\vspace{-3pt}
\end{table*}

\subsection{Related Work}\label{sec:literature}
The classic works on MFGs study the continuous-time setting where the equilibrium point simultaneously satisfies a Hamilton–Jacobi–Bellman equation on the optimality of the policy and a Fokker–Planck equation that describes the dynamics of the mean field, and have proposed optimal control techniques that provably find the solution \citep{huang2006large,huang2007large,lasry2007mean}. 
In discrete time, MFGs can be considered a generalization of MDPs and are widely solved using RL. Among the latest representative works, \citet{yang2018learning,carmona2021linear,perolat2021scaling} build upon policy optimization and \citet{anahtarci2020value,angiuli2022unified,angiuli2023convergence} consider valued-based methods. The algorithms proposed in these works, however, either do not come with convergence analysis or are only shown to converge asymptotically.

Finite-time convergent algorithms are recently developed under the contraction assumption or strict weak monotonicity, as discussed earlier in the section. 
Table~1 highlights our contribution in terms of assumptions, algorithm structure, and sample complexity. 
Most papers listed introduce a large regularization to the MFG and establish the convergence to a regularized equilibrium. Notably, if \citet{xie2021learning,mao2022mean,zaman2023oracle,yardim2023policy} could choose the regularization weight freely (note that they actually could not since the contraction condition only holds when the weight is sufficiently large as previously discussed), they can solve the original unregularized game by making the weight small enough. 
Solving the original game requires doubled complexities, i.e., become $\widetilde{\Ocal}(\epsilon^{-10})$, $\widetilde{\Ocal}(\epsilon^{-8})$, or $\widetilde{\Ocal}(\epsilon^{-4})$ to the original solution.

Among works based on the strict weak monotonicity, \citet{zhang2024learning} proposes a mirror descent algorithm with a complexity matching that of \texttt{ASAC-MFG}, but is less convenient to implement due to its nested-loop structure and the requirement to pre-generate and store offline samples.
\citet{perrin2020fictitious} proposes a continuous-time algorithm and \citet{geist2021concave} studies a deterministic gradient algorithm not based on samples, making their complexities not directly comparable.

Finally, we note the separate line of works \citep{guo2024mf,mandal2023performative} that reformulate the MFG policy optimization problem as a constrained program with convex constraints and a bounded objective. The simple projected gradient descent algorithm provably solves the constrained program, leading to a solution of the MFG. However, a finite-time convergence guarantee is not established, unless again a sufficiently large regularization is added.


\section{FORMULATION}\label{sec:formulation}
We study MFGs in the stationary and infinite-horizon average-reward setting. An MFG is characterized by~$M=(\Scal,\Acal,\Pcal,r)$, 
where $\Scal$ and $\Acal$ denote the \textit{finite} state and action spaces. 
From the perspective of a single representative agent, the state transition depends not only on its own action but also on the aggregate behavior of all other agents. 
The aggregate behavior is described by the mean field~$\mu\in\Delta_{\Scal}$\footnote{We use $\Delta_{\Scal}$ and $\Delta_{\Acal}$ to denote the probability simplex over the state and action spaces.}, which measures the percentage of population in each state.
The transition kernel of an MFG is represented by $\Pcal:\Scal\times\Acal\times \Delta_{\Scal}\rightarrow\Delta_{\Scal}$, where $\Pcal^{\mu}(s'\mid s,a)$ denotes the probability that the state of the representative agent transitions from $s$ to $s'$ when it takes action $a$ and mean field is $\mu$. The mean field also affects the reward function $r:\Scal\times\Acal\times\Delta_{\Scal}\rightarrow[0,1]$ -- the agent receives reward $r(s,a,\mu)$ when it takes action $a$ in state $s$ under mean field $\mu$. The agent \textit{does not} observe the mean field, and takes actions according to policy $\pi:\Scal\rightarrow\Delta_{\Acal}$, which can be represented as a table~$\Delta_{\Acal}^{\Scal}\subset\mathbb{R}^{|\Scal|\times|\Acal|}$.


Given a policy $\pi$ and mean field $\mu$, the sequentially generated states form a Markov chain with transition matrix $P^{\pi,\,\mu}\in\mathbb{R}^{|\Scal|\times|\Scal|}$, with $P^{\pi,\,\mu}_{s',s} = \textstyle\sum_{a\in\Acal}\Pcal^{\mu}(s'\mid s,a)\pi(a\mid s).$
We denote by $\nu^{\pi,\,\mu}\in\Delta_{\Scal}$ the stationary distribution of the Markov chain, which is the right singular vector of $P^{\pi,\,\mu}$ associated with singular value~$1$, i.e., $\nu^{\pi,\,\mu} = P^{\pi,\,\mu} \nu^{\pi,\,\mu}.$
When the mean field is $\mu$ and the agent generates actions according to $\pi$, the agent can expect to collect the cumulative reward $J(\pi,\mu)$
\begin{align}
J(\pi,\mu)
&\triangleq \textstyle\lim_{T\rightarrow\infty}\frac{1}{T}\mathbb{E}_{\pi,\Pcal^{\mu}}[ \sum_{t=0}^{T-1}r(s_t, a_t, \mu ) \mid s_0]\notag\\
&=\mathbb{E}_{s\sim \nu^{\pi,\,\mu},\,a\sim\pi(\cdot\mid s)}[r(s, a, \mu)].\label{eq:def_J}
\end{align}
If the mean field were fixed to a given $\mu$, the goal of the agent would be to find a policy $\pi$ that maximizes $J(\pi,\mu)$. However, when every agent in the infinite population follows the same policy as the representative agent, the mean field evolves as a function of $\pi$. We use $\mu^{\star}:\Delta_{\Acal}^{\Scal}\rightarrow\Delta_{\Scal}$ to denote the mapping from a policy to the induced mean field, which is the stationary distribution of states when the infinite number of players in the game all adopt policy $\pi$, i.e, $\mu^{\star}(\pi)= \nu^{\pi,\,\mu^{\star}(\pi)}$. 
The goal of the representative agent in an MFG is to find a policy optimal under the mean field induced by the policy. Mathematically, we want to find a pair of policy and mean field $(\bar{\pi}, \bar{\mu})$, known to always exist  \citep{cui2021approximately},  
as the solution to the system
\begin{numcases}{}
J(\bar{\pi},\bar{\mu})\geq J(\pi,\bar{\mu}),\quad\forall\pi\label{eq:obj:optimality}\\
\bar{\mu} = \mu^{\star}(\bar{\pi}).\label{eq:obj:consistency}
\end{numcases}
We assume that the induced mean field $\mu^{\star}(\pi)$ is unique for any $\pi$. This does not imply that the MFE $(\bar{\pi}, \bar{\mu})$ is unique.

\begin{definition}
 [$\epsilon-$MFE]
\label{def:epsilon_MFE}
The pair of policy and mean field $(\pi,\mu)$ is an $\epsilon$-Mean Field Equilibrium (MFE) if 
\begin{align}
J(\pi',\mu)-J(\pi,\mu)\leq\epsilon,\forall \pi',\text{  and  }\|\mu-\mu^{\star}(\pi)\|\leq\epsilon.\label{def:epsilon_MFE:eq1}
\end{align}
\end{definition}

If the pair $(\pi,\mu)$ satisfies \eqref{def:epsilon_MFE:eq1} with $\epsilon=0$, it is obviously an exact MFE as a solution to \eqref{eq:obj:optimality}-\eqref{eq:obj:consistency}.
The definition characterizes what it mathematically means when we say $(\pi,\mu)$ is an approximate solution to the MFG. We also make use of the differential value function $V^{\pi,\,\mu}\in\mathbb{R}^{|\Scal|}$ to quantify the relative value of each initial state
\begin{align*}
    V^{\pi,\,\mu}(s) 
    \hspace{-1pt}\triangleq\hspace{-1pt} \mathbb{E}_{\pi,\Pcal^{\mu}}\big[\textstyle\sum_{t=0}^{\infty} (r(s_t, a_t, \mu )-J(\pi,\mu) ) \hspace{-1pt} \mid\hspace{-1pt} s_0\hspace{-1pt}=\hspace{-1pt}s\big].
\end{align*}

\subsection{Herding MFGs}

In this section, we identify a novel class of solvable MFGs named the ``fully herding MFGs'' and show that it contains instances that do not satisfy the contraction or strict weak monotonicity conditions.

\begin{definition}[Herding MFG] \label{def:h-mfg}
Given $\kappa\geq0$, an MFG $M$ is in the $\kappa$-herding class, denoted as $M\in\Mcal_{\kappa}$, if there exists a constant $0<\rho<\infty$ such that $M$ satisfies the following ``herding condition'' for all $\pi,\pi'$
\begin{align}
&J(\pi,\mu^{\star}(\pi))-J(\pi',\mu^{\star}(\pi'))\notag\\
& \leq \rho\Big(J(\pi,\mu^{\star}(\pi))-J(\pi',\mu^{\star}(\pi))\Big)+\kappa\|\pi-\pi'\|.\label{eq:def_Delta}
\end{align}
\end{definition}
The MFG $M$ is further said to be a \textbf{fully herding MFG} if $M\in\Mcal_{\text{fully herding}}\triangleq\Mcal_0$ and a partially herding MFG if $M\in\Mcal_{\kappa}$ for some $0<\kappa<\infty$.

Conceptually, in MFGs with a small or zero $\kappa$, the representative agent receives a higher reward by ``following the crowd'' or displaying a ``herding'' behavior. 
The constant $\kappa$ quantifies the difficulty of an MFG.
We will show that the proposed algorithm \texttt{ASAC-MFG} provably finds a global equilibrium for fully herding MFGs ($\kappa=0$), whereas for MFGs in the partial herding class we can at least provably find an $\Ocal(\sqrt{\kappa})$-MFE (in the sense of Definition ~\ref{def:epsilon_MFE}).
Notably, all Lipschitz MFGs satisfy the herding condition with $\rho=1$ and $\kappa=LL_V$ in the worst case, where $L,L_V$ are Lipschitz constants introduced later. This means that the proposed method exhibits stable behavior and converges to an approximate (though non-ideal) solution in any MFG, whereas the existing algorithms lose stability guarantees when contraction/monotonicity fails to hold. Next, we introduce a few examples to clarify the herding condition.


\begin{figure}
	\begin{center}
		\includegraphics[width=.48\textwidth]{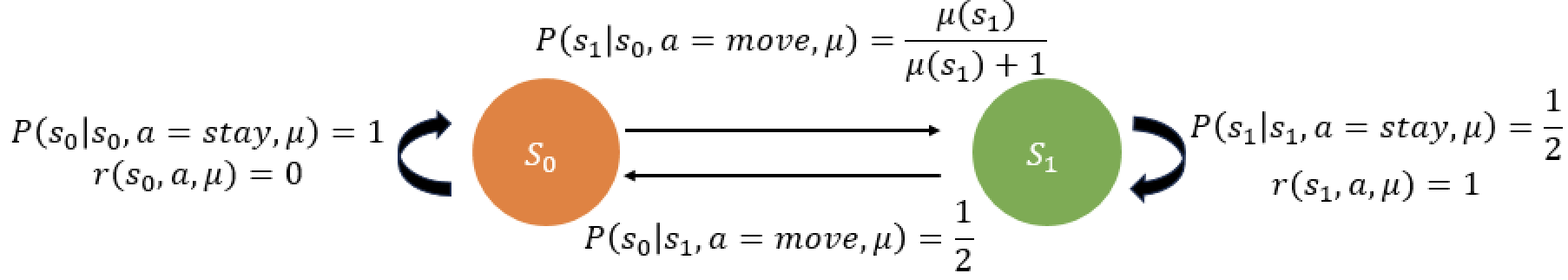}
	\end{center}
    \vspace{-10pt}
	\caption{Illustration of Example~\ref{example:twostate}}
    \vspace{-5pt}
	\label{fig:threestate}
\end{figure}

\begin{example}\label{example:assumption_Delta}	
Consider an MFG which satisfies $\forall \mu,\mu',s,a,s'$
\begin{align*}
r(s,a,\mu)\hspace{-1pt}=q\big(\mu(s)\big)^p\text{ and }
\Pcal(s'\hspace{-2pt}\mid \hspace{-2pt} s,a,\mu)\hspace{-1pt}=\hspace{-1pt}\Pcal(s'\hspace{-2pt}\mid\hspace{-2pt} s,a,\mu'),
\end{align*}
for some scalars $p,q>0$, i.e., the transition probability kernel is independent of the mean field. The MFG belongs to the fully herding class ($\kappa=0,\rho=p+1$) and is perfectly solvable by \texttt{ASAC-MFG}.
\end{example}
Consider an MFG instance satisfying the conditions in Example~\ref{example:assumption_Delta} with $|\Scal|\hspace{-2pt}=\hspace{-2pt}|\Acal|\hspace{-2pt}=\hspace{-2pt}2$, $q\hspace{-2pt}=\hspace{-2pt}1,p\hspace{-2pt}=\hspace{-2pt}1$, where the transition kernel is such that in either state $s\in\{s_1,s_2\}$, the action $a_1$ (resp. $a_2$) leads the next state to $s_1$ (resp. $s_2$) with probability $p=3/4$. There exist an infinite number of equilibria in this MFG. They occur at policies $\bar{\pi}_1$, $\bar{\pi}_2$ such that $\forall s$
	\[\bar{\pi}_1(a\mid s)\hspace{-1pt}=\hspace{-1pt}\begin{cases}
		1, &\text{if $a=a_1$}\\
		0, &\text{if $a=a_0$}
	\end{cases}\;\,\bar{\pi}_2(a\mid s)\hspace{-1pt}=\hspace{-1pt}\begin{cases}
		0, &\text{if $a=a_1$}\\
		1, &\text{if $a=a_0$}
	\end{cases}
	\]
with the induced mean field $\bar{\mu}_1=[3/4,1/4]^{\top}$, $\bar{\mu}_2=[1/4,3/4]^{\top}$, and at all policies that induce the mean field $[1/2,1/2]^{\top}$ (such as $\bar{\pi}_3(a\mid s)=1/2$ for all $s,a$). The MFE is not unique, implying that MFGs in Example~\ref{example:assumption_Delta} generally do not satisfy the contraction or strict weak monotonicity assumptions. 
The derivation of the mean field equilibria and the proof that MFGs in Example~\ref{example:assumption_Delta} are in the fully herding class can be found in Appendix~\ref{sec:details_example}.
We show another example of fully herding MFG in which the transition kernel is dependent on the mean field. 

\begin{example}\label{example:twostate}
Consider an MFG with two states, $s_0$ and $s_1$, and two actions ``move'' and ``stay''. From $s_0$, we can move to $s_1$ with probability $\frac{\mu_{s_1}}{1+\mu(s_1)}$ or stay in $s_0$ with probability $\frac{1}{1+\mu(s_1)}$ by taking action ``move''. Taking the action ``stay'' in state $s_0$ makes us stay in the state with probability 1. When in state $s_1$, we transition to $s_0$ or $s_1$ for the next state each with probability $1/2$ under any action.
In state $s_1$, we collect a reward of 1 regardless of action and mean field. 
The reward in any other situation is 0. 
This MFG has an infinite number of equilibria -- the optimal action to take in state $s_0$ is always ``move'', but any policy can be taken in state $s_1$.
This MFG satisfies \eqref{eq:def_Delta} with $\kappa=0,\rho=2$ but not the contraction or strict weak monotonicity condition.
\end{example}

Note that Examples~\ref{example:assumption_Delta} and \ref{example:twostate} are instances within the set $\Mcal_{\text{fully herding}}\backslash\Mcal_{\text{unique}}$ in Figure~\ref{fig:Venn} and not known to be provably solvable in the existing literature.

The herding class connects to practical problems in real life. An important example is crowd motion with a reward function that models ``attraction to the mean''. Specifically, the reward of the representative player is high for the states that match the mean field, which may take the form of $r(s,a,\mu)=\mu(s)$ and fall under Example 1. (See Section 3.1.1 of \citet{dayanikli2024machine}. Their reward for the ``attraction to the mean'' setting can be regarded as a variant of $r(s,a,\mu)=\mu(s)$ where the states are discretized continuous numbers.) Global carbon emission can be understood as a specific example of crowd motion with an ``attraction to the mean'' reward, in which the representative agent is a country that needs to determine its emission level. The agent has a tendency of following the emission level of the population (in this case, all other countries), since if it produces emissions lower than the average, ``it behaves as an opportunity cost of not producing more'' \citep{dayanikli2024multi}; on the other hand, if it produces emission higher than the average, ``it behaves as a reputation cost for polluting more'' \citep{dayanikli2024multi}.

\begin{algorithm}[!ht]
\caption{Accelerated Single-loop Actor Critic Algorithm for Mean Field Games (\texttt{ASAC-MFG})}
\label{alg:main}
\begin{algorithmic}[1]
\STATE{\textbf{Initialize:} policy parameter $\theta_0$, value function estimate $\hat{V}_0,\hat{J}_0$, mean field estimate $\hat{\mu}_0\in\Delta_{\Scal}$, gradient/operator estimates $f_0=0\in\mathbb{R}^{|\Scal||\Acal|},g_0^V=0\in\mathbb{R}^{|\Scal|},g_0^J=0\in\mathbb{R},h_0=0\in\mathbb{R}^{|\Scal|}$, initial state $s_0$ 
}
\FOR{iteration $k=0,1,2,...$}
\STATE{Take action $a_k \sim \pi_{\theta_k}(\cdot\mid s_k)$. Observe $r(s_k,a_k,\hat{\mu}_k)$ and $s_{k+1}\sim \Pcal^{\hat{\mu}_k}(\cdot\mid s_k,a_k)$}
\STATE{Policy (actor) update:
\vspace{-5pt}
\begin{align}
\theta_{k+1} = \theta_k + \alpha_k f_k. \label{alg:main:actor}
\end{align}
}
\vspace{-15pt}
\STATE{Mean field update:
\vspace{-5pt}
\begin{align}
\hat{\mu}_{k+1}=\Pi_{\Delta_{\Scal}} \big(\hat{\mu}_{k}+\xi_k h_k\big).\label{alg:main:meanfield}
\end{align}
}
\vspace{-10pt}
\STATE{Value function (critic) update:
\vspace{-5pt}
\begin{align}
\begin{aligned}
\hat{V}_{k+1} &= \Pi_{B_V}(\hat{V}_{k} + \beta_k g_k^V),\\
\hat{J}_{k+1}&=\Pi_{[0,1]}(\hat{J}_k+\beta_k g_k^J).
\end{aligned}\label{alg:main:critic}
\end{align}
}
\vspace{-5pt}
\STATE{Gradient/Operator estimate update:
\vspace{-5pt}
\begin{align*}
\begin{aligned}
f_{k+1}&=(1-\lambda_k)f_{k}+\lambda_k\nabla\log\pi_{\theta_k}(a_k\mid s_k)\\
&\hspace{5pt}\times(r(s_k,a_k,\hat{\mu}_k)+\hat{V}_k(s_{k+1})-\hat{V}_k(s_{k})),\\
g_{k+1}^V&=(1-\lambda_k)g_{k}^V+\lambda_k e_{s_k}\\
&\hspace{5pt}\times(r(s_k,a_k,\hat{\mu}_k)\hspace{-2pt}-\hspace{-2pt}\hat{J}_k\hspace{-2pt}+\hspace{-2pt}\hat{V}_k(s_{k+1})\hspace{-2pt}-\hspace{-2pt}\hat{V}_k(s_k)),\\
g_{k+1}^J&=(1-\lambda_k) g_{k}^J+\lambda_k c_J( r(s_k,a_k,\hat{\mu}_k)-\hat{J}_k),\\
h_{k+1}&=(1-\lambda_k)h_k+\lambda_k(e_{s_k}-\hat{\mu}_k).
\end{aligned}
\end{align*}
}
\vspace{-10pt}
\ENDFOR
\end{algorithmic}
\end{algorithm}

\section{ALGORITHM}\label{sec:algorithm}
Our algorithm solves MFGs from the perspective of direct policy optimization.
As we do not directly deal with the mean field optimality-consistency operator, we bypass the need to assume that it is contractive.
We see from \eqref{eq:obj:optimality} that if the optimal policy under $\bar{\mu}$ were unique and we knew $\bar{\mu}$, we could easily find $\bar{\pi}$ through policy optimization with the mean field fixed to $\bar{\mu}$. 
On the other hand, if we knew the equilibrium policy $\bar{\pi}$, we could obtain $\bar{\mu}$ by finding $\mu^{\star}(\bar{\pi})$. However, we do not know either $\bar{\pi}$ or $\bar{\mu}$ in reality and therefore consider the approach of simultaneous learning. 
Specifically, we maintain a parameter $\theta$ that encodes the policy $\pi_{\theta}$,
and an iterate $\hat{\mu}$ to estimate the mean field induced by the current policy, and improve $\theta$ and $\hat{\mu}$ with respect to each other by iteratively taking the steps
\begin{align}
\theta_{k+1} = \theta_{k} + \alpha_k \nabla_{\theta}J(\pi_{\theta_k},\hat{\mu}_k),\quad\hat{\mu}_{k+1}= \mu^{\star}(\pi_{\theta_k})
\label{eq:wishful_update_eq1}
\end{align}
where $k$ indexes the iteration and $\alpha_k$ is a step size.
By the policy gradient theorem \citep{sutton1999policy}, a closed-form expression for 
$\nabla_{\theta} J(\pi_{\theta},\mu)$ is
\begin{align}
\nabla_{\theta} J(\pi_{\theta},\mu) 
&=\mathbb{E}_{s\sim \nu^{\pi_{\theta},\mu},a\sim\pi_{\theta}(\cdot\mid s),s'\sim\Pcal^{\mu}(\cdot\mid s,a)}\notag\\
&\hspace{-35pt}\Big[(r(s,a,\mu)\hspace{-1pt}+\hspace{-1pt}V^{\pi_{\theta},\mu}(s')\hspace{-1pt}-\hspace{-1pt}V^{\pi_{\theta},\mu}(s))\nabla_{\theta}\log\pi_{\theta}(a\mid s)\Big].\notag
\end{align}
Without loss of generality, we use a softmax parameterization, 
i.e., the parameter $\theta\in\mathbb{R}^{|\Scal||\Acal|}$ represents the policy as\looseness=-1
\vspace{-2pt}
\[\pi_{\theta}(a \mid s)=\frac{\exp (\theta(s, a))}{\sum_{a' \in \Acal} \exp(\theta(s, a'))}.\]
\vspace{-10pt}

In large and/or unknown environments in the real life, performing \eqref{eq:wishful_update_eq1} poses challenges. The updates require the knowledge of $\mu^{\star}(\pi_{\theta_k})$ and value function $V^{\pi_{\theta_k},\,\mu^{\star}(\pi_{\theta_k})}$. Neither of these quantities can be exactly calculated without the exact knowledge of the transition model, which is usually unavailable. We propose learning $\mu^{\star}(\pi_{\theta_k})$ and $V^{\pi_{\theta_k},\,\mu^{\star}(\pi_{\theta_k})}$ simultaneously, with the policy and mean field iterates updated using a continuous path of samples from the MFG. 
We recognize that for any $\theta$
\begin{align}
&\mu^{\star}({\pi_{\theta}})
=\textstyle\lim_{T\rightarrow\infty}\frac{1}{T}\mathbb{E}_{\pi_{\theta},\Pcal^{\mu^{\star}({\pi_{\theta}})}}[\sum_{t=0}^{T-1}e_{s_t}\mid s_0],\label{eq:mu_rewrite}
\end{align}
where $e_s\in\mathbb{R}^{|\Scal|}$ is the indicator vector whose entry $s'$ is 1 if $s'=s$ and 0 otherwise.
To solve Eq.~\eqref{eq:mu_rewrite} with multi-time-scale stochastic approximation, we carry out
\begin{align}
\hat{\mu}_{k+1}=\hat{\mu}_{k}+\xi_{k}(e_{s_k}-\hat{\mu}_{k})\label{eq:mu_online}
\end{align}
iteratively for some step size $\xi_{k} \gg \alpha_k$. Due to the difference in time scales (step size), $\hat{\mu}_k$ becomes an increasingly accurate estimate of $\mu^{\star}({\pi_{\theta_k}})$ as the iterations proceed.

We know that $V^{\pi_{\theta_k},\,\hat{\mu}_k}$ satisfies the Bellman equation
\begin{align}
V^{\pi_{\theta_k},\,\hat{\mu}_k}(s)&= \textstyle\sum_{a}\pi_{\theta_k}(a\mid s)r(s,a,\hat{\mu}_k)-J(\pi_{\theta_k},\hat{\mu}_k) \notag\\
&\hspace{20pt}\textstyle+ \sum_{s'}P^{\pi_{\theta_k},\,\hat{\mu}_k}_{s',s}V^{\pi_{\theta_k},\,\hat{\mu}_k}(s'),\;\forall s.\label{eq:Bellman}
\end{align}
We introduce an auxiliary variable $\hat{V}\in\mathbb{R}^{|\Scal|}$ to estimate $V^{\pi_{\theta_k},\,\hat{\mu}_k}$, again by stochastic approximation.
The following update solves \eqref{eq:Bellman} under proper choices of $\beta_k$
\begin{align}
\hat{V}_{k+1}(s_k)&=\hat{V}_{k}(s_k)\label{eq:V_online}\\
&\hspace{-10pt}+\beta_{k}\big(r(s_k,a_k,\hat{\mu}_k)-\hat{J}_{k}+\hat{V}_{k}(s_{k+1})-\hat{V}_{k}(s_k)\big),\notag
\end{align}
where the unknown $J(\pi_{\theta},\mu^{\star}(\pi_{\theta}))$ is replaced with an estimate that itself is iteratively refined
\begin{align}
\hat{J}_{k+1}=\hat{J}_{k}+\beta_k(r(s_k,a_k,\hat{\mu}_k)-\hat{J}_{k}).\label{eq:J_online}
\end{align}
Combining Eqs.~\eqref{eq:mu_online}, \eqref{eq:V_online}, and \eqref{eq:J_online} with the $\theta$ update in \eqref{eq:wishful_update_eq1} results in a single-loop single-sample-path algorithm where in the slowest time scale we ascend the policy parameter $\theta_k$ along the gradient direction and the faster time scales are used to compute the quantities necessary for the gradient evaluation. 
While such an algorithm can be shown to converge to an MFE, the convergence does not take the best possible rate due to the coupling between iterates: $\theta_k$, $\hat{\mu}_k$, $\hat{V}_k$, and $\hat{J}_k$ directly affect each other's update, causing noise in any variable to immediately propagate to the others.  \citet{zeng2024fast} details the degradation in the algorithm complexity resulting from such coupling effect when \textit{two} variables are simultaneously updated. They further introduce a way of recovering the optimal complexity by modifying the algorithm with a denoising step. We adopt this technique and extend it to handle the three-time-scale coupling ($\alpha_k$, $\beta_k$, $\xi_k$) in our updates, which is more challenging to tackle than two timescales. The modification to the algorithm is simple -- we first estimate smoothed and denoised versions of the gradients before using them to update the policy, mean field, and value function iterates. 
We present the full details in Algorithm~\ref{alg:main}, in which the smoothed gradient estimates are $f_k$, $g_k^V$, $g_k^J$, and $h_k$ updated recursively according to line 8.

In \eqref{alg:main:meanfield}, $\Pi_{\Delta_{\Scal}}:\mathbb{R}^{|\Scal|}\rightarrow\Delta_{\Scal}$ denotes the projection to the simplex over $\Scal$.
In \eqref{alg:main:critic}, $\Pi_{B_V}:\mathbb{R}^{|\Scal|}\rightarrow\mathbb{R}^{|\Scal|}$ denotes the projection to the $\ell_2$-norm ball with radius $B_V$, and $\Pi_{[0,1]}:\mathbb{R}\rightarrow\mathbb{R}$ is the projection of a scalar to the range $[0,1]$. The projection operators guarantee the stability of the critic iterates in \eqref{alg:main:critic} and are a frequently used tool in the analysis of actor-critic algorithms in the literature \citep{wu2020finite,chen2024finite,panda2024critic}.

\section{ASSUMPTIONS \& FINITE-TIME ANALYSIS}\label{sec:main}
This section introduces the main technical assumptions made in this paper and presents the finite-time convergence of Algorithm~\ref{alg:main} to an MFE.

\begin{assump}[Uniform Geometric Ergodicity]
\label{assump:ergodic}
For any $\pi,\mu$, 
the Markov chain $\{s_k\}$ generated by $P^{\pi,\,\mu}$ according to $s_{k+1}\sim P^{\pi,\,\mu}(\cdot\mid s_k)$ is irreducible and aperiodic. In addition, there exist $C_{0}\geq 1$ and $C_1\in (0,1)$ such that
\begin{equation*}
\textstyle\sup_{s}d_{\text{TV}}\big(\mathbb{P}(s_k=\cdot \mid s_0=s),\nu^{\pi,\,\mu}(\cdot)\big)\leq C_{0}C_1^k,\, \forall k\geq 0,
\end{equation*}
where $d_{TV}$ denotes the total variation (TV) distance (see definition in Eq.~\eqref{eq:TV_def} in the appendix). 
\end{assump}
Assumption~\ref{assump:ergodic} requires that the $k_{\text{th}}$ sample of the Markov chain exponentially approaches the stationary distribution as $k$ goes up. In other words, the Markov chain generated under $P^{\pi,\,\mu}$ is geometrically ergodic for any $\pi,\mu$. This assumption is common in papers that study the complexity of sample-based single-loop RL algorithms \citep{wu2020finite,zeng2022finite,chen2024finite}.

\begin{assump}[Estimability of Induced Mean Field]\label{assump:nu}
There exists a constant $\delta\in(0,1)$ such that
\[\|\nu^{\pi,\,\mu_1}-\nu^{\pi,\,\mu_2}\|\leq\delta\|\mu_1-\mu_2\|, \quad\forall\pi,\mu_1,\mu_2.\]
\end{assump}
Assumption~\ref{assump:nu} is standard (for example, see Eq.(8) in \citet{guo2019learning}, Assumption 3 in \citet{xie2021learning}, Assumption 3 in \citet{mao2022mean}, on the contractive $\Gamma_2$ operator)
and can be viewed as an estimability condition on the induced mean field, whose validity depends only on the transition kernel $\Pcal$. 
The assumption says that for any $\pi$ the stationary distribution $\nu^{\pi,\,\mu}$ is contractive in $\mu$, and implies the uniqueness of $\mu^{\star}(\pi)$. It guarantees that to estimate the induced mean field of a policy $\pi$, we can start from any initial $\mu_0$, iteratively update it according to $\mu_k=\nu^{\pi,\mu_{k-1}}$, and have $\mu_k\rightarrow\mu^{\star}(\pi)$ as the iterations proceed. 
Note that this assumption is not to be confused with the contraction condition on the mean field optimality-consistency operator, which requires the existence of a $\delta\in(0,1)$ such that
\begin{align}
\|\nu^{\pi^{\star}(\mu_1),\mu_1}-\nu^{\pi^{\star}(\mu_2),\mu_2}\|\leq\delta\|\mu_1-\mu_2\|.\label{eq:assumption_contractive_mu}
\end{align}
\eqref{eq:assumption_contractive_mu} is a much stronger assumption with validity depending on both reward and transition kernel. While \eqref{eq:assumption_contractive_mu} is popular in prior works (Table 1), we do not need the assumption.

\begin{assump}[Lipschitz Continuity and Boundedness]\label{assump:Lipschitz}
Given two distributions $d_1,d_2$ over $\Scal$, policies $\pi_1,\pi_2$, and mean fields $\mu_1,\mu_2$, we draw samples according to
$s\sim d_1, s'\sim P^{\pi_1,\,\mu_1}(\cdot\mid s)$ and $\hat{s}\sim d_2,\hat{s}'\sim P^{\pi_2,\,\mu_2}(\cdot\mid \hat{s})$.
We assume that there exists a constant $L>0$ such that
\begin{align}
&|r(s,a,\mu_1)-r(s,a,\mu_2)|\leq L|\mu_1-\mu_2\|,\label{assump:Lipschitz:eq1}\\
&d_{TV}(\mathbb{P}(s'=\cdot), \mathbb{P}(\hat{s}'=\cdot)) \leq d_{TV}(d_1, d_2)\label{assump:Lipschitz:eq2}\\
&\hspace{100pt}+ L(\|\pi_1-\pi_2\|+\|\mu_1-\mu_2\|),\notag\\
&d_{TV}(\nu^{\pi_1,\,\mu_1}, \nu^{\pi_2,\,\mu_2})\hspace{-2pt}\leq\hspace{-2pt} L(\|\pi_1-\pi_2\|\hspace{-2pt}+\hspace{-2pt}\|\mu_1-\mu_2\|),\label{assump:Lipschitz:eq3}\\
&\|\mu^{\star}(\pi_1)-\mu^{\star}(\pi_2)\|\leq L\|\pi_1-\pi_2\|.\label{assump:Lipschitz:eq4}
\end{align}
In addition, there is a constant $B_V>0$ such that
$\|V^{\pi,\,\mu}\|\leq B_V$, for all $\pi,\mu$.
\end{assump}
Eq.~\eqref{assump:Lipschitz:eq1} states that the reward function is Lipschitz in the mean field. Eq.~\eqref{assump:Lipschitz:eq2} amounts to a regularity condition on the transition probability matrix $P^{\pi,\,\mu}$ as a function of $\pi$ and $\mu$ and can be shown to hold if the transition kernel $\Pcal^{\mu}$ is Lipschitz in $\mu$ (using an argument similar to \citet{wu2020finite}[Lemma B.2]).
Eqs.~\eqref{assump:Lipschitz:eq3} and \eqref{assump:Lipschitz:eq4} impose the Lipschitz continuity of the stationary distribution and induced mean field, which also can be shown to hold under Assumption~\ref{assump:nu} if the transition kernel is Lipschitz (see \citet{zou2019finite}[Lemma 3]). In this work, we directly assume Eqs.~\eqref{assump:Lipschitz:eq2}-\eqref{assump:Lipschitz:eq4} for simplicity.
%
All conditions in Assumption~\ref{assump:Lipschitz} are common in the literature of MFGs and RL \citep{wu2020finite,yardim2023policy,anahtarci2023q}.

\begin{assump}[Fisher Non-Degenerate Policy]
Let $\mathbb{F}(\theta)$ denote the Fisher information matrix under parameter $\theta$: 
\begin{align*}
\mathbb{F}(\theta)\hspace{-2pt}=\hspace{-2pt}\mathbb{E}_{s\sim \mu^{\star}(\pi_{\theta}),a\sim\pi_{\theta}(\cdot\mid s)}[\nabla_\theta \log \pi_\theta(a \hspace{-2pt}\mid \hspace{-2pt}s)\nabla_\theta \log \pi_\theta(a \hspace{-2pt}\mid\hspace{-2pt} s)^{\top}].
\end{align*}
There is a constant $\sigma>0$ such that $\mathbb{F}(\theta)-\sigma I_{|\Scal||\Acal|\times |\Scal||\Acal|}$ is positive definite $\forall\theta$.
\label{assump:Fisher}
\end{assump}
This assumption on Fisher non-degenerate policy implies a ``gradient domination'' condition -- for any $\mu$, every stationary point of the cumulative return $J(\pi_{\theta},\mu)$ is globally optimal.
This is again a standard assumption in the existing literature on policy optimization \citep{zhang2020global,liu2020improved,fatkhullin2023stochastic,ganesh2024variance}.
It is worth noting that we do not need Assumption~\ref{assump:Fisher} to establish the main theoretical result (Theorem~\ref{thm:main}) where the policy convergence is measured by its distance to a first-order stationary point. The assumption is only used to translate a stationary point to a globally optimal policy in Corollary~\ref{cor:main} under the average-reward formulation. 

\begin{remark}
Our algorithm and analysis can be extended to the discounted-reward setting, in which case Assumption~\ref{assump:Fisher} can be removed, as the gradient domination condition can be shown to hold under sufficient exploration (see Lemma 8 of \citet{mei2020global}). 
The same sample complexity can be established for \texttt{ASAC-MFG} under the discounted-reward formulation. 
The only major difference is the necessity of sampling from the discounted visitation measure, which requires using two sample trajectory (one for visitation measure, one for mean field). 
\end{remark}

\subsection{Main Results}\label{sec:analysis}

Each variable in Algorithm~\ref{alg:main} has a target to chase. The target of $\theta_k$ is a policy parameter optimal under its induced mean field, whereas $\hat{\mu}_k$ and $\hat{V}_k,\hat{J}_k$ aim to converge to the mean field induced by $\pi_{\theta_k}$ and the value functions under $\pi_{\theta_k},\hat{\mu}_k$. We quantify the gap between these variables and their targets by the convergence metrics below, and will shortly show that they all decay to zero.
\begin{align*}
\begin{gathered}
\varepsilon_k^{\pi} \triangleq \|\nabla_{\theta} J(\pi_{\theta_k},\mu)\hspace{-2pt}\mid_{\mu=\mu^{\star}(\pi_{\theta_k})}\hspace{-3pt}\|^2,\,
\varepsilon_k^{\mu} \triangleq \|\hat{\mu}_k\hspace{-2pt}-\hspace{-2pt}\mu^{\star}(\pi_{\theta_k})\|^2,\\
\varepsilon_k^{V}\triangleq \|\Pi_{\Ecal_\perp}(\hat{V}_k-V^{\pi_{\theta_k},\hat{\mu}_k})\|^2,\,
\varepsilon_k^{J} \triangleq (\hat{J}_k-J(\pi_{\theta_k},\hat{\mu}_k))^2.
\end{gathered}
\end{align*}
We would like $\hat{V}_k$ to converge to $V^{\pi_{\theta_k},\hat{\mu}_k}$ which solves the Bellman equation \eqref{eq:Bellman}. However, the solution is not unique. If $V\in\mathbb{R}^{|\Scal|}$ solves \eqref{eq:Bellman}, so does $V+c\1_{|\Scal|}$ for any scalar $c$. 
We denote by $\Ecal$ the subspace spanned by $\1_{|\Scal|}$ in $\mathbb{R}^{|\Scal|}$ and by $\Ecal_\perp$ its orthogonal complement, i.e., for any $V\in\Ecal_\perp$ we have $V^{\top}\1_{|\Scal|}=0$.
To make the convergence of the value function well-defined, we consider the metric $\varepsilon_k^{V}$ above where $\Pi_{\Ecal_\perp}$ is the orthogonal projection to $\Ecal_\perp$. 
It is easy to see
$\Pi_{\Ecal_\perp} = I_{|\Scal|\times|\Scal|}-\1_{|\Scal|}\1_{|\Scal|}^{\top}/|\Scal|$.

\begin{thm}\label{thm:main}
Consider the iterates generated by Algorithm~\ref{alg:main} on a $\kappa$-herding MFG, with the step sizes satisfying
\begin{align*}
\lambda_{k} \hspace{-2pt}=\hspace{-2pt} \frac{\lambda_0}{\sqrt{k+1}},\, \alpha_{k} \hspace{-2pt}=\hspace{-2pt} \frac{\alpha_0}{\sqrt{k+1}},\, \beta_{k} \hspace{-2pt}=\hspace{-2pt} \frac{\beta_0}{\sqrt{k+1}},\, \xi_{k} \hspace{-2pt}=\hspace{-2pt} \frac{\xi_0}{\sqrt{k+1}},
\end{align*}
where the constants $\lambda_0,\alpha_0,\beta_0,\xi_0$ are specified later in Appendix~\ref{sec:proof_thm:proof}.
Under Assumptions~\ref{assump:ergodic}-\ref{assump:nu}, we have for all $k\geq\tau_k$
\begin{align*}
\min_{t<k}\mathbb{E}[\varepsilon_t^{\pi}+\varepsilon_t^{\mu}+\varepsilon_t^{V}+\varepsilon_t^{J}]\leq \Ocal\left(\frac{\log^3(k+1)}{\sqrt{k+1}}+\kappa\right),
\end{align*}
where $\tau_k$ denotes the mixing time, which is an affine function of $\log(k+1)$ defined in Appendix~\ref{sec:mixing_time}. 
\end{thm}
Theorem~\ref{thm:main} states that all main variables of Algorithm~\ref{alg:main} converge to their learning targets with a rate of $\widetilde{\Ocal}(k^{-1/2})$ up to an error linear in $\kappa$, under a single trajectory of Markovian samples. 
Since Algorithm~\ref{alg:main} draws one sample in each iteration, this translates to a finite-sample complexity of the same order.
We defer the detailed proof of the theorem to Appendix~\ref{sec:proof_thm} but point out that the convergence rate is derived through a careful multi-time-scale analysis. The step sizes have the same dependency on $k$, but need to observe $\alpha_0\leq\xi_0\leq\beta_0\leq\lambda_0$. 


Our ultimate goal is to find an $\epsilon$-MFE in the sense of Definition~\ref{def:epsilon_MFE}. This requires us to connect the convergence of $\varepsilon_k^{\pi}$ (gradient norm convergence) to the optimality gap below
\begin{align}
\textstyle\max_{\pi}J(\pi,\mu^{\star}(\pi_{\theta_k}))-J(\pi_{\theta_k},\mu^{\star}(\pi_{\theta_k})).\label{eq:theta_gap}
\end{align}

Under Assumption~\ref{assump:Fisher} a ``gradient domination'' condition holds, which upper bounds \eqref{eq:theta_gap} by $\sqrt{\varepsilon_k^{\pi}}$. We take advantage of the condition to derive the following corollary.
\begin{cor}\label{cor:main}
Consider the iterates generated by Algorithm~\ref{alg:main} on a $\kappa$-herding MFG under the step sizes in Theorem~\ref{thm:main}. Under Assumptions~\ref{assump:ergodic}-\ref{assump:Fisher}, we have for all $k\geq\tau_k$
\begin{align*}
	&\textstyle\min_{t<k}\mathbb{E}\left[\max_{\pi}J(\pi,\mu^{\star}(\pi_{\theta_t}))-J(\pi_{\theta_t},\mu^{\star}(\pi_{\theta_t}))\right]\\
    &\hspace{100pt}\leq\widetilde{\Ocal}((k+1)^{-1/4})+\Ocal(\sqrt{\kappa}),\notag\\
	&\textstyle\min_{t<k}\mathbb{E}[\|\hat{\mu}_k-\mu^{\star}(\pi_{\theta_k})\|]\leq\widetilde{\Ocal}((k\hspace{-1pt}+\hspace{-1pt}1)^{-1/4})\hspace{-1pt}+\hspace{-1pt}\Ocal(\sqrt{\kappa}).
\end{align*}
\end{cor}
Corollary~\ref{cor:main} guarantees that within at most $\widetilde{\Ocal}(\epsilon^{-4})$ iterations Algorithm~\ref{alg:main} finds an $\epsilon$-MFE in the sense of Definition~\ref{def:epsilon_MFE} for MFGs in the fully herding class and an $(\epsilon+\Ocal(\sqrt{\kappa}))$-MFE for general herding MFGs.

\section{ACTOR-CRITIC ALGORITHM FOR MARKOV DECISION PROCESSES}\label{sec:MDP}
An average-reward MDP can be regarded as a degenerate average-reward MFG in which the mean field has no impact on the transition kernel or reward, i.e., $\Pcal^{\mu}(s'\mid s,a)=\Pcal(s'\mid s,a)$ and $r(s,a,\mu)=r(s,a)$.
Observing this connection, we recognize that Algorithm~\ref{alg:main} (with the mean field update \eqref{alg:main:meanfield} removed; details presented in Appendix~\ref{sec:appendix_MDP} and Algorithm~\ref{alg:MDP}) reduces to an online single-loop actor-critic algorithm that optimizes the following objective
\begin{align*}
\textstyle J_{\text{MDP}}(\pi)\triangleq\lim_{T\rightarrow\infty}\frac{1}{T}\mathbb{E}_{\pi,\Pcal}[\sum_{t=0}^{T-1}r(s_t, a_t) \mid s_0].
\end{align*}
There exist a series of works on this subject \citep{wu2020finite,olshevsky2023small,chen2024finite}, with the best-known complexity $\widetilde{\Ocal}(1/\sqrt{k})$ established in \citet{chen2024finite}. However, these prior works base their analyses on the unrealistic assumption that there exists a constant $\gamma\in(0,1)$ such that given a policy $\pi\in\Delta_{\Acal}^{\Scal}$, the following inequality holds for all $V\in\mathbb{R}^{|\Scal|}$
\begin{align}
V^{\top}\mathbb{E}_{s\sim \nu^{\pi},s'\sim P^{\pi}_{\cdot,s}}[e_s(e_{s'}-e_s)^{\top}]V\hspace{-2pt}\leq\hspace{-2pt}-\gamma\|V\|^2.
\label{assump:restrictive:eq1}
\end{align}
The assumption contradicts the common knowledge that the value function in average-reward MDP is non-unique as we discussed in Sec.\ref{sec:analysis}, and therefore never holds true. 
Fortunately, under Assumption~\ref{assump:ergodic}, it can be shown that the inequality \eqref{assump:restrictive:eq1} holds for all $V\in\Ecal_{\perp}$ (as opposed to $V\in\mathbb{R}^{|\Scal|}$). This result is stated in Lemma~\ref{lem:negative_drift} and the proof has been established in \citep{tsitsiklis1999average,zhang2021finite}.
This fact allows us to remove the assumption \eqref{assump:restrictive:eq1} in our analysis by treating the convergence of the value function in the space of $\Ecal_\perp$. 
Specifically, while the prior works consider the Lyapunov function
\[\Lcal_k=\mathbb{E}[\|\nabla_{\theta}J_{\text{MDP}}(\pi_{\theta_k})\|^2+\|\hat{V}_k-V_{\text{MDP}}^{\pi_{\theta_k}}\|^2],\]
with $V_{\text{MDP}}^{\pi}(s) \triangleq \mathbb{E}_{\pi_{\theta},\Pcal} [\sum_{t=0}^{\infty}(r(s_t, a_t)\hspace{-2pt}-\hspace{-2pt}J_{\text{MDP}}(\pi)) \hspace{-2pt}\mid\hspace{-2pt} s_0\hspace{-2pt}=\hspace{-2pt}s]$, we instead use $\|\Pi_{\Ecal_\perp}(\hat{V}_k-V_{\text{MDP}}^{\pi_{\theta_k}})\|^2$ to replace $\|\hat{V}_k-V_{\text{MDP}}^{\pi_{\theta_k}}\|^2$.
Note that we do not modify the algorithm to perform the projection to $\Ecal_\perp$ but enhance 
the analysis only.

\begin{cor}\label{cor:MDP}
Consider the policy $\pi_{\theta_k}$ generated by Algorithm~\ref{alg:MDP} with properly selected step sizes.
Under Assumptions~\ref{assump:ergodic} and \ref{assump:Lipschitz} (mapped to the context of single-agent MDP), we have for all $k\geq\tau_k$
\begin{align*}
\textstyle\min_{t<k}\mathbb{E}\big[\|\nabla_{\theta}J_{\text{MDP}}(\pi_{\theta_k})\|^2\big]\hspace{-1pt} \leq \hspace{-1pt} \Ocal\left(\log^3(k+1)/\sqrt{k\hspace{-1pt}+\hspace{-1pt}1}\right).
\end{align*}
\end{cor}
Corollary~\ref{cor:MDP} guarantees the best-iterate convergence of Algorithm~\ref{alg:MDP} to a stationary point of $J_{\text{MDP}}$, with a finite-time complexity of $\widetilde{\Ocal}(k^{-1/2})$. This matches the state-of-the-art bound in \citet{chen2024finite}, without making the restrictive assumption \eqref{assump:restrictive:eq1}.
More details on the problem formulation and algorithm in the context of MDP can be found in Appendix~\ref{sec:appendix_MDP}. The proof is presented in Appendix~\ref{sec:proof:cor:MDP}.


\section{NUMERICAL SIMULATIONS}\label{sec:simulations}

\begin{figure*}[!ht]
  \centering
  \includegraphics[width=.95\linewidth]{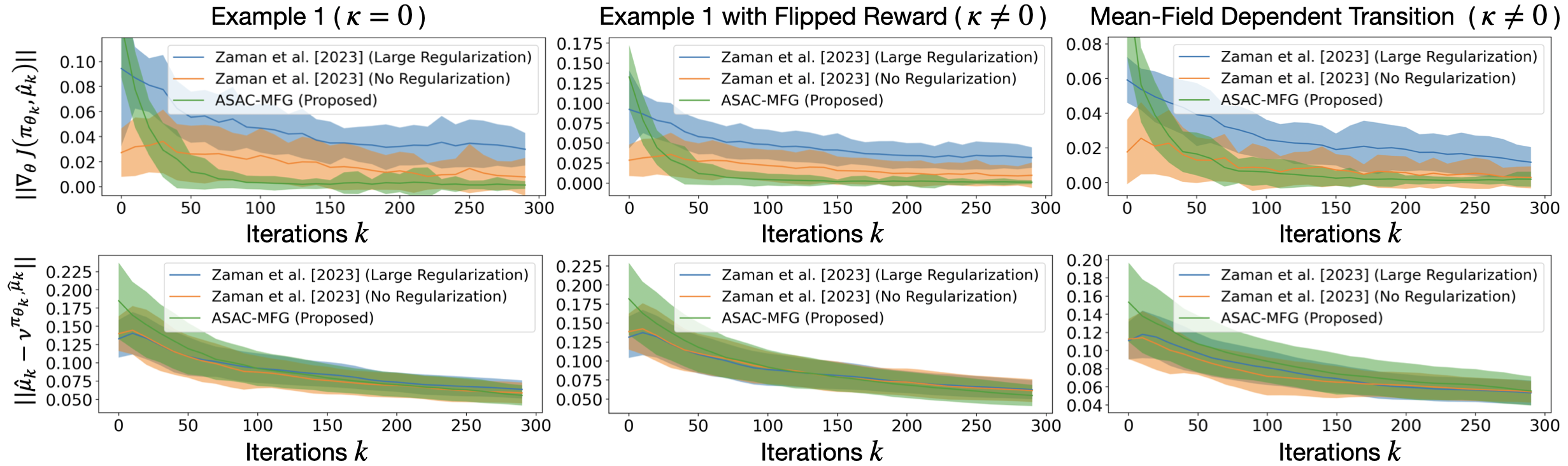}
  \vspace{-8pt}
  \caption{Algorithm performance in synthetic mean field games. 
  First row shows sub-optimality gap of policy under latest mean field estimate. Second row shows convergence of mean field estimate to mean field induced by latest policy iterate. 
  First column: Environment 1. Second column: Environment 2. Third column: Environment 3.
  }
  \vspace{-5pt}
  \label{fig:synthetic}
\end{figure*}

\begin{figure}[!ht]
  \centering
  \includegraphics[width=.8\linewidth]{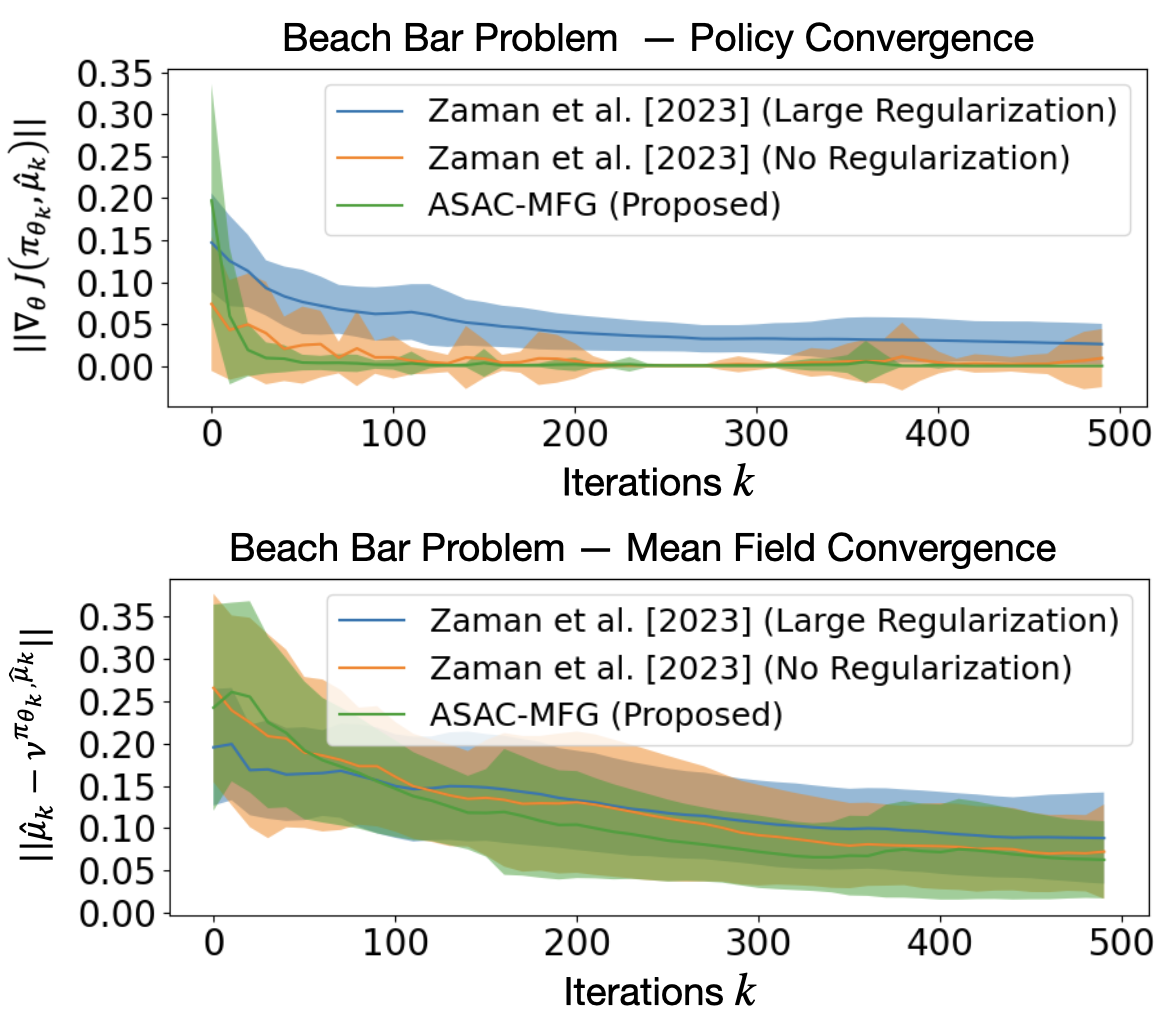}
  \vspace{-8pt}
  \caption{Algorithm performance in the beach bar problem.}
  \vspace{-10pt}
  \label{fig:beachbar}
\end{figure}

We verify the performance of the proposed algorithm through simulations on 1) small-scale synthetic MFGs, 2) the beach bar problem \citep{perrin2020fictitious}. 

First, we consider three environments of dimension $|\Scal|\hspace{-2pt}=\hspace{-2pt}|\Acal|\hspace{-2pt}=\hspace{-2pt}10$, all with randomly generated transition kernels. Environment 1 is taken from Example~\ref{example:assumption_Delta}, with the transition kernel independent of $\mu$ and the reward $r(s,a,\mu)=\mu(s)$ in expectation, which we know belongs to fully herding class and thus optimally solvable by \texttt{ASAC-MFG}. 
Environment 2 is generated the same way as Environment 1 except that the reward has a flipped sign, i.e., $r(s,a,\mu)=-\mu(s)$ in expectation, and does not satisfy \eqref{eq:def_Delta} with $\kappa=0$. Environment 3 has the same reward as Environment 1 and a random mean-field-dependent transition kernel.\footnote{More details of the experimental are in Appendix~\ref{sec:simulations_appendix}.} With high probability the environment also does not satisfy \eqref{eq:def_Delta} with $\kappa=0$.

Because the equilibria are unknown, we measure the policy convergence by $\|\nabla_{\theta}J(\pi_{\theta_k},\hat{\mu}_k)\|$ 
and the mean field convergence by $\|\hat{\mu}_k-\nu^{\pi_k,\hat{\mu}_k}\|$ as a proxy for $\|\hat{\mu}_k-\mu^{\star}(\pi_{\theta_k})\|$.

We compare \texttt{ASAC-MFG} with the algorithm proposed in~\citet{zaman2020reinforcement} as the information oracles are similar, enabling a fair comparison. We consider two variations of their algorithm: 1) with regularization large enough that the contraction assumption holds on the regularized game, and 2) with no regularization, which breaks the assumption.

As shown in Figure~\ref{fig:synthetic}, all algorithms have their mean field iterates converge to the mean field induced by the latest policy iterate, while the convergence of the policy varies. \texttt{ASAC-MFG} and \citet{zaman2023oracle} with no regularization enjoy convergence to a global MFE. However, \texttt{ASAC-MFG} converges at a faster rate, which we believe can be attributed to the efficiency of single-loop updates. 
\texttt{ASAC-MFG} is also superior in that the convergence path has a smaller variance.
The blue curve shows that while the algorithm in \citet{zaman2023oracle} with sufficiently large regularization may converge to a solution of the regularized problem, the persistent bias caused by the regularization prevents it from finding a solution to the original game. The theoretical result in Sec.\ref{sec:analysis} guarantees the convergence of \texttt{ASAC-MFG} up to an error proportional to~$\sqrt{\kappa}$, which is non-zero for Environments 2 and 3. However, we observe that in the simulations \texttt{ASAC-MFG} consistently converges to an equilibrium across all environments. This suggests that the herding condition with $\kappa=0$ may only be a sufficient condition for the solvability of an MFG, a subject worth further investigating in the future. 

We also apply the proposed algorithm to the beach bar problem, which is a common test case in the MFG literature. We take the formulation from \citet{perrin2020fictitious} with $|\Scal|=5$. The problem cannot be verified to satisfy the exact herding condition, but Figure~\ref{fig:beachbar} shows that \texttt{ASAC-MFG} still converges to an MFE, with a rate faster than that of the algorithm from \citet{zaman2023oracle}. The observation consistently matches that from Figure~\ref{fig:synthetic}.

\section{CONCLUSION}
\vspace{-5pt}

We made several important contributions to the literature on MFGs that are worth re-emphasizing: (i) We proposed a fast policy optimization algorithm for solving MFGs. 
Being the first of its kind, the algorithm is single-loop and uses a single trajectory of samples.
(ii) We identified a class of MFGs -- satisfying a novel herding condition -- that can be optimally solved by the proposed algorithm. This expands the class of solvable MFGs, including MFGs with more than one equilibrium. It is worth noting that the proposed algorithm is general purpose, in that it can be applied to MFGs not in the proposed class to obtain stable iterates that provably converge to an approximate equilibrium (Corollary~\ref{cor:main}). (iii) The current analysis for actor-critic algorithms in average-reward MDPs requires an assumption that cannot possibly hold. By recognizing that an MFG reduces to a standard MDP with transition kernel and reward independent of the mean field, we showed that our main result leads to a finite-sample analysis of an actor-critic algorithm for average-reward MDPs that matches the state-of-the-art complexity without the untenable assumption.

While the mean-field considered in this paper is a distribution over the states, we expect the algorithm and the analysis to go through with minimal modifications when the mean-field is a more general distribution over states and actions. 
An important future work is to characterize the explicit connections between the different known solvable classes of MFGs (cf.~\ref{fig:Venn}).

\section*{Disclaimer}
This paper was prepared for informational purposes by the Artificial Intelligence Research group of JPMorgan Chase \& Co. and its affiliates ("JP Morgan'') and is not a product of the Research Department of JP Morgan. JP Morgan makes no representation and warranty whatsoever and disclaims all liability, for the completeness, accuracy or reliability of the information contained herein. This document is not intended as investment research or investment advice, or a recommendation, offer or solicitation for the purchase or sale of any security, financial instrument, financial product or service, or to be used in any way for evaluating the merits of participating in any transaction, and shall not constitute a solicitation under any jurisdiction or to any person, if such solicitation under such jurisdiction or to such person would be unlawful.

\bibliographystyle{plainnat} 
\bibliography{references}

\clearpage

\section*{Checklist}



 \begin{enumerate}

 \item For all models and algorithms presented, check if you include:
 \begin{enumerate}
   \item A clear description of the mathematical setting, assumptions, algorithm, and/or model. \textcolor{red}{Yes}
   \item An analysis of the properties and complexity (time, space, sample size) of any algorithm. \textcolor{red}{Yes}
   \item (Optional) Anonymized source code, with specification of all dependencies, including external libraries. \textcolor{red}{Yes}
 \end{enumerate}

 \item For any theoretical claim, check if you include:
 \begin{enumerate}
   \item Statements of the full set of assumptions of all theoretical results. \textcolor{red}{Yes}
   \item Complete proofs of all theoretical results. \textcolor{red}{Yes}
   \item Clear explanations of any assumptions. \textcolor{red}{Yes}   
 \end{enumerate}

 \item For all figures and tables that present empirical results, check if you include:
 \begin{enumerate}
   \item The code, data, and instructions needed to reproduce the main experimental results (either in the supplemental material or as a URL). \textcolor{red}{Yes}
   \item All the training details (e.g., data splits, hyperparameters, how they were chosen). \textcolor{red}{Not Applicable. No training involved for the problem we study.} 
   \item A clear definition of the specific measure or statistics and error bars (e.g., with respect to the random seed after running experiments multiple times). \textcolor{red}{Yes. We plot both mean and standard deviation in all figures.}
   \item A description of the computing infrastructure used. (e.g., type of GPUs, internal cluster, or cloud provider). \textcolor{red}{Not Applicable. Experiments are all small-scale and no more than 30 minutes to run on a standard computer.}
 \end{enumerate}

 \item If you are using existing assets (e.g., code, data, models) or curating/releasing new assets, check if you include:
 \begin{enumerate}
   \item Citations of the creator If your work uses existing assets. \textcolor{red}{Not Applicable.}
   \item The license information of the assets, if applicable. \textcolor{red}{Not Applicable.}
   \item New assets either in the supplemental material or as a URL, if applicable. \textcolor{red}{Not Applicable.}
   \item Information about consent from data providers/curators. \textcolor{red}{Not Applicable.}
   \item Discussion of sensible content if applicable, e.g., personally identifiable information or offensive content. \textcolor{red}{Not Applicable.}
 \end{enumerate}

 \item If you used crowdsourcing or conducted research with human subjects, check if you include:
 \begin{enumerate}
   \item The full text of instructions given to participants and screenshots. \textcolor{red}{Not Applicable.}
   \item Descriptions of potential participant risks, with links to Institutional Review Board (IRB) approvals if applicable. \textcolor{red}{Not Applicable.}
   \item The estimated hourly wage paid to participants and the total amount spent on participant compensation. \textcolor{red}{Not Applicable.}
 \end{enumerate}

 \end{enumerate}

\clearpage
\appendix

\onecolumn
\tableofcontents

\section{Notations and Frequently Used Identities}
We introduce a few shorthand notations frequently used in the analysis.
First, we define
\begin{align}
\begin{aligned}
F(\theta,V,\mu,s,a,s')&\triangleq (r(s,a,\mu) + V(s')-V(s))\nabla_{\theta}\log\pi_{\theta}(a\mid s),\\
G^V(V,J,\mu,s,a,s') &\triangleq (r(s,a,\mu) -J + V(s') - V(s))e_{s},\\
G^J(J,\mu,s,a) &\triangleq c_J(r(s,a,\mu)-J),\\
G(V,J,\mu,s,a,s') &\triangleq \left[\begin{array}{c}
G^V(V,J,\mu,s,a,s') \\
G^J(J,\mu,s,a) 
\end{array}\right]=\left[\begin{array}{c}
(r(s,a,\mu) -J + V(s') - V(s))e_{s} \\
c_J (r(s,a,\mu) - J) 
\end{array}\right],\\
H(\mu,s)&\triangleq e_s-\mu.
\end{aligned}
\label{eq:def_FGH}
\end{align}
Then, the update of $f_k$, $g_k^V$, $g_k^J$, and $h_k$ in Algorithm~\ref{alg:main} can be expressed as
\begin{gather*}
f_{k+1} = (1-\lambda_k)f_{k}+\lambda_k F(\theta_k,\hat{V}_k,\hat{\mu}_k,s_k,a_k,s_{k+1}),\\
g_{k+1}^V = (1-\lambda_k)g_{k}^V+\lambda_k G^V(\hat{V}_k,\hat{J}_k,\hat{\mu}_k,s_k,a_k,s_{k+1}),\\
g_{k+1}^J = (1-\lambda_k)g_{k}^J+\lambda_k G^J(\hat{J}_k,\hat{\mu}_k,s_k,a_k),\\
h_{k+1} = (1-\lambda_k)h_{k}+\lambda_k H(\hat{\mu}_k,s_k).
\end{gather*}
Denote $g_k=[(g_k^V)^{\top},g_k^J]^{\top}$. The update of $g_k$ is
\[g_{k+1} = \left[\begin{array}{c}
g_{k+1}^V \\
g_{k+1}^J
\end{array}\right] = (1-\lambda_k)g_{k} + \lambda_k G(\hat{V}_k,\hat{J}_k,\hat{\mu}_k,s_k,a_k,s_{k+1}).\]

We also define 
\begin{align}
\begin{aligned}
\bar{F}(\theta,V,\mu)&\triangleq\mathbb{E}_{s\sim \nu^{\pi_{\theta},\,\mu},a\sim\pi_{\theta}(\cdot\mid s),s'\sim\Pcal^{\mu}(\cdot\mid s,a)}[F(\theta,V,\mu,s,a,s')],\\
\bar{G}^V(\theta,V,J,\mu)&\triangleq\mathbb{E}_{s\sim \nu^{\pi_{\theta},\,\mu},a\sim\pi_{\theta}(\cdot\mid s),s'\sim\Pcal^{\mu}(\cdot\mid s,a)}[G^V(V,J,\mu,s,a,s')],\\
\bar{G}^J(\theta,J,\mu)&\triangleq\mathbb{E}_{s\sim \nu^{\pi_{\theta},\,\mu},a\sim\pi_{\theta}(\cdot\mid s)}[G^J(J,\mu,s,a)],\\
\bar{G}(\theta,V,J,\mu)&\triangleq\mathbb{E}_{s\sim \nu^{\pi_{\theta},\,\mu},a\sim\pi_{\theta}(\cdot\mid s),s'\sim\Pcal^{\mu}(\cdot\mid s,a)}[G(V,J,\mu,s,a,s')]=\left[\begin{array}{c}
\bar{G}^V(\theta,V,J,\mu) \\
\bar{G}^J(\theta,J,\mu) 
\end{array}\right],\\
\bar{H}(\theta,\mu)&\triangleq\mathbb{E}_{s\sim \nu^{\pi_{\theta},\,\mu}}[H(\mu,s)]=\mathbb{E}_{s\sim \nu^{\pi_{\theta},\,\mu}}[e_s-\mu].
\end{aligned}
\label{eq:def_FGH_aggregate}
\end{align}

We measure the convergence of auxiliary variables $f_k$, $g_k^V$, $g_k^J$, and $h_k$ by
\begin{gather*}
\Delta f_k \triangleq f_k-\bar{F}(\theta_k,\hat{V}_k,\hat{\mu}_k), \quad \Delta g_k^V \triangleq g_k^V-\bar{G}^V(\theta_k,\hat{V}_k,\hat{J}_k,\hat{\mu}_k), \\
\Delta g_k^J \triangleq g_k^J-\bar{G}^J(\theta_k,\hat{J}_k,\hat{\mu}_k),\quad \Delta h_k \triangleq h_k-\bar{H}(\theta_k,\hat{\mu}_k),
\end{gather*}
and denote 
\[\Delta g_k = \left[\begin{array}{c}
\Delta g_k^V \\
\Delta g_k^J
\end{array}\right] = g_k-\bar{G}(\theta_k,\hat{V}_k,\hat{J}_k,\hat{\mu}_k).\]

We use $\ell(\pi)$ to denote the cumulative reward collected by policy $\pi$ under the induced mean field $\mu^{\star}(\pi)$
\[\ell(\pi)\triangleq J(\pi,\mu^{\star}(\pi)).\]
This is well-defined since $\mu^{\star}(\pi)$ is unique.

We denote by $\Fcal_k=\{s_0,a_0,s_1,a_1\cdots,s_k,a_k\}$ denote the filtration (set of all randomness information) up to iteration $k$. Given two probability distributions $\phi_1$ and $\phi_2$ over space $\Xcal$, their TV distance is defined as
\begin{equation}
d_{\text{TV}}(\phi_1,\phi_2)\hspace{-2pt}=\hspace{-2pt}\frac{1}{2} \textstyle\sup _{\psi: \Xcal \rightarrow[-1,1]}\left|\int \psi d \phi_1\hspace{-2pt}-\hspace{-2pt}\int \psi d \phi_2\right|.
\label{eq:TV_def}
\end{equation}

Under Assumptions~\ref{assump:ergodic} and \ref{assump:Lipschitz}, it can be shown using an argument similar to Lemma B.1 of \citet{wu2020finite} that there exists a constant $L_{TV}$ depending only on $|\Acal|$, $L$, $C_0$, and $C_1$ such that for all $\pi_1,\pi_2,\mu_1,\mu_2$
\begin{align}
d_{TV}(\nu^{\pi_1,\,\mu_1}\otimes\pi_1\otimes\Pcal^{\mu_1},\nu^{\pi_2,\,\mu_2}\otimes\pi_2\otimes\Pcal^{\mu_2})\leq L_{TV}(\|\pi_1-\pi_2\|+\|\mu_1-\mu_2\|).\label{eq:TV_stationary_Lipschitz}
\end{align}

Without loss of generality, we assume $L\geq 1$, a condition that we will sometimes use to simplify and combine terms.

\subsection{Mixing Time}\label{sec:mixing_time}

An immediate consequence of Assumption~\ref{assump:ergodic} is that the Markov chain under any policy and mean field has a geometric mixing time.
%
\begin{definition}\label{def:mixing_time}
Consider a Markov chain $\{\hat{s}_k\}$ generated according to $\hat{s}_k\sim P^{\pi,\,\mu}(\cdot\mid \hat{s}_{k-1})$, for which $\nu^{\pi,\,\mu}$ is the stationary distribution. For any $c>0$, the $c$-mixing time of the Markov chain is
\[
    \tau^{\pi,\,\mu}(c) \triangleq \min\big\{k\in\mathbb{N}:\textstyle\sup_{s}d_{TV}\big(\mathbb{P}(\hat{s}_k=\cdot\mid \hat{s}_0=s),\,\nu^{\pi,\,\mu}(\cdot)\big)\leq c\big\}.
\]
\end{definition}
The mixing time measures time for the samples of the Markov chain to approach its stationary distribution in TV distance. We define $\tau_k\triangleq\sup_{\pi,\,\mu}\tau^{\pi,\,\mu}(\alpha_k)$ as the time when the TV distance drops below $\alpha_k$, where $\alpha_k$ is a step size for the policy parameter update in Algorithm~\ref{alg:main}.
Under Assumption~\ref{assump:ergodic}, it is obvious that there exists a constant $C$ as a function of $C_0, C_1$ such that 
\[\tau_k \leq C\log\left(1/\alpha_k\right)=C\log(\frac{(k+1)^{1/2}}{\alpha_0})=\frac{C}{2}\log(k+1)-C\log(\alpha_0).\]

\subsection{Supporting Lemmas}

The value function $V^{\pi_{\theta},\,\mu}$ is Lipschitz in both $\theta$ and $\mu$, as shown in the lemma below.
\begin{lem}\label{lem:Lipschitz_V}
Under Assumption~\ref{assump:Lipschitz}, there exist a bounded constant $L_V\geq 1$ such that for any policy parameter $\theta_1,\theta_2$ and mean field $\mu_1,\mu_2$, we have
\begin{gather*}
\|\Pi_{\Ecal_\perp}(V^{\pi_{\theta_1},\,\mu_1}-V^{\pi_{\theta_2},\,\mu_2})\| \leq L_V\left(\|\theta_1-\theta_2\|+\|\mu_1-\mu_2\|\right),\\
\|J(\pi_{\theta_1},\,\mu_1)-J(\pi_{\theta_2},\,\mu_2)\| \leq L_V\left(\|\theta_1-\theta_2\|+\|\mu_1-\mu_2\|\right),\\
\|\nabla_{\theta}J(\pi_{\theta_1},\,\mu_1)-\nabla_{\theta}J(\pi_{\theta_2},\,\mu_2)\| \leq L_V\left(\|\theta_1-\theta_2\|+\|\mu_1-\mu_2\|\right),\\
\|\nabla_{\mu}J(\pi_{\theta_1},\,\mu_1)-\nabla_{\mu}J(\pi_{\theta_2},\,\mu_2)\| \leq L_V\left(\|\theta_1-\theta_2\|+\|\mu_1-\mu_2\|\right).
\end{gather*}
\end{lem}

We establish the boundedness of the operators $F$, $G$, and $H$. 
\begin{lem}\label{lem:boundedness}
For any $\theta\in\mathbb{R}^{|\Scal||\Acal|}$, $V\in\mathbb{R}^{|\Scal|}$ with norm bounded by $B_V$, $J\in[0,1]$, $\mu\in\Delta_{\Scal}$, and $s,a,s'$, we have
\begin{gather*}
\|F(\theta,V,\mu,s,a,s')\|\leq B_F,
\|G(V,J,\mu,s,a,s')\|\leq B_G,
\|H(\mu,s)\|\leq B_H,
\end{gather*}
where $B_F=B_V+1$, $B_G=2(B_V+c_J+2)$, $B_H=2$.
\end{lem}

Since $f_k$, $g_k^V$, $g_k^J$, and $h_k$ are simply convex combination with the operators $F$, $G^V$, $G^J$, and $H$, Lemma~\ref{lem:boundedness} implies for all $k$
\begin{align*}
\|f_k\|\leq B_F, \quad \|g_k^V\|\leq B_G, \quad|g_k^J|\leq B_G, \quad\|h_k\|\leq B_H.
\end{align*}

We also establish the Lipschitz continuity of a few important operators. 
\begin{lem}\label{lem:Lipschitz_operators}
We have for any $\theta_1,\theta_2\in\mathbb{R}^{|\Scal||\Acal|}$, $\mu_1,\mu_2\in\Delta_{\Scal}$, $V_1,V_2\in\mathbb{R}^{|\Scal|}$, and $J_1,J_2\in\mathbb{R}$
\begin{align*}
&\|\bar{F}(\theta_1,V_1,\mu_1)-\bar{F}(\theta_2,V_2,\mu_2)\|\leq L_F\Big(\|\theta_1-\theta_2\|+\|\Pi_{\Ecal_\perp}(V_1-V_2)\|+\|\mu_1-\mu_2\|\Big)\\
&\|\bar{G}(\theta_1,V_1,J_1,\mu_1)-\bar{G}(\theta_2,V_2,J_2,\mu_2)\|\notag\\
&\hspace{100pt}\leq L_G\Big(\|\theta_1-\theta_2\|+\|\Pi_{\Ecal_\perp}(V_1-V_2)\|+|J_1-J_2|+\|\mu_1-\mu_2\|\Big),\\
&\|\bar{H}(\theta_1,\mu_1)-\bar{H}(\theta_2,\mu_2)\|\leq L_H\Big(\|\theta_1-\theta_2\|+\|\mu_1-\mu_2\|\Big),
\end{align*}
where the constants are $L_F=10B_V+L+2B_F L_{TV}+5$, $L_G=2B_G L_{TV}+(L+1)(c_J+1)+2$, and $L_H=L+1$.
\end{lem}

As a result of Lemma~\ref{lem:Lipschitz_operators}, we can establish the following bounds on the energy of the auxiliary variables $f_k$, $g_k$, and $h_k$.
\begin{lem}\label{lem:bound_f_g}
We have for any $k\geq 0$
\begin{align*}
\|f_k\|&\leq\|\Delta f_k\|+L_F\sqrt{\varepsilon_k^{V}}+L_F(L_V+1)\sqrt{\varepsilon_k^{\mu}}+\sqrt{\varepsilon_k^{\pi}},\\
\|g_k\|&\leq\|\Delta g_k\|+L_G\sqrt{\varepsilon_k^V}+L_G\sqrt{\varepsilon_k^J},\\
\|h_k\|&\leq\|\Delta h_k\|+L_H\sqrt{\epsilon_k^{\mu}}.
\end{align*}
\end{lem}

Also as a consequence of Assumption~\ref{assump:ergodic}, the following lemma holds which states that the Bellman backup operator of the value function is almost everywhere contractive (except along the direction of the all-one vector). This lemma is adapted from \citet{zhang2021finite}[Lemma 2] and \citet{tsitsiklis1999average}[Lemma 7].
\begin{lem}\label{lem:negative_drift}
Recall the definition of $\Ecal_\perp$ in Sec.\ref{sec:analysis}.
There exists a constant $\gamma\in(0,1)$ such that for any $\theta,\mu$ and $V\in\Ecal_\perp$
\begin{align*}
V^{\top}\mathbb{E}_{s\sim \nu^{\pi_{\theta},\,\mu},a\sim\pi_{\theta}(\cdot\mid s),s'\sim\Pcal^{\mu}(\cdot\mid s,a)}[e_s(e_{s'}-e_s)^{\top}]V\leq-\gamma\|V\|^2.
\end{align*}
\end{lem}

\section{Proof of Main Theorem}\label{sec:proof_thm}

\subsection{Intermediate Results}

The proof of Theorem~\ref{thm:main} relies critically on the iteration-wise convergence of policy iterate $\theta_k$, mean field iterate $\hat{\mu}_k$, value function estimate $\hat{V}_k$, $\hat{J}_k$, and auxiliary variables $f_k$, $h_k$, and $g_k$, which we bound individually in the propositions below.

\subsubsection{Convergence of Policy Iterate}

\begin{prop}\label{prop:policy_conv}
Under Assumptions~\ref{assump:ergodic}-\ref{assump:Lipschitz}, we have
\begin{align*}
\ell(\pi_{\theta_{k}})-\ell(\pi_{\theta_{k+1}}) 
&\leq -\frac{\rho \alpha_k}{2}\varepsilon_k^{\pi}+\rho \alpha_k\|\Delta f_k\|^2\notag\\
&\hspace{20pt}+\rho L_F^2\alpha_k(\varepsilon_k^{V}+\varepsilon_k^{\mu})+\frac{\rho L_V B_F^2\alpha_k^2}{2}+B_F\alpha_k\kappa.
\end{align*}
\end{prop}

\begin{prop}\label{prop:f_conv}
Under Assumptions~\ref{assump:ergodic}-\ref{assump:Lipschitz}, we have for all $k\geq \tau_k$
\begin{align*}
&\mathbb{E}[\|\Delta f_{k+1}\|^2]\notag\\
&\leq (1-\lambda_k)\mathbb{E}[\|\Delta f_{k}\|^2]+(-\frac{\lambda_k}{2}+\lambda_k^2+\frac{48L_F^2 \alpha_k^2}{\lambda_k})\mathbb{E}[\|\Delta f_{k}\|^2]\notag\\
&\hspace{20pt}+\frac{36L_F^2\beta_k^2}{\lambda_k}\mathbb{E}[\|\Delta g_k\|^2]+\frac{24L_F^2 L_H^2 \xi_k^2}{\lambda_k}\mathbb{E}[\|\Delta h_k\|^2]+\frac{48L_F^2 \alpha_k^2}{\lambda_k}\mathbb{E}[\varepsilon_k^{\pi}]+\frac{216L_F^4 L_V^2 \xi_k^2}{\lambda_k}\mathbb{E}[\varepsilon_k^{\mu}]\notag\\
&\hspace{20pt}+\frac{96L_F^4 L_G^2\beta_k^2}{\lambda_k}\mathbb{E}[\varepsilon_k^{V}]+\frac{48L_F^2 L_G^2 \beta_k^2}{\lambda_k}\mathbb{E}[\varepsilon_k^{J}]+(28L+2|\Acal|) L_F L_{TV} B_F^3 B_G B_H^2\tau_k^2\lambda_k\lambda_{k-\tau_k}.
\end{align*}
\end{prop}

The proofs of Propositions~\ref{prop:policy_conv} and \ref{prop:f_conv} can be found in Sec.\ref{sec:policy_conv:proof} and \ref{sec:f_conv:proof}.

\subsubsection{Convergence of Mean Field Estimate}

\begin{prop}\label{prop:meanfield_conv}
Under Assumptions~\ref{assump:ergodic}-\ref{assump:nu}, we have for all $k$
\begin{align*}
\varepsilon_{k+1}^{\mu} \leq (1-\frac{(1-\delta)\xi_k}{8})\varepsilon_k^{\mu}+\frac{8\xi_k}{1-\delta}\|\Delta h_k\|^2+\frac{32L^2\alpha_k^2}{(1-\delta)\xi_k}\left(\|\Delta f_k\|^2+L_F^2\varepsilon_k^{V}+\varepsilon_k^{\pi}\right)+9L^2 B_F^2 B_H^2\xi_k^2.
\end{align*}
\end{prop}

\begin{prop}\label{prop:h_conv}
Under Assumptions~\ref{assump:ergodic}-\ref{assump:Lipschitz}, we have for all $k\geq\tau_k$
\begin{align*}
&\mathbb{E}[\|\Delta h_{k+1}\|^2]\notag\\
&\leq(1-\lambda_k)\mathbb{E}[\|\Delta h_k\|^2]+(-\frac{\lambda_k}{2}+\lambda_k^2+\frac{16L_H^2\xi_k^2}{\lambda_k})\mathbb{E}[\|\Delta h_k\|^2]+\frac{32L_H^2 \alpha_k^2}{\lambda_k}\mathbb{E}[\|\Delta f_k\|^2]\notag\\
&\hspace{20pt}+\frac{32L_H^2 L_F^2\alpha_k^2}{\lambda_k}\mathbb{E}[\varepsilon_k^{V}]+\frac{144 L_F^2 L_V^2 L_H^4 \xi_k^2}{\lambda_k}\mathbb{E}[\varepsilon_k^{\mu}]+\frac{32L_H^2 \alpha_k^2}{\lambda_k}\mathbb{E}[\varepsilon_k^{\pi}] + 24L B_F B_H^2\tau_k^2\lambda_k\lambda_{k-\tau_k}.
\end{align*}
\end{prop}

The proofs of Propositions~\ref{prop:meanfield_conv} and \ref{prop:h_conv} can be found in Sec.\ref{sec:meanfield_conv:proof} and \ref{sec:h_conv:proof}.

\subsubsection{Convergence of Valuation Function Estimate}

\begin{prop}\label{prop:V_conv}
Under Assumptions~\ref{assump:ergodic}-\ref{assump:Lipschitz},
\begin{align*}
\varepsilon_{k+1}^{V}+\varepsilon_{k+1}^{J}&\leq (1-\frac{\gamma \beta_k}{4}) (\varepsilon_k^{V}+\varepsilon_k^{J})+\frac{128L_V^2\alpha_k^2}{\gamma\beta_k} \|\Delta f_k\|^2+\frac{8\beta_k}{\gamma}\|\Delta g_k\|^2+\frac{64L_V^2\xi_k^2}{\gamma\beta_k} \|\Delta h_k\|^2\notag\\
&\hspace{20pt}+\frac{128L_V^2\alpha_k^2}{\gamma\beta_k} (L_F^2 \varepsilon_k^{V}+\varepsilon_k^{\pi})+\frac{192L_V^2\xi_k^2}{\gamma\beta_k} \varepsilon_k^{\mu}+28L_V^2 B_F^2 B_G^2 B_H^2 \beta_k^2.
\end{align*}
\end{prop}

\begin{prop}\label{prop:g_conv}
Under Assumptions~\ref{assump:ergodic}-\ref{assump:Lipschitz}, we have for all $k\geq \tau_k$
\begin{align*}
\mathbb{E}[\|\Delta g_{k+1}\|^2]&\leq(1-\lambda_k)\mathbb{E}[\|\Delta g_{k}\|^2]+(-\frac{\lambda_k}{2}+\lambda_k^2+\frac{72|\Scal|L_G^2 \beta_k^2}{\lambda_k})\mathbb{E}[\|\Delta g_{k}\|^2]+\frac{48L_G^2 \alpha_k^2}{\lambda_k}\mathbb{E}[\|\Delta f_k\|^2]\notag\\
&\hspace{20pt}+\frac{24L_G^2 \xi_k^2}{\lambda_k}\mathbb{E}[\|\Delta h_k\|^2]+\frac{48L_G^2 \alpha_k^2}{\lambda_k}\mathbb{E}[\varepsilon_{k}^{\pi}]+\frac{216L_F^2L_G^2L_H^2L_V^2\xi_k^2}{\lambda_k}\mathbb{E}[\varepsilon_{k}^{\mu}]\notag\\
&\hspace{20pt}+\frac{120|\Scal|L_F^2 L_G^4\beta_k^2}{\lambda_k}\mathbb{E}[\varepsilon_k^{V}+\varepsilon_k^{J}]+(30L+2|\Acal|) L_F L_{TV} B_F B_G^2 B_H \tau_k^2\lambda_k\lambda_{k-\tau_k}.
\end{align*}
\end{prop}

The proofs of Propositions~\ref{prop:V_conv} and \ref{prop:g_conv} can be found in Sec.\ref{sec:V_conv:proof} and \ref{sec:g_conv:proof}.

\subsection{Proof of Theorem~\ref{thm:main}}\label{sec:proof_thm:proof}

The exact requirements on $\lambda_0,\alpha_0,\beta_0,\xi_0$ include $c_J\geq1/\gamma$, $\alpha_0\leq \xi_0 \leq \beta_0 \leq \lambda_0$, and
\begin{align}\label{eq:step_size_details}
\begin{gathered}
\alpha_0\leq\min\left\{\frac{1}{192(L_F^2+L_G^2+L_H^2+L_V^2+L^2/(1-\delta)+\rho)}\lambda_0,\,C_{\beta}\beta_0,\,C_{\xi}\xi_0\right\},\\
\xi_0\leq\min\left\{\frac{\lambda_0}{64(L_H^2 L_F^2+L_G^2+L_V^2/\gamma+1/(1-\delta))},\,\frac{(1-\delta)\gamma\beta_0}{6912(L_F^4 L_V^2+L_F^2 L_G^2 L_H^2 L_V^2+L_F^2 L_H^4 L_V^2+L_V^2)}\right\},\\
\beta_0\leq\min\left\{\frac{\lambda_0}{72|\Scal|L_G^2+36L_F^2+8/\gamma},\,\frac{\gamma}{4L_G^2},\frac{1-\delta}{2L_H^2}\right\},\quad\lambda_0\leq\frac{1}{4},
\end{gathered}
\end{align}
where $C_{\xi}=\min\{\frac{(1-\delta)}{32 \rho L_F^2},\frac{1-\delta}{4\rho },\frac{L_H}{2L_F L_V},\frac{1-\delta}{16L L_F L_V}\}$ and 
\begin{align*}
C_{\beta}&=\min\Big\{\frac{\gamma}{4},\frac{\rho \gamma}{512(L_F^2+L_G^2+L_H^2+L_V^2+L^2/(1-\delta)},\\
&\hspace{20pt}\sqrt{\frac{\gamma}{3456|\Scal|(L_F^4 L_G^4+L_F^2L_H^2+\rho L_F^2+L_V^2/\gamma+L^2 L_F^2(1-\delta))}},\frac{\gamma}{2\rho }\Big\}.
\end{align*}
We note that such parameters can always chosen with no conflict in any MFG.

We consider the potential function
\begin{align*}
\Lcal_k = \mathbb{E}[\|\Delta f_k\|^2+\|\Delta g_k\|^2+\|\Delta h_k\|^2-\ell(\pi_{\theta_k})+\varepsilon_k^V+\varepsilon_k^J+\varepsilon_k^{\mu}].
\end{align*}

Collecting the bounds from Propositions~\ref{prop:policy_conv}-\ref{prop:g_conv}, we have for all $k\geq\tau_k$
\begin{align}
&\Lcal_{k+1} \notag\\
&= \mathbb{E}[\|\Delta f_{k+1}\|^2+\|\Delta g_{k+1}\|^2+\|\Delta h_{k+1}\|^2-\ell(\pi_{\theta_{k+1}})+\varepsilon_{k+1}^V+\varepsilon_{k+1}^J+\varepsilon_{k+1}^{\mu}]\notag\\
&\leq (1-\lambda_k)\mathbb{E}[\|\Delta f_{k}\|^2]+(-\frac{\lambda_k}{2}+\lambda_k^2+\frac{48L_F^2 \alpha_k^2}{\lambda_k})\mathbb{E}[\|\Delta f_{k}\|^2]\notag\\
&\hspace{20pt}+\frac{36L_F^2\beta_k^2}{\lambda_k}\mathbb{E}[\|\Delta g_k\|^2]+\frac{24L_F^2 L_H^2 \xi_k^2}{\lambda_k}\mathbb{E}[\|\Delta h_k\|^2]+\frac{48L_F^2 \alpha_k^2}{\lambda_k}\mathbb{E}[\varepsilon_k^{\pi}]+\frac{216L_F^4 L_V^2 \xi_k^2}{\lambda_k}\mathbb{E}[\varepsilon_k^{\mu}]\notag\\
&\hspace{20pt}+\frac{96L_F^4 L_G^2\beta_k^2}{\lambda_k}\mathbb{E}[\varepsilon_k^{V}]+\frac{48L_F^2 L_G^2 \beta_k^2}{\lambda_k}\mathbb{E}[\varepsilon_k^{J}]+(28L+2|\Acal|) L_F L_{TV} B_F^3 B_G B_H^2\tau_k^2\lambda_k\lambda_{k-\tau_k}\notag\\
&\hspace{20pt}+(1-\lambda_k)\mathbb{E}[\|\Delta g_{k}\|^2]+(-\frac{\lambda_k}{2}+\lambda_k^2+\frac{72|\Scal|L_G^2 \beta_k^2}{\lambda_k})\mathbb{E}[\|\Delta g_{k}\|^2]+\frac{48L_G^2 \alpha_k^2}{\lambda_k}\mathbb{E}[\|\Delta f_k\|^2]\notag\\
&\hspace{20pt}+\frac{24L_G^2 \xi_k^2}{\lambda_k}\mathbb{E}[\|\Delta h_k\|^2]+\frac{48L_G^2 \alpha_k^2}{\lambda_k}\mathbb{E}[\varepsilon_{k}^{\pi}]+\frac{216L_F^2L_G^2L_H^2L_V^2\xi_k^2}{\lambda_k}\mathbb{E}[\varepsilon_{k}^{\mu}]\notag\\
&\hspace{20pt}+\frac{120|\Scal|L_F^2 L_G^4\beta_k^2}{\lambda_k}\mathbb{E}[\varepsilon_k^{V}+\varepsilon_k^{J}]+(30L+2|\Acal|) L_F L_{TV} B_F B_G^2 B_H \tau_k^2\lambda_k\lambda_{k-\tau_k}\notag\\
&\hspace{20pt}+(1-\lambda_k)\mathbb{E}[\|\Delta h_k\|^2]+(-\frac{\lambda_k}{2}+\lambda_k^2+\frac{16L_H^2\xi_k^2}{\lambda_k})\mathbb{E}[\|\Delta h_k\|^2]+\frac{32L_H^2 \alpha_k^2}{\lambda_k}\mathbb{E}[\|\Delta f_k\|^2]\notag\\
&\hspace{20pt}+\hspace{-2pt}\frac{32L_H^2 L_F^2\alpha_k^2}{\lambda_k}\mathbb{E}[\varepsilon_k^{V}]\hspace{-2pt}+\hspace{-2pt}\frac{144 L_F^2 L_V^2 L_H^4 \xi_k^2}{\lambda_k}\mathbb{E}[\varepsilon_k^{\mu}]\hspace{-2pt}+\hspace{-2pt}\frac{32L_H^2 \alpha_k^2}{\lambda_k}\mathbb{E}[\varepsilon_k^{\pi}] + 24L B_F B_H^2\tau_k^2\lambda_k\lambda_{k-\tau_k}\notag\\
&\hspace{20pt}-\mathbb{E}[\ell(\pi_{\theta_k})]-\frac{\rho \alpha_k}{2}\mathbb{E}[\varepsilon_k^{\pi}]+\rho \alpha_k\mathbb{E}[\|\Delta f_k\|^2]\notag\\
&\hspace{20pt}+\rho L_F^2\alpha_k\mathbb{E}[\varepsilon_k^{V}+\varepsilon_k^{\mu}]+\frac{\rho L_V B_F^2\alpha_k^2}{2}+B_F\alpha_k\kappa\notag\\
&\hspace{20pt}+(1-\frac{\gamma \beta_k}{4}) \mathbb{E}[\varepsilon_k^{V}+\varepsilon_k^{J}]+\frac{128L_V^2\alpha_k^2}{\gamma\beta_k} \mathbb{E}[\|\Delta f_k\|^2]+\frac{8\beta_k}{\gamma}\mathbb{E}[\|\Delta g_k\|^2]+\frac{64L_V^2\xi_k^2}{\gamma\beta_k} \mathbb{E}[\|\Delta h_k\|^2]\notag\\
&\hspace{20pt}+\frac{128L_V^2\alpha_k^2}{\gamma\beta_k} (L_F^2 \mathbb{E}[\varepsilon_k^{V}]+\mathbb{E}[\varepsilon_k^{\pi}])+\frac{192L_V^2\xi_k^2}{\gamma\beta_k} \mathbb{E}[\varepsilon_k^{\mu}]+28L_V^2 B_F^2 B_G^2 B_H^2 \beta_k^2\notag\\
&\hspace{20pt}+\hspace{-2pt}(1\hspace{-2pt}-\hspace{-2pt}\frac{(1\hspace{-2pt}-\hspace{-2pt}\delta)\xi_k}{8})\mathbb{E}[\varepsilon_k^{\mu}]\hspace{-2pt}+\hspace{-2pt}\frac{8\xi_k}{1\hspace{-2pt}-\hspace{-2pt}\delta}\mathbb{E}[\|\Delta h_k\|^2]\hspace{-2pt}+\hspace{-2pt}\frac{32L^2\alpha_k^2}{(1\hspace{-2pt}-\hspace{-2pt}\delta)\xi_k}\mathbb{E}[\|\Delta f_k\|^2\hspace{-2pt}+\hspace{-2pt}L_F^2\varepsilon_k^{V}\hspace{-2pt}+\hspace{-2pt}\varepsilon_k^{\pi}]\hspace{-2pt}+\hspace{-2pt}9L^2 B_F^2 B_H^2\xi_k^2\notag\\
&\leq (1-\lambda_k)\mathbb{E}[\|\Delta f_{k}\|^2+\|\Delta g_{k}\|^2+\|\Delta h_k\|^2]-\mathbb{E}[\ell(\pi_{\theta_k})]-\frac{\rho \alpha_k}{4}\mathbb{E}[\varepsilon_k^{\pi}]\notag\\
&\hspace{20pt}+(1-\frac{\gamma \beta_k}{8}) \mathbb{E}[\varepsilon_k^{V}+\varepsilon_k^{J}]+(1-\frac{(1-\delta)\xi_k}{16})\mathbb{E}[\varepsilon_k^{\mu}]+B_F\alpha_k\kappa\notag\\
&\hspace{20pt}+(28L+2|\Acal|) L_F L_{TV} B_F^3 B_G B_H^2\tau_k^2\lambda_k\lambda_{k-\tau_k}+(30L+2|\Acal|) L_F L_{TV} B_F B_G^2 B_H \tau_k^2\lambda_k\lambda_{k-\tau_k}\notag\\
&\hspace{20pt}+ 24LB_F B_H^2\tau_k^2\lambda_k\lambda_{k-\tau_k}+\frac{\rho L_V B_F^2\alpha_k^2}{2}+28L_V^2 B_F^2 B_G^2 B_H^2 \beta_k^2+9L^2 B_F^2 B_H^2\xi_k^2\notag\\
&\hspace{20pt}+\underbrace{(-\frac{\lambda_k}{2}+\lambda_k^2+\frac{48L_F^2 \alpha_k^2}{\lambda_k}+\frac{48L_G^2\alpha_k^2}{\lambda_k}+\frac{32L_H^2\alpha_k^2}{\lambda_k}+\rho \alpha_k+\frac{128L_V^2\alpha_k^2}{\gamma\beta_k}+\frac{32L^2\alpha_k^2}{(1-\delta)\lambda_k})}_{A_1}\mathbb{E}[\|\Delta f_{k}\|^2]\notag\\
&\hspace{20pt}+\underbrace{(-\frac{\lambda_k}{2}+\lambda_k^2+\frac{72|\Scal|L_G^2 \beta_k^2}{\lambda_k}+\frac{36L_F^2 \beta_k^2}{\lambda_k}+\frac{8\beta_k}{\gamma})}_{A_2}\mathbb{E}[\|\Delta g_{k}\|^2]\notag\\
&\hspace{20pt}+\underbrace{(-\frac{\lambda_k}{2}+\lambda_k^2+\frac{16L_H^2\xi_k^2}{\lambda_k}+\frac{24L_F^2L_H^2\xi_k^2}{\lambda_k}+\frac{24L_G^2\xi_k^2}{\lambda_k}+\frac{64L_V^2\xi_k^2}{\gamma\lambda_k}+\frac{8\xi_k}{1-\delta})}_{A_3}\mathbb{E}[\|\Delta h_k\|^2]\notag\\
&\hspace{20pt}+\underbrace{(-\frac{\rho \alpha_k}{4}+\frac{48L_F^2\alpha_k^2}{\lambda_k}+\frac{48L_G^2\alpha_k^2}{\lambda_k}+\frac{32L_H^2\alpha_k^2}{\lambda_k}+\frac{128L_V^2\alpha_k^2}{\gamma\beta_k}+\frac{32L^2\alpha_k^2}{(1-\delta)\lambda_k})}_{A_4}\mathbb{E}[\varepsilon_k^{\pi}]\notag\\
&\hspace{20pt}+\underbrace{(-\frac{\gamma \beta_k}{8}\hspace{-2pt}+\hspace{-2pt}\frac{96L_F^4L_G^2\beta_k^2}{\lambda_k}+\frac{120|\Scal|L_F^2 L_G^4\beta_k^2}{\lambda_k}\hspace{-2pt}+\hspace{-2pt}\frac{32L_F^2 L_H^2\alpha_k^2}{\lambda_k}+\rho L_F^2\alpha_k\hspace{-2pt}+\hspace{-2pt}\frac{128L_V^2\alpha_k^2}{\gamma\beta_k}\hspace{-2pt}+\hspace{-2pt}\frac{32L^2 L_F^2 \alpha_k^2}{(1-\delta)\lambda_k})}_{A_5} \mathbb{E}[\varepsilon_k^{V}\hspace{-2pt}+\hspace{-2pt}\varepsilon_k^{J}]\notag\\
&\hspace{20pt}+\underbrace{(-\frac{(1-\delta)\xi_k}{16}+\frac{216L_F^4 L_V^2\xi_k^2}{\lambda_k}+\frac{216L_F^2L_G^2L_H^2 L_V^2\xi_k^2}{\lambda_k}+\frac{144L_F^2 L_H^4 L_V^2\xi_k^2}{\lambda_k}+\rho L_F^2\alpha_k+\frac{192L_V^2\xi_k^2}{\gamma\beta_k})}_{A_6}\mathbb{E}[\varepsilon_k^{\mu}].\label{thm:main:proof_eq1}
\end{align}

We show that the terms $A_1$-$A_6$ are all non-positive under the step size conditions in \eqref{eq:step_size_details}. First, under the step size condition $\alpha_k\leq\frac{\gamma}{4}\beta_k$, $\lambda_k\leq 1/4$, and $\alpha_k\leq(192(L_F^2+L_G^2+L_H^2+L_V^2+L^2/(1-\delta)+\rho))^{-1}\lambda_k$
\begin{align}
A_1&=-\frac{\lambda_k}{2}+\lambda_k^2+\frac{48L_F^2 \alpha_k^2}{\lambda_k}+\frac{48L_G^2\alpha_k^2}{\lambda_k}+\frac{32L_H^2\alpha_k^2}{\lambda_k}+\rho \alpha_k+\frac{128L_V^2\alpha_k^2}{\gamma\beta_k}+\frac{32L^2\alpha_k^2}{(1-\delta)\lambda_k}\notag\\
&\leq -\frac{\lambda_k}{4}+\frac{48(L_F^2+L_G^2+L_H^2+L^2/(1-\delta)) \alpha_k^2}{\lambda_k}+\rho \alpha_k+32L_V^2\alpha_k\notag\\
&\leq -\frac{\lambda_k}{4}+48(L_F^2+L_G^2+L_H^2+L_V^2+L^2/(1-\delta)+\rho) \alpha_k\notag\\
&\leq 0.\label{thm:main:proof_eq2}
\end{align}

Next, under the step size condition $\lambda_k\leq 1/4$ and $\beta_k\leq(72|\Scal|L_G^2+36L_F^2+8/\gamma)^{-1}\lambda_k$
\begin{align}
A_2 &= -\frac{\lambda_k}{2}+\lambda_k^2+\frac{72|\Scal|L_G^2 \beta_k^2}{\lambda_k}+\frac{36L_F^2 \beta_k^2}{\lambda_k}+\frac{8\beta_k}{\gamma}\notag\\
&\leq-\frac{\lambda_k}{4}+(72|\Scal|L_G^2+36L_F^2+8/\gamma)\beta_k\notag\\
&\leq 0.
\end{align}

Next, under the step size condition $\lambda_k\leq 1/4$ and $\xi_k\leq(64(L_H^2 L_F^2+L_G^2+L_V^2/\gamma+1/(1-\delta)))^{-1}\lambda_k$
\begin{align}
A_3 &= -\frac{\lambda_k}{2}+\lambda_k^2+\frac{16L_H^2\xi_k^2}{\lambda_k}+\frac{24L_F^2L_H^2\xi_k^2}{\lambda_k}+\frac{24L_G^2\xi_k^2}{\lambda_k}+\frac{64L_V^2\xi_k^2}{\gamma\lambda_k}+\frac{8\xi_k}{1-\delta}\notag\\
&\leq -\frac{\lambda_k}{4}+64(L_H^2 L_F^2+L_G^2+L_V^2/\gamma+1/(1-\delta))\xi_k\notag\\
&\leq 0.
\end{align}

Next, we have
\begin{align}
A_4 &= -\frac{\rho \alpha_k}{4}+\frac{48L_F^2\alpha_k^2}{\lambda_k}+\frac{48L_G^2\alpha_k^2}{\lambda_k}+\frac{32L_H^2\alpha_k^2}{\lambda_k}+\frac{128L_V^2\alpha_k^2}{\gamma\beta_k}+\frac{32L^2\alpha_k^2}{(1-\delta)\lambda_k}\notag\\
&\leq-\frac{\rho \alpha_k}{4}+\frac{128}{\gamma}(L_F^2+L_G^2+L_H^2+L_V^2+L^2/(1-\delta))\frac{\alpha_k^2}{\beta_k}\notag\\
&\leq 0,
\end{align}
under the step size condition
\[\alpha_k\leq\frac{\rho \gamma}{512(L_F^2+L_G^2+L_H^2+L_V^2+L^2/(1-\delta)}\beta_k.\]

Then,
\begin{align}
A_5 &= -\frac{\gamma \beta_k}{8}+\frac{96L_F^4L_G^2\beta_k^2}{\lambda_k}+\frac{120|\Scal|L_F^2 L_G^4\beta_k^2}{\lambda_k}+\frac{32L_F^2 L_H^2\alpha_k^2}{\lambda_k}\notag\\
&\hspace{20pt}+\rho L_F^2\alpha_k+\frac{128L_V^2\alpha_k^2}{\gamma\beta_k}+\frac{32L^2 L_F^2 \alpha_k^2}{(1-\delta)\lambda_k}\notag\\
&\leq -\frac{\gamma \beta_k}{8}+432|\Scal|(L_F^4 L_G^4+L_F^2L_H^2+\rho L_F^2+L_V^2/\gamma+L^2 L_F^2(1-\delta))\frac{\alpha_k^2}{\beta_k}\notag\\
&\leq 0,
\end{align}
due to the condition
\[\alpha_k\leq\sqrt{\frac{\gamma}{3456|\Scal|(L_F^4 L_G^4+L_F^2L_H^2+\rho L_F^2+L_V^2/\gamma+L^2 L_F^2(1-\delta))}}\beta_k.\]

Finally, as a result of $\alpha_k\leq \frac{(1-\delta)}{32\rho L_F^2}\xi_k$ and $\xi_k\leq\frac{(1-\delta)\gamma}{6912(L_F^4 L_V^2+L_F^2 L_G^2 L_H^2 L_V^2+L_F^2 L_H^4 L_V^2+L_V^2)}\beta_k$
\begin{align}
A_6 &= -\frac{(1-\delta)\xi_k}{16}+\frac{216L_F^4 L_V^2\xi_k^2}{\lambda_k}+\frac{216L_F^2L_G^2L_H^2 L_V^2\xi_k^2}{\lambda_k}\notag\\
&\hspace{20pt}+\frac{144L_F^2 L_H^4 L_V^2\xi_k^2}{\lambda_k}+\rho L_F^2\alpha_k+\frac{192L_V^2\xi_k^2}{\gamma\beta_k}\notag\\
&\leq-\frac{(1-\delta)\xi_k}{32}+\frac{216}{\gamma}(L_F^4 L_V^2+L_F^2 L_G^2 L_H^2 L_V^2+L_F^2 L_H^4 L_V^2+L_V^2)\frac{\xi_k^2}{\beta_k}\notag\\
&\leq 0.\label{thm:main:proof_eq3}
\end{align}

Plugging \eqref{thm:main:proof_eq2}-\eqref{thm:main:proof_eq3} into \eqref{thm:main:proof_eq1}, we have for all $k\geq\tau_k$
\begin{align}
&\Lcal_{k+1} \notag\\
&\leq (1-\lambda_k)\mathbb{E}[\|\Delta f_{k}\|^2+\|\Delta g_{k}\|^2+\|\Delta h_k\|^2]-\mathbb{E}[\ell(\pi_{\theta_k})]-\frac{\rho \alpha_k}{4}\mathbb{E}[\varepsilon_k^{\pi}]\notag\\
&\hspace{20pt}+(1-\frac{\gamma \beta_k}{8}) \mathbb{E}[\varepsilon_k^{V}+\varepsilon_k^{J}]+(1-\frac{(1-\delta)\xi_k}{16})\mathbb{E}[\varepsilon_k^{\mu}]+B_F\alpha_k\kappa\notag\\
&\hspace{20pt}+(28L\hspace{-2pt}+\hspace{-2pt}2|\Acal|) L_F L_{TV} B_F^3 B_G B_H^2\tau_k^2\lambda_k\lambda_{k-\tau_k}+(30L\hspace{-2pt}+\hspace{-2pt}2|\Acal|) L_F L_{TV} B_F B_G^2 B_H \tau_k^2\lambda_k\lambda_{k-\tau_k}\notag\\
&\hspace{20pt}+ 24LB_F B_H^2\tau_k^2\lambda_k\lambda_{k-\tau_k}+\frac{\rho L_V B_F^2\alpha_k^2}{2}+28L_V^2 B_F^2 B_G^2 B_H^2 \beta_k^2+9L^2 B_F^2 B_H^2\xi_k^2\notag\\
&\leq \Lcal_k-\min\left\{\frac{\rho \alpha_k}{4},\frac{\gamma\beta_k}{8},\frac{(1\hspace{-2pt}-\hspace{-2pt}\delta)\xi_k}{16}\right\}\mathbb{E}[\varepsilon_k^{\pi}+\varepsilon_k^{\mu}+\varepsilon_k^{V}+\varepsilon_k^{J}]+B_F\alpha_k\kappa+\Ocal(\frac{\log^2(k+1)}{k+1})\notag\\
&\leq \Lcal_k-\frac{\rho \alpha_k}{4}\mathbb{E}[\varepsilon_k^{\pi}+\varepsilon_k^{\mu}+\varepsilon_k^{V}+\varepsilon_k^{J}]+B_F\alpha_k\kappa+\Ocal(\frac{\log^2(k+1)}{k+1}),
\end{align}
where the last inequality follows from the step size condition $\alpha_k\leq\frac{\gamma}{2\rho }\beta_k$ and $\alpha_k\leq\frac{1-\delta}{4\rho }\xi_k$.

Re-arranging the terms and summing over iterations, we have
\begin{align*}
\sum_{t=\tau_k}^{k-1}\alpha_t\mathbb{E}[\varepsilon_t^{\pi}+\varepsilon_t^{\mu}+\varepsilon_t^{V}+\varepsilon_t^{J}] &\leq \frac{4}{\rho}\sum_{t=\tau_k}^{k-1}(\Lcal_t-\Lcal_{t+1})+B_F\kappa\sum_{t=\tau_k}^{k-1}\alpha_t+\sum_{t=\tau_k}^{k-1}\Ocal(\frac{\log^2(t+1)}{t+1})\notag\\
&\leq \frac{4}{\rho}(\Lcal_{\tau_k}+1)+B_F\kappa\sum_{t=\tau_k}^{k-1}\alpha_t+\Ocal(\log^3(k+1)),
\end{align*}
where the second inequality follows from $-\Lcal_{k+1}\leq-\ell(\pi_{\theta_{k+1}})\leq1$ and the well-known relationship that
\[\sum_{t=\tau_k}^{k-1}\frac{1}{t+1}\leq\sum_{t=0}^{k-1}\frac{1}{t+1}\leq2\log(k+1).\]

Due to $\tau_k\leq\Ocal(\log(k+1))$, it is also a standard result that (for example, see \citet{zeng2021two}[Lemma 3])
\[\sum_{t=\tau_k}^{k-1}\alpha_t=\sum_{t=\tau_k}^{k-1}\frac{\alpha_0}{\sqrt{t+1}}=\Theta(k+1).\]

Dividing both sides of the inequality by $\sum_{t=\tau_k}^{k-1}\alpha_t$, we get
\begin{align*}
\min_{t<k}\mathbb{E}[\varepsilon_t^{\pi}+\varepsilon_t^{\mu}+\varepsilon_t^{V}+\varepsilon_t^{J}]&\leq\frac{\sum_{t=\tau_k}^{k-1}\alpha_t\mathbb{E}[\varepsilon_t^{\pi}+\varepsilon_t^{\mu}+\varepsilon_t^{V}+\varepsilon_t^{J}]}{\sum_{t=\tau_k}^{k-1}\alpha_t}\notag\\
&\leq \Ocal(\frac{1}{\sqrt{k+1}})\left(\frac{4}{\rho}(\Lcal_{\tau_k}+1)+\Ocal(\log^3(k+1))\right)+B_F\kappa.
\end{align*}

Since the updates of all iterates in Algorithm~\ref{alg:main} are bounded, $\Lcal_{\tau_k}\leq\Ocal(\tau_k)\leq\Ocal(\log(k+1))$. As a result, we eventually have
\begin{align*}
\min_{\tau_k\leq t<k}\mathbb{E}[\varepsilon_t^{\pi}+\varepsilon_t^{\mu}+\varepsilon_t^{V}+\varepsilon_t^{J}]\leq \Ocal\left(\frac{\log^3(k+1)}{\sqrt{k+1}}\right)+\Ocal(\kappa).
\end{align*}

\qed

\section{Proof of Corollaries}

\subsection{Proof of Corollary~\ref{cor:main}}

As a result of Assumption~\ref{assump:Fisher}, we have the following gradient domination condition, which is adapted from Lemma 19 of \citet{ganesh2024variance}.

\begin{lem}\label{lem:gradient_domination}
Under Assumption~\ref{assump:Fisher}, we have the following gradient domination condition for any policy parameter $\theta$ and mean field $\mu$
\begin{align*}
\max_{\bar{\pi}}J(\bar{\pi},\mu)-J(\pi_{\theta},\mu)\leq\frac{1}{\sigma}\|\nabla_{\theta} J(\pi_{\theta},\mu)\|.
\end{align*}
\end{lem}

Since $\varepsilon_t^{\pi},\varepsilon_t^{\mu},\varepsilon_t^{V},\varepsilon_t^{J}$ are all non-negative, we have
\begin{gather*}
\min_{\tau_k\leq t<k}\mathbb{E}\left[\|\nabla_{\theta}J(\pi_{\theta_t},\mu)\mid_{\mu=\mu^{\star}(\pi_{\theta_t})}\|^2\right]\hspace{-2pt}\leq\hspace{-2pt}\Ocal\left(\frac{\log^3(k+1)}{\sqrt{k+1}}\right)\hspace{-2pt}+\hspace{-2pt}\Ocal(\kappa)\hspace{-2pt}=\hspace{-2pt}\widetilde{\Ocal}\left(\frac{\log^3(k+1)}{\sqrt{k+1}}\right)\hspace{-2pt}+\hspace{-2pt}\Ocal(\kappa),\notag\\
\min_{\tau_k\leq t<k}\mathbb{E}[\|\hat{\mu}_k-\mu^{\star}(\pi_{\theta_k})\|^2]\leq\Ocal\left(\frac{\log^3(k+1)}{\sqrt{k+1}}\right)+\Ocal(\kappa)=\widetilde{\Ocal}\left(\frac{\log^3(k+1)}{\sqrt{k+1}}\right)+\Ocal(\kappa).
\end{gather*}

Applying Lemma~\ref{lem:gradient_domination} with $\theta=\theta_t$ and $\mu=\mu^{\star}(\pi_{\theta_t})$, 
\begin{align*}
\max_{\pi}J(\pi,\mu^{\star}(\pi_{\theta_t}))-J(\pi_{\theta},\mu^{\star}(\pi_{\theta_t}))\leq\frac{1}{\sigma}\|\nabla_{\theta} J(\pi_{\theta_t},\mu)\mid_{\mu=\mu^{\star}(\pi_{\theta_t})}\|.
\end{align*}

By Jensen's inequality,
\begin{align*}
&\left(\min_{\tau_k\leq t<k}\mathbb{E}\left[\max_{\pi}J(\pi,\mu^{\star}(\pi_{\theta_t}))-J(\pi_{\theta_t},\mu^{\star}(\pi_{\theta_t}))\right]\right)^2\notag\\
&\leq\min_{\tau_k\leq t<k}\mathbb{E}\left[\left(\max_{\pi}J(\pi,\mu^{\star}(\pi_{\theta_t}))-J(\pi_{\theta_t},\mu^{\star}(\pi_{\theta_t}))\right)^2\right]\notag\\
&\leq \frac{1}{\sigma^2}\min_{\tau_k\leq t<k}\mathbb{E}\left[\|\nabla_{\theta}J(\pi_{\theta_t},\mu)\mid_{\mu=\mu^{\star}(\pi_{\theta_t})}\|^2\right]\notag\\
&\leq\widetilde{\Ocal}\left(\frac{1}{\sqrt{k+1}}\right)+\Ocal(\kappa).
\end{align*}
Taking square root on both sides of this inequality leads to the claimed result on the convergence of the policy.

Similarly, we have
\begin{align*}
\min_{\tau_k\leq t<k}\mathbb{E}[\|\hat{\mu}_k-\mu^{\star}(\pi_{\theta_k})\|]&\leq\sqrt{\min_{\tau_k\leq t<k}\mathbb{E}[\|\hat{\mu}_k-\mu^{\star}(\pi_{\theta_k})\|^2]}\notag\\
&\leq\sqrt{\widetilde{\Ocal}\left(\frac{\log^3(k+1)}{\sqrt{k+1}}\right)+\Ocal(\kappa)}\notag\\
&\leq \widetilde{\Ocal}\left(\frac{1}{(k+1)^{1/4}}\right)+\Ocal(\sqrt{\kappa}).
\end{align*}

\qed

\subsection{Proof of Corollary~\ref{cor:MDP}}\label{sec:proof:cor:MDP}

In the context of single-agent MDP we define
\begin{align*}
\varepsilon_k^{\pi}=\|\nabla_{\theta}J_{\text{MDP}}(\pi_{\theta_k})\|^2,\quad\varepsilon_k^{V}=\|\Pi_{\Ecal_\perp}(\hat{V}_k-V_{MDP}^{\pi_{\theta_k}})\|^2,\quad\varepsilon_k^{J}=(\hat{J}_k-J_{\text{MDP}}(\pi_{\theta_k}))^2.
\end{align*}

Theorem~\ref{thm:main} implies that
\begin{align*}
    \min_{\tau_k\leq t<k}\mathbb{E}[\varepsilon_t^{\pi}+\varepsilon_t^{V}+\varepsilon_t^{J}]\leq \Ocal\left(\frac{\log^3(k+1)}{\sqrt{k+1}}\right)+\Ocal(\kappa).
\end{align*}

As neither the transition kernel nor the reward function depends on the mean field, \eqref{eq:def_Delta} trivially holds with $\rho=1,\kappa=0$. Therefore, the claimed result holds by recognizing that $\varepsilon_k^{V}$ and $\varepsilon_k^J$ are non-negative for any $k$.

\qed
\section{Proof of Propositions}

\subsection{Proof of Proposition~\ref{prop:policy_conv}}\label{sec:policy_conv:proof}

By the $L_V$-Lipschitz continuity of the function $J$
\begin{align}
&J(\pi_{\theta_{k}},\mu^{\star}(\pi_{\theta_{k}}))-J(\pi_{\theta_{k+1}},\mu^{\star}(\pi_{\theta_{k}})) \notag\\
&\leq -\langle \nabla_{\theta}J(\pi_{\theta_k},\mu)\mid_{\mu=\mu^{\star}(\pi_{\theta_{k}})},\theta_{k+1}-\theta_k\rangle+\frac{L_V}{2}\|\theta_{k+1} - \theta_k\|^2\notag\\
&= -\alpha_k\langle \nabla_{\theta}J(\pi_{\theta_k},\mu)\mid_{\mu=\mu^{\star}(\pi_{\theta_{k}})},f_k \rangle+\frac{L_V \alpha_k^2}{2}\|f_k\|^2\notag\\
&= \hspace{-2pt}-\alpha_k\langle \nabla_{\theta}J(\pi_{\theta_k},\mu)\hspace{-2pt}\mid_{\mu=\mu^{\star}(\pi_{\theta_{k}})},\Delta f_k\rangle\hspace{-2pt}-\hspace{-2pt}\alpha_k\langle \nabla_{\theta}J(\pi_{\theta_k},\mu)\hspace{-2pt}\mid_{\mu=\mu^{\star}(\pi_{\theta_{k}})},\hspace{-1pt}\bar{F}(\theta_k,\hat{V}_k,\hat{\mu}_k)\rangle\hspace{-2pt}+\hspace{-2pt}\frac{L_V\alpha_k^2}{2}\|f_k\|^2\notag\\
&= -\alpha_k\langle \nabla_{\theta}J(\pi_{\theta_k},\mu)\mid_{\mu=\mu^{\star}(\pi_{\theta_{k}})}, \Delta f_k\rangle-\alpha_k\| \nabla_{\theta}J(\pi_{\theta_k},\mu)\mid_{\mu=\mu^{\star}(\pi_{\theta_k})}\|^2\notag\\
&\hspace{20pt}+\alpha_k\langle \nabla_{\theta}J(\pi_{\theta_k},\mu)\hspace{-2pt}\mid_{\mu=\mu^{\star}(\pi_{\theta_{k}})},\bar{F}(\theta_k,V^{\pi_{\theta_k},\,\mu^{\star}(\pi_{\theta_k})},\mu^{\star}(\pi_{\theta_k}))\hspace{-2pt}-\hspace{-2pt}\bar{F}(\theta_k,\hat{V}_k,\hat{\mu}_k)\rangle\hspace{-2pt}+\hspace{-2pt}\frac{L_V\alpha_k^2}{2}\|f_k\|^2\notag\\
&\leq -\alpha_k\| \nabla_{\theta}J(\pi_{\theta_k},\mu)\mid_{\mu=\mu^{\star}(\pi_{\theta_k})} \|^2 -\alpha_k\langle \nabla_{\theta}J(\pi_{\theta_k},\mu)\mid_{\mu=\mu^{\star}(\pi_{\theta_{k}})},\Delta f_k\rangle\notag\\
&\hspace{20pt}+\alpha_k\langle \nabla_{\theta}J(\pi_{\theta_k},\mu)\mid_{\mu=\mu^{\star}(\pi_{\theta_{k}})},\bar{F}(\theta_k,V^{\pi_{\theta_k},\,\mu^{\star}(\pi_{\theta_k})},\mu^{\star}(\pi_{\theta_k}))-\bar{F}(\theta_k,\hat{V}_k,\hat{\mu}_k)\rangle\hspace{-2pt}+\hspace{-2pt}\frac{L_V B_F^2\alpha_k^2}{2},\label{prop:policy_conv:proof_eq1}
\end{align}
where the third equation follows from $\nabla_{\theta}J(\pi_{\theta},\mu)\mid_{\mu=\mu^{\star}(\pi_{\theta})}\,=\bar{F}(\theta,V^{\pi_{\theta},\,\mu^{\star}(\pi_{\theta})},\mu^{\star}(\pi_{\theta}))$ for any $\theta$.

To bound the second term on the right hand side of \eqref{prop:policy_conv:proof_eq1}, we use the fact that $\langle\vec{a},\vec{b}\rangle\leq \frac{c}{2}\|\vec{a}\|^2+\frac{1}{2c}\|\vec{b}\|^2$ for any vectors $\vec{a},\vec{b}$ and scalar $c>0$
\begin{align}
-\alpha_k\langle \nabla_{\theta}J(\pi_{\theta_k},\mu)\mid_{\mu=\mu^{\star}(\pi_{\theta_{k}})},\Delta f_k\rangle&\leq \frac{\alpha_k}{4}\|\nabla_{\theta}J(\pi_{\theta_k},\mu)\mid_{\mu=\mu^{\star}(\pi_{\theta_{k}})}\|^2+\alpha_k\|\Delta f_k\|^2.
\label{prop:policy_conv:proof_eq2}
\end{align}

Similarly, for the third term of \eqref{prop:policy_conv:proof_eq1}, we have
\begin{align}
&\alpha_k\langle \nabla_{\theta}J(\pi_{\theta_k},\mu)\mid_{\mu=\mu^{\star}(\pi_{\theta_{k}})},\mu^{\star}(\pi_{\theta_k}))-\bar{F}(\theta_k,\hat{V}_k,\hat{\mu}_k)\rangle\notag\\
&\leq\frac{\alpha_k}{4}\|\nabla_{\theta}J(\pi_{\theta_k},\mu)\mid_{\mu=\mu^{\star}(\pi_{\theta_{k}})}\|^2+\alpha_k\|\bar{F}(\theta_k,V^{\pi_{\theta_k},\,\mu^{\star}(\pi_{\theta_k})},\mu^{\star}(\pi_{\theta_k}))-\bar{F}(\theta_k,\hat{V}_k,\hat{\mu}_k)\|^2\notag\\
&\leq \frac{\alpha_k}{4}\|\nabla_{\theta}J(\pi_{\theta_k},\mu)\mid_{\mu=\mu^{\star}(\pi_{\theta_{k}})}\|^2\hspace{-2pt}+\hspace{-2pt}L_F^2\alpha_k\|\Pi_{\Ecal_\perp}(\hat{V}_k-V^{\pi_{\theta_k},\mu^{\star}(\pi_{\theta_k})})\|^2+L_F^2\alpha_k\|\hat{\mu}_k-\mu^{\star}(\pi_{\theta_k})\|^2\notag\\
&= \frac{\alpha_k}{4}\|\nabla_{\theta}J(\pi_{\theta_k},\mu)\mid_{\mu=\mu^{\star}(\pi_{\theta_k})}\|^2+L_F^2\alpha_k(\varepsilon_k^{V}+\varepsilon_k^{\mu}).\label{prop:policy_conv:proof_eq3}
\end{align}

Plugging \eqref{prop:policy_conv:proof_eq2}-\eqref{prop:policy_conv:proof_eq3} into \eqref{prop:policy_conv:proof_eq1}, we have
\begin{align}
&J(\pi_{\theta_{k}},\mu^{\star}(\pi_{\theta_{k}}))-J(\pi_{\theta_{k+1}},\mu^{\star}(\pi_{\theta_{k}}))\notag\\
&\leq -\alpha_k\| \nabla_{\theta}J(\pi_{\theta_k},\mu)\mid_{\mu=\mu^{\star}(\pi_{\theta_k})} \|^2 -\alpha_k\langle \nabla_{\theta}J(\pi_{\theta_k},\mu)\mid_{\mu=\mu^{\star}(\pi_{\theta_{k}})},\Delta f_k\rangle\notag\\
&\hspace{20pt}+\alpha_k\langle \nabla_{\theta}J(\pi_{\theta_k},\mu)\mid_{\mu=\mu^{\star}(\pi_{\theta_{k}})},\bar{F}(\theta_k,V^{\pi_{\theta_k},\,\mu^{\star}(\pi_{\theta_k})},\mu^{\star}(\pi_{\theta_k}))-\bar{F}(\theta_k,\hat{V}_k,\hat{\mu}_k)\rangle+\frac{L_V B_F^2\alpha_k^2}{2}\notag\\
&\leq -\alpha_k\| \nabla_{\theta}J(\pi_{\theta_k},\mu)\mid_{\mu=\mu^{\star}(\pi_{\theta_k})}\|^2+\frac{\alpha_k}{4}\|\nabla_{\theta}J(\pi_{\theta_k},\mu)\mid_{\mu=\mu^{\star}(\pi_{\theta_{k}})}\|^2+\alpha_k\|\Delta f_k\|^2\notag\\
&\hspace{20pt}+\frac{\alpha_k}{4}\|\nabla_{\theta}J(\pi_{\theta_k},\mu)\mid_{\mu=\mu^{\star}(\pi_{\theta_k})}\|^2+L_F^2\alpha_k(\varepsilon_k^{V}+\varepsilon_k^{\mu})+\frac{L_V B_F^2\alpha_k^2}{2}\notag\\
&\leq -\frac{\alpha_k}{2}\|\nabla_{\theta}J(\pi_{\theta_k},\mu)\mid_{\mu=\mu^{\star}(\pi_{\theta_k})}\|^2+\alpha_k\|\Delta f_k\|^2+L_F^2\alpha_k(\varepsilon_k^{V}+\varepsilon_k^{\mu})+\frac{L_V B_F^2\alpha_k^2}{2}.
\label{prop:policy_conv:proof_eq4}
\end{align}

By \eqref{eq:def_Delta}, we have
\begin{align*}
&J(\pi_{\theta_{k}},\mu^{\star}(\pi_{\theta_{k}}))-J(\pi_{\theta_{k+1}},\mu^{\star}(\pi_{\theta_{k+1}}))\notag\\
&\leq 
\rho\Big(-\frac{\alpha_k}{2}\|\nabla_{\theta}J(\pi_{\theta_k},\mu)\mid_{\mu=\mu^{\star}(\pi_{\theta_k})}\|^2+\alpha_k\|\Delta f_k\|^2+L_F^2\alpha_k(\varepsilon_k^{V}+\varepsilon_k^{\mu})+\frac{L_V B_F^2\alpha_k^2}{2}\Big)+B_F\alpha_k\kappa\notag\\
&\leq -\frac{\rho\alpha_k}{2}\|\nabla_{\theta}J(\pi_{\theta_k},\mu)\mid_{\mu=\mu^{\star}(\pi_{\theta_k})}\|^2+\rho\alpha_k\|\Delta f_k\|^2+\rho L_F^2\alpha_k(\varepsilon_k^{V}+\varepsilon_k^{\mu})+\frac{\rho L_V B_F^2\alpha_k^2}{2}+B_F\alpha_k\kappa.
\end{align*}

\qed

\subsection{Proof of Proposition~\ref{prop:f_conv}}\label{sec:f_conv:proof}

The proof of Proposition~\ref{prop:f_conv} relies on the lemma below. We defer the proof of the lemma to Sec.\ref{sec:f_markovian:proof}.

\begin{lem}\label{lem:f_markovian}
We have for all $k\geq\tau_k$
\begin{align*}
\mathbb{E}[\langle\Delta f_k, F(\theta_{k},\hat{V}_{k},\hat{\mu}_{k},s_k,a_k,s_{k+1})\hspace{-2pt}-\hspace{-2pt}\bar{F}(\theta_{k},\hat{V}_{k},\hat{\mu}_{k})\rangle]\leq(20L\hspace{-2pt}+\hspace{-2pt}2|\Acal|) L_F L_{TV} B_F^3 B_G B_H^2\tau_k^2\lambda_{k-\tau_k}.
\end{align*}
\end{lem}

By the update rule of $f_k$, 
\begin{align*}
\Delta f_{k+1}
&=f_{k+1}-\bar{F}(\theta_{k+1},\hat{V}_{k+1},\hat{\mu}_{k+1})\notag\\
&=(1-\lambda_k)f_k+\lambda_k F(\theta_{k},\hat{V}_{k},\hat{\mu}_{k},s_k,a_k,s_{k+1})-\bar{F}(\theta_{k+1},\hat{V}_{k+1},\hat{\mu}_{k+1})\notag\\
&=(1-\lambda_k)f_k+\lambda_k \bar{F}(\theta_{k},\hat{V}_{k},\hat{\mu}_{k})-\bar{F}(\theta_{k+1},\hat{V}_{k+1},\hat{\mu}_{k+1})\notag\\
&\hspace{20pt}+\lambda_k\Big( F(\theta_{k},\hat{V}_{k},\hat{\mu}_{k},s_k,a_k,s_{k+1})-\bar{F}(\theta_{k},\hat{V}_{k},\hat{\mu}_{k})\Big)\notag\\
&=(1-\lambda_k)\Delta f_k+ \Big(\bar{F}(\theta_{k},\hat{V}_{k},\hat{\mu}_{k})-\bar{F}(\theta_{k+1},\hat{V}_{k+1},\hat{\mu}_{k+1})\Big)\notag\\
&\hspace{20pt}+\lambda_k\Big( F(\theta_{k},\hat{V}_{k},\hat{\mu}_{k},s_k,a_k,s_{k+1})-\bar{F}(\theta_{k},\hat{V}_{k},\hat{\mu}_{k})\Big).
\end{align*}

Taking the norm, we have
\begin{align}
&\|\Delta f_{k+1}\|^2\notag\\
&=(1-\lambda_k)^2\|\Delta f_{k}\|^2+\|\bar{F}(\theta_{k},\hat{V}_{k},\hat{\mu}_{k})-\bar{F}(\theta_{k+1},\hat{V}_{k+1},\hat{\mu}_{k+1})\|^2\notag\\
&\hspace{20pt}+\lambda_k^2\|F(\theta_{k},\hat{V}_{k},\hat{\mu}_{k},s_k,a_k,s_{k+1})-\bar{F}(\theta_{k},\hat{V}_{k},\hat{\mu}_{k})\|^2\notag\\
&\hspace{20pt}+(1-\lambda_k)\langle\Delta f_k,\bar{F}(\theta_{k},\hat{V}_{k},\hat{\mu}_{k})-\bar{F}(\theta_{k+1},\hat{V}_{k+1},\hat{\mu}_{k+1})\rangle\notag\\
&\hspace{20pt}+(1-\lambda_k)\lambda_k\langle\Delta f_k,F(\theta_{k},\hat{V}_{k},\hat{\mu}_{k},s_k,a_k,s_{k+1})-\bar{F}(\theta_{k},\hat{V}_{k},\hat{\mu}_{k})\rangle\notag\\
&\hspace{20pt}+\lambda_k\langle \bar{F}(\theta_{k},\hat{V}_{k},\hat{\mu}_{k})-\bar{F}(\theta_{k+1},\hat{V}_{k+1},\hat{\mu}_{k+1}), F(\theta_{k},\hat{V}_{k},\hat{\mu}_{k},s_k,a_k,s_{k+1})-\bar{F}(\theta_{k},\hat{V}_{k},\hat{\mu}_{k})\rangle\notag\\
&\leq (1-\lambda_k)^2\|\Delta f_{k}\|^2+2\|\bar{F}(\theta_{k},\hat{V}_{k},\hat{\mu}_{k})-\bar{F}(\theta_{k+1},\hat{V}_{k+1},\hat{\mu}_{k+1})\|^2\notag\\
&\hspace{20pt}+2\lambda_k^2\|F(\theta_{k},\hat{V}_{k},\hat{\mu}_{k},s_k,a_k,s_{k+1})-\bar{F}(\theta_{k},\hat{V}_{k},\hat{\mu}_{k})\|^2\notag\\
&\hspace{20pt}+\frac{\lambda_k}{2}\|\Delta f_k\|^2+\frac{2}{\lambda_k}\|\bar{F}(\theta_{k},\hat{V}_{k},\hat{\mu}_{k})-\bar{F}(\theta_{k+1},\hat{V}_{k+1},\hat{\mu}_{k+1})\|^2\notag\\
&\hspace{20pt}+(1-\lambda_k)\lambda_k\langle\Delta f_k,F(\theta_{k},\hat{V}_{k},\hat{\mu}_{k},s_k,a_k,s_{k+1})-\bar{F}(\theta_{k},\hat{V}_{k},\hat{\mu}_{k})\rangle\notag\\
&\leq (1-\lambda_k)\|\Delta f_{k}\|^2+(-\frac{\lambda_k}{2}+\lambda_k^2)\|\Delta f_{k}\|^2+\frac{4}{\lambda_k}\|\bar{F}(\theta_{k},\hat{V}_{k},\hat{\mu}_{k})-\bar{F}(\theta_{k+1},\hat{V}_{k+1},\hat{\mu}_{k+1})\|^2\notag\\
&\hspace{20pt}+8B_F^2\lambda_k^2+(1-\lambda_k)\lambda_k\langle\Delta f_k,F(\theta_{k},\hat{V}_{k},\hat{\mu}_{k},s_k,a_k,s_{k+1})-\bar{F}(\theta_{k},\hat{V}_{k},\hat{\mu}_{k})\rangle,\label{prop:f_conv:proof_eq1}
\end{align}
where the final inequality follows from the step size condition $\lambda_k\leq 1$ and the boundedness of operator $F$ which implies
\begin{align*}
\|\bar{F}(\theta_{k},\hat{V}_{k},\hat{\mu}_{k})-\bar{F}(\theta_{k+1},\hat{V}_{k+1},\hat{\mu}_{k+1})\|\leq2B_F.
\end{align*}

Taking the expectation, we can simplify \eqref{prop:f_conv:proof_eq1} as
\begin{align}
&\mathbb{E}[\|\Delta f_{k+1}\|^2]\notag\\
&\leq \mathbb{E}\Big[(1-\lambda_k)\|\Delta f_{k}\|^2+(-\frac{\lambda_k}{2}+\lambda_k^2)\|\Delta f_{k}\|^2+\frac{4}{\lambda_k}\|\bar{F}(\theta_{k},\hat{V}_{k},\hat{\mu}_{k})-\bar{F}(\theta_{k+1},\hat{V}_{k+1},\hat{\mu}_{k+1})\|^2\notag\\
&\hspace{20pt}+8B_F^2\lambda_k^2+(1-\lambda_k)\lambda_k\langle\Delta f_k,F(\theta_{k},\hat{V}_{k},\hat{\mu}_{k},s_k,a_k,s_{k+1})-\bar{F}(\theta_{k},\hat{V}_{k},\hat{\mu}_{k})\rangle\Big]\notag\\
&\leq (1-\lambda_k)\mathbb{E}[\|\Delta f_{k}\|^2]+(-\frac{\lambda_k}{2}+\lambda_k^2)\mathbb{E}[\|\Delta f_{k}\|^2]+8B_F^2\lambda_k^2\notag\\
&\hspace{20pt}+\frac{4L_F^2}{\lambda_k}\mathbb{E}[\left(\|\theta_{k}-\theta_{k+1}\|+\|\hat{V}_{k}-\hat{V}_{k+1}\|+\|\hat{\mu}_{k}-\hat{\mu}_{k+1}\|\right)^2]\notag\\
&\hspace{20pt}+(1-\lambda_k)\lambda_k\cdot(20L+2|\Acal|) L_F L_{TV} B_F^3 B_G B_H^2\tau_k^2\lambda_{k-\tau_k}\notag\\
&\leq (1-\lambda_k)\mathbb{E}[\|\Delta f_{k}\|^2]+(-\frac{\lambda_k}{2}+\lambda_k^2)\mathbb{E}[\|\Delta f_{k}\|^2]+(28L\hspace{-2pt}+\hspace{-2pt}2|\Acal|) L_F L_{TV} B_F^3 B_G B_H^2\tau_k^2\lambda_k\lambda_{k-\tau_k}\notag\\
&\hspace{20pt}+\frac{4L_F^2}{\lambda_k}\mathbb{E}[\left(\alpha_k\|f_k\|+\beta_k\|g_k\|+\xi_k\|h_k\|\right)^2]\notag\\
&\leq (1-\lambda_k)\mathbb{E}[\|\Delta f_{k}\|^2]+(-\frac{\lambda_k}{2}+\lambda_k^2)\mathbb{E}[\|\Delta f_{k}\|^2]+(28L\hspace{-2pt}+\hspace{-2pt}2|\Acal|) L_F L_{TV} B_F^3 B_G B_H^2\tau_k^2\lambda_k\lambda_{k-\tau_k}\notag\\
&\hspace{20pt}+\mathbb{E}\Big[\frac{12L_F^2 \alpha_k}{\lambda_k}\left(\|\Delta f_k\|+L_F\sqrt{\varepsilon_k^{V}}+L_F(L_V+1)\sqrt{\varepsilon_k^{\mu}}+\sqrt{\varepsilon_k^{\pi}}\right)^2\notag\\
&\hspace{20pt}+\frac{12L_F^2 \beta_k}{\lambda_k}\left(\|\Delta g_k\|+L_G\sqrt{\varepsilon_k^{V}}+L_G\sqrt{\varepsilon_k^{J}}\right)^2+\frac{12L_F^2 \xi_k}{\lambda_k}\left(L_H\|\Delta h_k\|+\sqrt{\varepsilon_k^{\mu}}\right)^2\Big],\label{prop:f_conv:proof_eq2}
\end{align}
where the second inequality plugs in the result of Lemma~\ref{lem:f_markovian} and bounds $\|\bar{F}(\theta_{k},\hat{V}_{k},\hat{\mu}_{k})-\bar{F}(\theta_{k+1},\hat{V}_{k+1},\hat{\mu}_{k+1})\|^2$ using the Lipschitz condition established in Lemma~\ref{lem:Lipschitz_operators}.

The sum of the last three terms can be bounded as
\begin{align}
&\frac{12L_F^2 \alpha_k}{\lambda_k}\left(\|\Delta f_k\|+L_F\sqrt{\varepsilon_k^{V}}+L_F(L_V+1)\sqrt{\varepsilon_k^{\mu}}+\sqrt{\varepsilon_k^{\pi}}\right)^2\notag\\
&\hspace{20pt}+\frac{12L_F^2 \beta_k}{\lambda_k}\left(\|\Delta g_k\|+L_G\sqrt{\varepsilon_k^{V}}+L_G\sqrt{\varepsilon_k^{J}}\right)^2+\frac{12L_F^2 \xi_k}{\lambda_k}\left(L_H\|\Delta h_k\|+\sqrt{\varepsilon_k^{\mu}}\right)^2\notag\\
&\leq \frac{48L_F^2 \alpha_k^2}{\lambda_k}\|\Delta f_k\|^2+\frac{48L_F^4 \alpha_k^2}{\lambda_k}\varepsilon_k^{V}+\frac{192L_F^4 L_V^2 \alpha_k^2}{\lambda_k}\varepsilon_k^{\mu}+\frac{48L_F^2 \alpha_k^2}{\lambda_k}\varepsilon_k^{\pi}\notag\\
&\hspace{20pt}+\frac{36L_F^2\beta_k^2}{\lambda_k}\|\Delta g_k\|^2+\frac{48L_F^2 L_G^2 \beta_k^2}{\lambda_k}\varepsilon_k^{V}+\frac{48L_F^2 L_G^2 \beta_k^2}{\lambda_k}\varepsilon_k^{J}\notag\\
&\hspace{20pt}+\frac{24L_F^2 L_H^2 \xi_k^2}{\lambda_k}\|\Delta h_k\|^2+\frac{24L_F^2 \xi_k^2}{\lambda_k}\varepsilon_k^{\mu}\notag\\
&\leq \frac{48L_F^2 \alpha_k^2}{\lambda_k}\|\Delta f_k\|^2+\frac{36L_F^2\beta_k^2}{\lambda_k}\|\Delta g_k\|^2+\frac{24L_F^2 L_H^2 \xi_k^2}{\lambda_k}\|\Delta h_k\|^2+\frac{48L_F^2 \alpha_k^2}{\lambda_k}\varepsilon_k^{\pi}\notag\\
&\hspace{20pt}+\frac{216L_F^4 L_V^2 \xi_k^2}{\lambda_k}\varepsilon_k^{\mu}+\frac{96L_F^4 L_G^2\beta_k^2}{\lambda_k}\varepsilon_k^{V}+\frac{48L_F^2 L_G^2 \beta_k^2}{\lambda_k}\varepsilon_k^{J}.\label{prop:f_conv:proof_eq3}
\end{align}

Combining \eqref{prop:f_conv:proof_eq2} and \eqref{prop:f_conv:proof_eq3}, we get
\begin{align*}
&\mathbb{E}[\|\Delta f_{k+1}\|^2]\notag\\
&\leq (1-\lambda_k)\mathbb{E}[\|\Delta f_{k}\|^2]+(-\frac{\lambda_k}{2}+\lambda_k^2)\mathbb{E}[\|\Delta f_{k}\|^2]+(28L+2|\Acal|) L_F L_{TV} B_F^3 B_G B_H^2\tau_k^2\lambda_k\lambda_{k-\tau_k}\notag\\
&\hspace{20pt}+\mathbb{E}\Big[\frac{48L_F^2 \alpha_k^2}{\lambda_k}\|\Delta f_k\|^2+\frac{36L_F^2\beta_k^2}{\lambda_k}\|\Delta g_k\|^2+\frac{24L_F^2 L_H^2 \xi_k^2}{\lambda_k}\|\Delta h_k\|^2+\frac{48L_F^2 \alpha_k^2}{\lambda_k}\varepsilon_k^{\pi}\notag\\
&\hspace{20pt}+\frac{216L_F^4 L_V^2 \xi_k^2}{\lambda_k}\varepsilon_k^{\mu}+\frac{96L_F^4 L_G^2\beta_k^2}{\lambda_k}\varepsilon_k^{V}+\frac{48L_F^2 L_G^2 \beta_k^2}{\lambda_k}\varepsilon_k^{J}\Big]\notag\\
&= (1-\lambda_k)\mathbb{E}[\|\Delta f_{k}\|^2]+(-\frac{\lambda_k}{2}+\lambda_k^2+\frac{48L_F^2 \alpha_k^2}{\lambda_k})\mathbb{E}[\|\Delta f_{k}\|^2]\notag\\
&\hspace{20pt}+\frac{36L_F^2\beta_k^2}{\lambda_k}\mathbb{E}[\|\Delta g_k\|^2]+\frac{24L_F^2 L_H^2 \xi_k^2}{\lambda_k}\mathbb{E}[\|\Delta h_k\|^2]+\frac{48L_F^2 \alpha_k^2}{\lambda_k}\mathbb{E}[\varepsilon_k^{\pi}]+\frac{216L_F^4 L_V^2 \xi_k^2}{\lambda_k}\mathbb{E}[\varepsilon_k^{\mu}]\notag\\
&\hspace{20pt}+\frac{96L_F^4 L_G^2\beta_k^2}{\lambda_k}\mathbb{E}[\varepsilon_k^{V}]+\frac{48L_F^2 L_G^2 \beta_k^2}{\lambda_k}\mathbb{E}[\varepsilon_k^{J}]+(28L+2|\Acal|) L_F L_{TV} B_F^3 B_G B_H^2\tau_k^2\lambda_k\lambda_{k-\tau_k}.
\end{align*}

\qed

\subsection{Proof of Proposition~\ref{prop:meanfield_conv}}\label{sec:meanfield_conv:proof}

We first introduce the following lemma, which will be used in the proof of Proposition~\ref{prop:meanfield_conv}. The proof of Lemma~\ref{lem:H_strongmonotone} is presented in Sec.\ref{sec:H_strongmonotone:proof}.
\begin{lem}\label{lem:H_strongmonotone}
Under Assumption~\ref{assump:nu}, we have for any policy parameter $\theta$ and mean field $\mu$
\begin{align*}
\langle \mu-\mu^{\star}(\pi_{\theta}),\bar{H}(\theta,\mu)-\bar{H}(\theta,\mu^{\star}(\pi_{\theta}))\rangle\leq-(1-\delta)\|\mu-\mu^{\star}(\pi_{\theta})\|^2.
\end{align*}
\end{lem}

By the definition of $\varepsilon_{k}^{\mu}$,
\begin{align}
\varepsilon_{k+1}^{\mu}&=\|\hat{\mu}_{k+1}-\mu^{\star}(\pi_{\theta_{k+1}})\|^2\notag\\
&=\|\Pi_{\Delta_{\Scal}}(\hat{\mu}_k+\xi_k h_k)-\mu^{\star}(\pi_{\theta_{k+1}})\|^2\notag\\
&\leq\|\hat{\mu}_k+\xi_k h_k-\mu^{\star}(\pi_{\theta_{k+1}})\|^2\notag\\
&=\|\hat{\mu}_k-\mu^{\star}(\pi_{\theta_{k}})+\xi_k \Delta h_k+\xi_k \bar{H}(\theta_k,\hat{\mu}_k)-\left(\mu^{\star}(\pi_{\theta_{k+1}})-\mu^{\star}(\pi_{\theta_{k}})\right)\|^2\notag\\
&=\|\hat{\mu}_k-\mu^{\star}(\pi_{\theta_{k}})+\xi_k \bar{H}(\theta_k,\hat{\mu}_k)\|^2+\xi_k^2\|\Delta h_k\|^2+\|\mu^{\star}(\pi_{\theta_{k+1}})-\mu^{\star}(\pi_{\theta_{k}})\|^2\notag\\
&\hspace{20pt}+2\xi_k\langle \hat{\mu}_k-\mu^{\star}(\pi_{\theta_{k}})+\xi_k \bar{H}(\theta_k,\hat{\mu}_k),\Delta h_k\rangle+2\langle\hat{\mu}_k-\mu^{\star}(\pi_{\theta_{k}})+\xi_k \bar{H}(\theta_k,\hat{\mu}_k),\mu^{\star}(\pi_{\theta_{k+1}})-\mu^{\star}(\pi_{\theta_{k}})\rangle\notag\\
&\hspace{20pt}+2\xi_k\langle\Delta h_k,\mu^{\star}(\pi_{\theta_{k+1}})-\mu^{\star}(\pi_{\theta_{k}})\rangle,\label{prop:meanfield_conv:proof_eq1}
\end{align}
where the first inequality is due to the fact that projection to a convex set is a non-expansive operator.

To bound the first term of \eqref{prop:meanfield_conv:proof_eq1}, 
\begin{align}
&\|\hat{\mu}_k-\mu^{\star}(\pi_{\theta_{k}})+\xi_k \bar{H}(\theta_k,\hat{\mu}_k)\|^2\notag\\
&=\|\hat{\mu}_k-\mu^{\star}(\pi_{\theta_{k}})+\xi_k \big(\bar{H}(\theta_k,\hat{\mu}_k)-\bar{H}(\theta_k,\mu^{\star}(\pi_{\theta_{k}}))\big)\|^2\notag\\
&=\|\hat{\mu}_k-\mu^{\star}(\pi_{\theta_{k}})\|^2+\xi_k^2\|\bar{H}(\theta_k,\hat{\mu}_k)-\bar{H}(\theta_k,\mu^{\star}(\pi_{\theta_{k}}))\|^2\notag\\
&\hspace{20pt}+2\xi_k\langle \hat{\mu}_k-\mu^{\star}(\pi_{\theta_{k}}), \bar{H}(\theta_k,\hat{\mu}_k)-\bar{H}(\theta_k,\mu^{\star}(\pi_{\theta_{k}}))\rangle\notag\\
&\leq \|\hat{\mu}_k-\mu^{\star}(\pi_{\theta_{k}})\|^2+L_H^2\xi_k^2 \|\hat{\mu}_k-\mu^{\star}(\pi_{\theta_{k}})\|^2 -(1-\delta)\xi_k\|\hat{\mu}_k-\mu^{\star}(\pi_{\theta_{k}})\|^2\notag\\
&\leq (1-\frac{(1-\delta)\xi_k}{2})\varepsilon_k^{\mu},\label{prop:meanfield_conv:proof_eq2}
\end{align}
where the first equation uses $\bar{H}(\theta,\mu^{\star}(\pi_{\theta}))=0$ for any $\theta$, the first inequality is a result of Lemma~\ref{lem:H_strongmonotone} and the Lipschitz continuity of $\bar{H}$, and the second inequality follows from the step size condition $\xi_k\leq\beta_k\leq\frac{1-\delta}{2L_H^2}$.

We next treat the second and third term of \eqref{prop:meanfield_conv:proof_eq1} using the fact that $\|h_k\|\leq B_H$, $\|\bar{H}(\theta_k,\hat{\mu}_k)\|\leq B_H$, $\|f_k\|\leq B_F$ and that the operator $\mu^{\star}$ is Lipschitz
\begin{align}
\xi_k^2\|\Delta h_k\|^2+\|\mu^{\star}(\pi_{\theta_{k+1}})-\mu^{\star}(\pi_{\theta_{k}})\|^2&\leq 2\xi_k^2\|h_k\|^2+2\xi_k^2\|\bar{H}(\theta_k,\hat{\mu}_k)\|^2+L\|\pi_{\theta_{k+1}}-\pi_{\theta_{k}}\|^2\notag\\
&\leq 4B_H^2\xi_k^2 + L^2\|f_k\|^2\notag\\
&\leq 4B_H^2\xi_k^2 + L^2 B_F^2\alpha_k^2.
\end{align}

The fourth term of \eqref{prop:meanfield_conv:proof_eq1} can be bounded leveraging the result in \eqref{prop:meanfield_conv:proof_eq2} as follows
\begin{align}
&2\xi_k\langle \hat{\mu}_k-\mu^{\star}(\pi_{\theta_{k}})+\xi_k \bar{H}(\theta_k,\hat{\mu}_k),\Delta h_k\rangle\notag\\
&\leq \frac{(1-\delta)\xi_k}{8}\|\hat{\mu}_k-\mu^{\star}(\pi_{\theta_{k}})+\xi_k \bar{H}(\theta_k,\hat{\mu}_k)\|^2+\frac{8\xi_k}{1-\delta}\|\Delta h_k\|^2\notag\\
&\leq \frac{(1-\delta)\xi_k}{8} \cdot (1-\frac{(1-\delta)\xi_k}{2})\varepsilon_k^{\mu}+\frac{8\xi_k}{1-\delta}\|\Delta h_k\|^2\notag\\
&\leq \frac{(1-\delta)\xi_k}{8}\varepsilon_k^{\mu}+\frac{8\xi_k}{1-\delta}\|\Delta h_k\|^2.
\end{align}

Similarly, for the fifth term of \eqref{prop:meanfield_conv:proof_eq1}, we have
\begin{align}
&2\langle\hat{\mu}_k-\mu^{\star}(\pi_{\theta_{k}})+\xi_k \bar{H}(\theta_k,\hat{\mu}_k),\mu^{\star}(\pi_{\theta_{k+1}})-\mu^{\star}(\pi_{\theta_{k}})\rangle \notag\\
&\leq \frac{(1-\delta)\xi_k}{8}\|\hat{\mu}_k-\mu^{\star}(\pi_{\theta_{k}})+\xi_k \bar{H}(\theta_k,\hat{\mu}_k)\|^2+\frac{8}{(1-\delta)\xi_k}\|\mu^{\star}(\pi_{\theta_{k+1}})-\mu^{\star}(\pi_{\theta_{k}})\|^2\notag\\
&\leq \frac{(1-\delta)\xi_k}{8}\varepsilon_k^{\mu}+\frac{8L^2}{(1-\delta)\xi_k}\|\pi_{\theta_{k+1}}-\pi_{\theta_{k}}\|^2\notag\\
&\leq \frac{(1-\delta)\xi_k}{8}\varepsilon_k^{\mu}+\frac{8L^2\alpha_k^2}{(1-\delta)\xi_k}\|f_k\|^2\notag\\
&\leq \frac{(1-\delta)\xi_k}{8}\varepsilon_k^{\mu}+\frac{8L^2\alpha_k^2}{(1-\delta)\xi_k}\left(\|\Delta f_k\|+L_F\sqrt{\varepsilon_k^{V}}+L_F(L_V+1)\sqrt{\varepsilon_k^{\mu}}+\sqrt{\varepsilon_k^{\pi}}\right)^2\notag\\
&\leq \frac{(1-\delta)\xi_k}{8}\varepsilon_k^{\mu}+\frac{32L^2\alpha_k^2}{(1-\delta)\xi_k}\left(\|\Delta f_k\|^2+L_F^2\varepsilon_k^{V}+4L_F^2 L_V^2\varepsilon_k^{\mu}+\varepsilon_k^{\pi}\right),
\end{align}
where the fourth inequality follows from Lemma~\ref{lem:bound_f_g}.

The final term of \eqref{prop:meanfield_conv:proof_eq1} can be bounded simply with the Cauchy-Schwarz inequality
\begin{align}
2\xi_k\langle\Delta h_k,\mu^{\star}(\pi_{\theta_{k+1}})-\mu^{\star}(\pi_{\theta_{k}})\rangle&\leq2\xi_k\|\Delta h_k\|\|\mu^{\star}(\pi_{\theta_{k+1}})-\mu^{\star}(\pi_{\theta_{k}})\|\notag\\
&\leq 4B_H\xi_k\cdot L\|\pi_{\theta_{k+1}}-\pi_{\theta_{k}}\|\notag\\
&\leq 4L B_F B_H\alpha_k\xi_k.
\label{prop:meanfield_conv:proof_eq3}
\end{align}

Plugging \eqref{prop:meanfield_conv:proof_eq2}-\eqref{prop:meanfield_conv:proof_eq3} into \eqref{prop:meanfield_conv:proof_eq1}, we get
\begin{align*}
\varepsilon_{k+1}^{\mu}&\leq (1-\frac{(1-\delta)\xi_k}{2})\varepsilon_k^{\mu}+4B_H^2\xi_k^2 + L^2 B_F^2\alpha_k^2+\frac{(1-\delta)\xi_k}{8}\varepsilon_k^{\mu}+\frac{8\xi_k}{1-\delta}\|\Delta h_k\|^2\notag\\
&\hspace{20pt}+\frac{(1-\delta)\xi_k}{8}\varepsilon_k^{\mu}+\frac{32L^2\alpha_k^2}{(1-\delta)\xi_k}\left(\|\Delta f_k\|^2+L_F^2\varepsilon_k^{V}+4L_F^2 L_V^2 \varepsilon_k^{\mu}+\varepsilon_k^{\pi}\right)+4L B_F B_H\alpha_k\xi_k\notag\\
&\leq (1-\frac{(1-\delta)\xi_k}{8})\varepsilon_k^{\mu}+\frac{8\xi_k}{1-\delta}\|\Delta h_k\|^2+4B_H^2\xi_k^2 + L^2 B_F^2\alpha_k^2\hspace{-2pt}+\hspace{-2pt}\frac{32L^2\alpha_k^2}{(1-\delta)\xi_k}\left(\|\Delta f_k\|^2\hspace{-2pt}+\hspace{-2pt}L_F^2\varepsilon_k^{V}\hspace{-2pt}+\hspace{-2pt}\varepsilon_k^{\pi}\right)\notag\\
&\hspace{20pt}+4L B_F B_H\alpha_k\xi_k+(-\frac{(1-\delta)\xi_k}{8}+\frac{128L^2 L_F^2L_V^2\alpha_k^2}{(1-\delta)\xi_k})\varepsilon_k^{\mu}\notag\\
&\leq (1-\frac{(1-\delta)\xi_k}{8})\varepsilon_k^{\mu}+\frac{8\xi_k}{1-\delta}\|\Delta h_k\|^2+\frac{32L^2\alpha_k^2}{(1-\delta)\xi_k}\left(\|\Delta f_k\|^2+L_F^2\varepsilon_k^{V}+\varepsilon_k^{\pi}\right)+9L^2 B_F^2 B_H^2\xi_k^2,
\end{align*}
where the last inequality is a result of the step size condition $\alpha_k\leq\xi_k$ and $\alpha_k\leq\frac{1-\delta}{16L L_F L_V}\xi_k$.

\qed

\subsection{Proof of Proposition~\ref{prop:h_conv}}\label{sec:h_conv:proof}

The proof of Proposition~\ref{prop:h_conv} uses an intermediate result established in the lemma below. We defer the proof of the lemma to Sec.\ref{sec:h_markovian:proof}.

\begin{lem}\label{lem:h_markovian}
We have for all $k\geq\tau_k$
\begin{align*}
\mathbb{E}[\langle\Delta h_k, e_{s_k}-\mathbb{E}_{s\sim\nu^{\pi_{\theta_k},\,\hat{\mu}_k}}[e_s]\rangle]\leq16L B_F B_H^2\tau_k^2\lambda_{k-\tau_k}.
\end{align*}
\end{lem}

By the update rule of $h_k$, 
\begin{align*}
\Delta h_{k+1}&=h_{k+1}-\bar{H}(\theta_{k+1},\hat{\mu}_{k+1})\notag\\
&=(1-\lambda_k)h_k+\lambda_k(e_{s_k}-\hat{\mu}_k)-\bar{H}(\theta_{k+1},\hat{\mu}_{k+1})\notag\\
&=(1-\lambda_k)h_k+\lambda_k \bar{H}(\theta_{k},\hat{\mu}_{k})-\bar{H}(\theta_{k+1},\hat{\mu}_{k+1})+\lambda_k\Big((e_{s_k}-\hat{\mu}_k)-\bar{H}(\theta_{k},\hat{\mu}_{k})\Big)\notag\\
&=(1-\lambda_k)\Delta h_k+ \Big(\bar{H}(\theta_{k},\hat{\mu}_{k})-\bar{H}(\theta_{k+1},\hat{\mu}_{k+1})\Big)+\lambda_k\Big((e_{s_k}-\hat{\mu}_k)-\bar{H}(\theta_{k},\hat{\mu}_{k})\Big).
\end{align*}
This implies
\begin{align*}
&\|\Delta h_{k+1}\|^2\notag\\
&=(1-\lambda_k)^2\|\Delta h_k\|^2+\|\bar{H}(\theta_{k},\hat{\mu}_{k})-\bar{H}(\theta_{k+1},\hat{\mu}_{k+1})\|^2+\lambda_k^2\|(e_{s_k}-\hat{\mu}_k)-\bar{H}(\theta_{k},\hat{\mu}_{k})\|^2\notag\\
&\hspace{20pt}+(1-\lambda_k)\langle\Delta h_k,\bar{H}(\theta_{k},\hat{\mu}_{k})-\bar{H}(\theta_{k+1},\hat{\mu}_{k+1})\rangle \notag\\
&\hspace{20pt}+ (1-\lambda_k)\lambda_k\langle\Delta h_k,(e_{s_k}-\hat{\mu}_k)-\bar{H}(\theta_{k},\hat{\mu}_{k})\rangle\notag\\
&\hspace{20pt}+\lambda_k\langle \bar{H}(\theta_{k},\hat{\mu}_{k})-\bar{H}(\theta_{k+1},\hat{\mu}_{k+1}),(e_{s_k}-\hat{\mu}_k)-\bar{H}(\theta_{k},\hat{\mu}_{k})\rangle\notag\\
&\leq (1-\lambda_k)^2\|\Delta h_k\|^2+2\|\bar{H}(\theta_{k},\hat{\mu}_{k})-\bar{H}(\theta_{k+1},\hat{\mu}_{k+1})\|^2+2\lambda_k^2\|(e_{s_k}-\hat{\mu}_k)-\bar{H}(\theta_{k},\hat{\mu}_{k})\|^2\notag\\
&\hspace{20pt}+\frac{\lambda_k}{2}\|\Delta h_k\|^2+\frac{2}{\lambda_k}\|\bar{H}(\theta_{k},\hat{\mu}_{k})-\bar{H}(\theta_{k+1},\hat{\mu}_{k+1})\|^2\notag\\
&\hspace{20pt}+(1-\lambda_k)\lambda_k\langle\Delta h_k, e_{s_k}-\mathbb{E}_{s\sim\nu^{\pi_{\theta_k},\,\hat{\mu}_k}}[e_s]\rangle\notag\\
&\leq(1-\lambda_k)\|\Delta h_k\|^2+(-\frac{\lambda_k}{2}+\lambda_k^2)\|\Delta h_k\|^2+\frac{4}{\lambda_k}\|\bar{H}(\theta_{k},\hat{\mu}_{k})-\bar{H}(\theta_{k+1},\hat{\mu}_{k+1})\|^2\notag\\
&\hspace{20pt}+(1-\lambda_k)\lambda_k\langle\Delta h_k, e_{s_k}-\mathbb{E}_{s\sim\nu^{\pi_{\theta_k},\,\hat{\mu}_k}}[e_s]\rangle+8B_H\lambda_k^2,
\end{align*}
where the final inequality follows from the step size choice $\lambda_k\leq 1$.
Taking the expectation and applying Lemma~\ref{lem:h_markovian} and the Lipschitz continuity of operator $\bar{H}$, we further have
\begin{align*}
&\mathbb{E}[\|\Delta h_{k+1}\|^2]\notag\\
&\leq(1\hspace{-2pt}-\hspace{-2pt}\lambda_k)\mathbb{E}[\|\Delta h_k\|^2] + (-\frac{\lambda_k}{2}\hspace{-2pt}+\hspace{-2pt}\lambda_k^2)\mathbb{E}[\|\Delta h_k\|^2]+\frac{4}{\lambda_k}\mathbb{E}[\left(L_H\|\theta_{k}-\theta_{k+1}\|+L_H\|\hat{\mu}_{k}-\hat{\mu}_{k+1}\|\right)^2]\notag\\
&\hspace{20pt}+(1-\lambda_k)\lambda_k\cdot 16L B_F B_H^2\tau_k^2\lambda_{k-\tau_k}+8B_H\lambda_k^2\notag\\
&\leq(1-\lambda_k)\mathbb{E}[\|\Delta h_k\|^2]+(-\frac{\lambda_k}{2}+\lambda_k^2)\mathbb{E}[\|\Delta h_k\|^2]+\frac{8L_H^2}{\lambda_k}\mathbb{E}[\alpha_k^2\|f_k\|^2+\xi_k^2\|h_k\|^2]\notag\\
&\hspace{20pt}+ 16L B_F B_H^2\tau_k^2\lambda_k\lambda_{k-\tau_k}+8B_H\lambda_k^2\notag\\
&\leq(1-\lambda_k)\mathbb{E}[\|\Delta h_k\|^2]+(-\frac{\lambda_k}{2}+\lambda_k^2)\mathbb{E}[\|\Delta h_k\|^2]\notag\\
&\hspace{20pt}+\frac{8L_H^2\alpha_k^2}{\lambda_k}\mathbb{E}[\big(\|\Delta f_k\|+L_F\sqrt{\varepsilon_k^{V}}+L_F(L_V+1)\sqrt{\varepsilon_k^{\mu}}+\sqrt{\varepsilon_k^{\pi}}\big)^2]\notag\\
&\hspace{20pt} +\frac{8L_H^2\xi_k^2}{\lambda_k}\mathbb{E}[\big(\|\Delta h_k\|+L_H\sqrt{\epsilon_k^{\mu}}\big)^2] + 16L B_F B_H^2\tau_k^2\lambda_k\lambda_{k-\tau_k}+8B_H\lambda_k^2\notag\\
&\leq(1-\lambda_k)\mathbb{E}[\|\Delta h_k\|^2]+(-\frac{\lambda_k}{2}+\lambda_k^2)\mathbb{E}[\|\Delta h_k\|^2]\notag\\
&\hspace{20pt}+\frac{32L_H^2\alpha_k^2}{\lambda_k}\mathbb{E}[\|\Delta f_k\|^2+L_F^2\varepsilon_k^{V}+4L_F^2 L_V^2\varepsilon_k^{\mu}+\varepsilon_k^{\pi}]\notag\\
&\hspace{20pt}+\frac{16L_H^2\xi_k^2}{\lambda_k}\mathbb{E}[\|\Delta h_k\|^2+L_H^2\epsilon_k^{\mu}] + 24L B_F B_H^2\tau_k^2\lambda_k\lambda_{k-\tau_k}\notag\\
&\leq(1-\lambda_k)\mathbb{E}[\|\Delta h_k\|^2]+(-\frac{\lambda_k}{2}+\lambda_k^2+\frac{16L_H^2\xi_k^2}{\lambda_k})\mathbb{E}[\|\Delta h_k\|^2]+\frac{32L_H^2 \alpha_k^2}{\lambda_k}\mathbb{E}[\|\Delta f_k\|^2]\notag\\
&\hspace{20pt}+\frac{32L_H^2 L_F^2\alpha_k^2}{\lambda_k}\mathbb{E}[\varepsilon_k^{V}]+\frac{144 L_F^2 L_V^2 L_H^4 \xi_k^2}{\lambda_k}\mathbb{E}[\varepsilon_k^{\mu}]+\frac{32L_H^2 \alpha_k^2}{\lambda_k}\mathbb{E}[\varepsilon_k^{\pi}] + 24L B_F B_H^2\tau_k^2\lambda_k\lambda_{k-\tau_k},
\end{align*}
where the third inequality bounds $\|f_k\|$ and $\|h_k\|$ with Lemma~\ref{lem:bound_f_g}. The step size condition $\alpha_k\leq\xi_k$ is used a few times to simplify and combine terms.

\qed

\subsection{Proof of Proposition~\ref{prop:V_conv}}\label{sec:V_conv:proof}

We use the following lemma in our analysis. The proof of the lemma is deferred to Sec.\ref{sec:G_strongmonotone:proof}.
\begin{lem}\label{lem:G_strongmonotone}
Under Assumption~\ref{assump:ergodic}, it holds for any $\theta$, $\mu$, and $V$ that
\begin{align*}
\left\langle \left[\begin{array}{c}
\Pi_{\Ecal_\perp}(V \hspace{-2pt}-\hspace{-2pt} V^{\pi_{\theta},\,\mu})\\
J - J(\pi_{\theta},\mu)
\end{array}\right],\left[\begin{array}{c}
\Pi_{\Ecal_\perp}\bar{G}^V(\theta,V,J,\mu)\\
\bar{G}^J(\theta,J,\mu)
\end{array}\right]\right\rangle \hspace{-2pt}\leq\hspace{-2pt} -\frac{\gamma}{2}(\|\Pi_{\Ecal_\perp}(V \hspace{-2pt}-\hspace{-2pt} V^{\pi_{\theta},\,\mu})\|^2 \hspace{-2pt}+\hspace{-2pt} (J \hspace{-2pt}-\hspace{-2pt} J(\pi_{\theta},\mu))^2),
\end{align*}
where $\gamma\in(0,1)$ is the discount factor in Lemma~\ref{lem:negative_drift}.
\end{lem}

By the definition of $\varepsilon_{k}^{V}$,
\begin{align}
&\varepsilon_{k+1}^{V}+\varepsilon_{k+1}^{J}\notag\\
&=\left\|\left[\begin{array}{c}
\Pi_{\Ecal_\perp}(\hat{V}_{k+1}-V^{\pi_{\theta_{k+1}},\,\hat{\mu}_{k+1}})\\
\hat{J}_{k+1} - J(\pi_{\theta_{k+1}},\hat{\mu}_{k+1})
\end{array}\right]\right\|^2\notag\\
&=
\left\|\left[\begin{array}{c}
\Pi_{\Ecal_\perp}(\Pi_{B_V}(\hat{V}_{k}+\beta_k g_k^V)-V^{\pi_{\theta_{k+1}},\,\hat{\mu}_{k+1}})\\
\Pi_{[0,1]}(\hat{J}_{k} +\beta_k g_k^J) - J(\pi_{\theta_{k+1}},\hat{\mu}_{k+1})
\end{array}\right]\right\|^2\notag\\
&\leq
\left\|\left[\begin{array}{c}
\Pi_{\Ecal_\perp}(\hat{V}_{k}+\beta_k g_k^V-V^{\pi_{\theta_{k+1}},\,\hat{\mu}_{k+1}})\\
\hat{J}_{k} +\beta_k g_k^J - J(\pi_{\theta_{k+1}},\hat{\mu}_{k+1})
\end{array}\right]\right\|^2\notag\\
&=\left\|\left[\begin{array}{c}\Pi_{\Ecal_\perp}\hspace{-4pt}\left(\hat{V}_{k}-V^{\pi_{\theta_{k}},\,\hat{\mu}_{k}}+\beta_k \bar{G}^V(\theta_k,\hat{V}_k,\hat{J}_k,\hat{\mu}_k)+\beta_k \Delta g_k^V -\left(V^{\pi_{\theta_{k+1}},\,\hat{\mu}_{k+1}}\hspace{-2pt}-\hspace{-2pt}V^{\pi_{\theta_{k}},\,\hat{\mu}_{k}}\right)\right)\\
\hat{J}_k-J(\pi_{\theta_{k}},\hat{\mu}_{k})+\beta_k \bar{G}^J(\theta_k,\hat{J}_k,\hat{\mu}_k)+\beta_k \Delta g_k^J -(J(\pi_{\theta_{k+1}},\hat{\mu}_{k+1})-J(\pi_{\theta_{k}},\hat{\mu}_{k}))
\end{array}\right]\right\|^2\notag\\
&\leq\left\|\left[\begin{array}{c}\Pi_{\Ecal_\perp}(\hat{V}_{k}-V^{\pi_{\theta_{k}},\,\hat{\mu}_{k}})\\
\hat{J}_k-J(\pi_{\theta_{k}},\hat{\mu}_{k})
\end{array}\right]+\beta_k \left[\begin{array}{c}
\Pi_{\Ecal_\perp}\bar{G}^V(\theta_k,\hat{V}_k,\hat{J}_k,\hat{\mu}_k)\\
\bar{G}^J(\theta_k,\hat{J}_k,\hat{\mu}_k)
\end{array}\right]\right\|^2+\beta_k^2\|\Delta g_k\|^2\notag\\
&\hspace{20pt}+\left\|\left[\begin{array}{c}
\Pi_{\Ecal_\perp}(V^{\pi_{\theta_{k+1}},\,\hat{\mu}_{k+1}}-V^{\pi_{\theta_{k}},\,\hat{\mu}_{k}})\\
J(\pi_{\theta_{k+1}},\hat{\mu}_{k+1})-J(\pi_{\theta_{k}},\hat{\mu}_{k})
\end{array}\right]\right\|^2\notag\\
&\hspace{20pt}+2\beta_k\left\langle \left[\begin{array}{c}\Pi_{\Ecal_\perp}(\hat{V}_{k}-V^{\pi_{\theta_{k}},\,\hat{\mu}_{k}})\\
\hat{J}_k-J(\pi_{\theta_{k}},\hat{\mu}_{k})
\end{array}\right]+\beta_k \left[\begin{array}{c}
\Pi_{\Ecal_\perp}\bar{G}^V(\theta_k,\hat{V}_k,\hat{J}_k,\hat{\mu}_k)\\
\bar{G}^J(\theta_k,\hat{J}_k,\hat{\mu}_k)
\end{array}\right],\left[\begin{array}{c}\Pi_{\Ecal_\perp}\Delta g_k^V\\
\Delta g_k^J\end{array}
\right]\right\rangle\notag\\
&\hspace{20pt}+2\left\langle\left[\begin{array}{c}\Pi_{\Ecal_\perp}(\hat{V}_{k}-V^{\pi_{\theta_{k}},\,\hat{\mu}_{k}})\\
\hat{J}_k-J(\pi_{\theta_{k}},\hat{\mu}_{k})
\end{array}\right]+\beta_k \left[\begin{array}{c}
\Pi_{\Ecal_\perp}\bar{G}^V(\theta_k,\hat{V}_k,\hat{J}_k,\hat{\mu}_k)\\
\bar{G}^J(\theta_k,\hat{J}_k,\hat{\mu}_k)
\end{array}\right],\left[\begin{array}{c}
\Pi_{\Ecal_\perp}(V^{\pi_{\theta_{k+1}},\,\hat{\mu}_{k+1}}-V^{\pi_{\theta_{k}},\,\hat{\mu}_{k}})\\
J(\pi_{\theta_{k+1}},\hat{\mu}_{k+1})-J(\pi_{\theta_{k}},\hat{\mu}_{k})
\end{array}\right]\right\rangle\notag\\
&\hspace{20pt}+2\beta_k\left\langle\Delta g_k, \left[\begin{array}{c}
\Pi_{\Ecal_\perp}(V^{\pi_{\theta_{k+1}},\,\hat{\mu}_{k+1}}-V^{\pi_{\theta_{k}},\,\hat{\mu}_{k}})\\
J(\pi_{\theta_{k+1}},\hat{\mu}_{k+1})-J(\pi_{\theta_{k}},\hat{\mu}_{k})
\end{array}\right] \right\rangle,\label{prop:V_conv:proof_eq1}
\end{align}
where the last inequality follows from the fact that $\Pi_{\Ecal_\perp}$ has all singular values smaller than or equal to 1.

To bound the first term of \eqref{prop:V_conv:proof_eq1}, 
\begin{align}
&\left\|\left[\begin{array}{c}\Pi_{\Ecal_\perp}(\hat{V}_{k}-V^{\pi_{\theta_{k}},\,\hat{\mu}_{k}})\\
\hat{J}_k-J(\pi_{\theta_{k}},\hat{\mu}_{k})
\end{array}\right]+\beta_k \left[\begin{array}{c}
\Pi_{\Ecal_\perp}\bar{G}^V(\theta_k,\hat{V}_k,\hat{J}_k,\hat{\mu}_k)\\
\bar{G}^J(\theta_k,\hat{J}_k,\hat{\mu}_k)
\end{array}\right]\right\|^2\notag\\
&\leq \|\Pi_{\Ecal_\perp}(\hat{V}_{k}-V^{\pi_{\theta_{k}},\,\hat{\mu}_{k}})\|^2+(\hat{J}_k-J(\pi_{\theta_{k}},\hat{\mu}_{k}))^2+\beta_k^2\|\Pi_{\Ecal_\perp}\bar{G}^V(\theta_k,\hat{V}_k,\hat{J}_k,\hat{\mu}_k)\|^2\notag\\
&\hspace{20pt}+\beta_k\Big(\bar{G}^J(\theta_k,\hat{J}_k,\hat{\mu}_k)\Big)^2+2\beta_k\left\langle\left[\begin{array}{c}\Pi_{\Ecal_\perp}(\hat{V}_{k}-V^{\pi_{\theta_{k}},\,\hat{\mu}_{k}})\\
\hat{J}_k-J(\pi_{\theta_{k}},\hat{\mu}_{k})
\end{array}\right],\left[\begin{array}{c}
\Pi_{\Ecal_\perp}\bar{G}^V(\theta_k,\hat{V}_k,\hat{J}_k,\hat{\mu}_k)\\
\bar{G}^J(\theta_k,\hat{J}_k,\hat{\mu}_k)
\end{array}\right]\right\rangle\notag\\
&\leq \|\Pi_{\Ecal_\perp}(\hat{V}_{k}-V^{\pi_{\theta_{k}},\,\hat{\mu}_{k}})\|^2+(\hat{J}_k-J(\pi_{\theta_{k}},\hat{\mu}_{k}))^2+\beta_k^2\|\bar{G}(\theta_k,\hat{V}_k,\hat{J}_k,\hat{\mu}_k)\|^2\notag\\
&\hspace{20pt}-\gamma \beta_k \|\Pi_{\Ecal_\perp}(\hat{V}_{k}-V^{\pi_{\theta_{k}},\,\hat{\mu}_{k}})\|^2 -\gamma \beta_k (\hat{J}_k-J(\pi_{\theta_{k}},\hat{\mu}_{k}))^2\notag\\
&= (1-\gamma \beta_k)\|\Pi_{\Ecal_\perp}(\hat{V}_{k}-V^{\pi_{\theta_{k}},\,\hat{\mu}_{k}})\|^2+(1-\gamma \beta_k)(\hat{J}_k-J(\pi_{\theta_{k}},\hat{\mu}_{k}))^2\notag\\
&\hspace{20pt}+\beta_k^2\|\bar{G}(\theta_k,\hat{V}_k,\hat{J}_k,\hat{\mu}_k)-\bar{G}(\theta_k,V^{\pi_{\theta_{k}},\,\hat{\mu}_{k}},J(\pi_{\theta_{k}},\,\hat{\mu}_{k}),\hat{\mu}_k)\|^2\notag\\
&\leq (1-\gamma \beta_k) \|\Pi_{\Ecal_\perp}(\hat{V}_{k}-V^{\pi_{\theta_{k}},\,\hat{\mu}_{k}})\|^2 +(1 -\gamma \beta_k) (\hat{J}_k-J(\pi_{\theta_{k}},\hat{\mu}_{k}))^2\notag\\
&\hspace{20pt}+L_G^2\beta_k^2\left(\|\Pi_{\Ecal_\perp}(\hat{V}_{k}-V^{\pi_{\theta_{k}},\,\hat{\mu}_{k}})\|+|\hat{J}_{k}-J(\pi_{\theta_{k}},\,\hat{\mu}_{k})|\right)^2\notag\\
&\leq (1-\gamma \beta_k + 2L_G^2\beta_k^2) \|\Pi_{\Ecal_\perp}(\hat{V}_{k}-V^{\pi_{\theta_{k}},\,\hat{\mu}_{k}})\|^2 +(1 -\gamma \beta_k + 2L_G^2\beta_k^2) (\hat{J}_k-J(\pi_{\theta_{k}},\hat{\mu}_{k}))^2\notag\\
&\leq (1-\frac{\gamma \beta_k}{2}) (\varepsilon_k^{V}+\varepsilon_k^{J}),\label{prop:V_conv:proof_eq2}
\end{align}
where the second inequality applies Lemma~\ref{lem:G_strongmonotone}, the first equation uses the $\bar{G}(\theta,V^{\pi_{\theta},\,\mu},J(\pi_{\theta},\mu),\mu)=0$ for any $\theta,\mu$, third inequality follows from the Lipschitz continuity of operator $\bar{G}$ established in Lemma~\ref{lem:Lipschitz_operators}, and the final inequality follows from the step size condition $\beta_k\leq\frac{\gamma}{4L_G^2}$.

To treat the second and third term of \eqref{prop:V_conv:proof_eq1}, we use the boundedness of $\Delta g_k$ and the Lipschitz continuity conditions from Lemma~\ref{lem:Lipschitz_V}
\begin{align}
&\beta_k^2\|\Delta g_k\|^2+\left\|\left[\begin{array}{c}
\Pi_{\Ecal_\perp}(V^{\pi_{\theta_{k+1}},\,\hat{\mu}_{k+1}}-V^{\pi_{\theta_{k}},\,\hat{\mu}_{k}})\\
J(\pi_{\theta_{k+1}},\hat{\mu}_{k+1})-J(\pi_{\theta_{k}},\hat{\mu}_{k})
\end{array}\right]\right\|^2\notag\\
&\leq \beta_k^2\|\Delta g_k\|^2+\|\Pi_{\Ecal_\perp}(V^{\pi_{\theta_{k+1}},\,\hat{\mu}_{k+1}}-V^{\pi_{\theta_{k}},\,\hat{\mu}_{k}})\|^2+(J(\pi_{\theta_{k+1}},\hat{\mu}_{k+1})-J(\pi_{\theta_{k}},\hat{\mu}_{k}))^2\notag\\
&\leq 4B_G^2 \beta_k^2 + \left(L_V\|\theta_{k+1}-\theta_{k}\|+L_V\|\hat{\mu}_{k+1}-\hat{\mu}_{k}\|\right)^2+\left(L_V\|\theta_{k+1}-\theta_{k}\|+L_V\|\hat{\mu}_{k+1}-\hat{\mu}_{k}\|\right)^2\notag\\
&= 4B_G^2 \beta_k^2 + 2L_V^2\left(\alpha_k\|f_k\|+\xi_k\|h_k\|\right)^2\notag\\
&= 4B_G^2 \beta_k^2 + 2L_V^2\xi_k^2\left(B_F+B_H\right)^2\notag\\
&\leq 4L_V^2(B_F^2+B_G^2+B_H^2)\beta_k^2,
\end{align}
where we combine terms using the step size condition $\alpha_k\leq\xi_k\leq\beta_k$.

The fourth term of \eqref{prop:V_conv:proof_eq1} can be bounded leveraging the result in \eqref{prop:V_conv:proof_eq2} as follows
\begin{align}
&2\beta_k\left\langle \left[\begin{array}{c}\Pi_{\Ecal_\perp}(\hat{V}_{k}-V^{\pi_{\theta_{k}},\,\hat{\mu}_{k}})\\
\hat{J}_k-J(\pi_{\theta_{k}},\hat{\mu}_{k})
\end{array}\right]+\beta_k \left[\begin{array}{c}
\Pi_{\Ecal_\perp}\bar{G}^V(\theta_k,\hat{V}_k,\hat{J}_k,\hat{\mu}_k)\\
\bar{G}^J(\theta_k,\hat{J}_k,\hat{\mu}_k)
\end{array}\right],\left[\begin{array}{c}\Pi_{\Ecal_\perp}\Delta g_k^V\\
\Delta g_k^J\end{array}
\right]\right\rangle\notag\\
&\leq\frac{\gamma\beta_k}{8}\left\|\left[\begin{array}{c}\Pi_{\Ecal_\perp}(\hat{V}_{k}-V^{\pi_{\theta_{k}},\,\hat{\mu}_{k}})\\
\hat{J}_k-J(\pi_{\theta_{k}},\hat{\mu}_{k})
\end{array}\right]+\beta_k \left[\begin{array}{c}
\Pi_{\Ecal_\perp}\bar{G}^V(\theta_k,\hat{V}_k,\hat{J}_k,\hat{\mu}_k)\\
\bar{G}^J(\theta_k,\hat{J}_k,\hat{\mu}_k)
\end{array}\right]\right\|^2\hspace{-2pt}+\hspace{-2pt}\frac{8\beta_k}{\gamma}\left\|\left[\begin{array}{c}\Pi_{\Ecal_\perp}\Delta g_k^V\\
g_k^J\end{array}
\right]\right\|^2\notag\\
&\leq \frac{\gamma\beta_k}{8}(1-\frac{\gamma \beta_k}{2}) (\varepsilon_k^{V}+\varepsilon_k^{J})+\frac{8\beta_k}{\gamma}\|\Delta g_k\|^2\notag\\
&\leq \frac{\gamma\beta_k}{8}(\varepsilon_k^{V}+\varepsilon_k^{J})+\frac{8\beta_k}{\gamma}\|\Delta g_k\|^2.
\end{align}

Similarly, for the fifth term of \eqref{prop:V_conv:proof_eq1}, we have
\begin{align}
&2\left\langle\left[\begin{array}{c}\Pi_{\Ecal_\perp}(\hat{V}_{k}-V^{\pi_{\theta_{k}},\,\hat{\mu}_{k}})\\
\hat{J}_k-J(\pi_{\theta_{k}},\hat{\mu}_{k})
\end{array}\right]+\beta_k \left[\begin{array}{c}
\Pi_{\Ecal_\perp}\bar{G}^V(\theta_k,\hat{V}_k,\hat{J}_k,\hat{\mu}_k)\\
\bar{G}^J(\theta_k,\hat{J}_k,\hat{\mu}_k)
\end{array}\right],\left[\begin{array}{c}
\Pi_{\Ecal_\perp}(V^{\pi_{\theta_{k+1}},\,\hat{\mu}_{k+1}}-V^{\pi_{\theta_{k}},\,\hat{\mu}_{k}})\\
J(\pi_{\theta_{k+1}},\hat{\mu}_{k+1})-J(\pi_{\theta_{k}},\hat{\mu}_{k})
\end{array}\right]\right\rangle\notag\\
&\leq \frac{\gamma\beta_k}{8}\left\|\left[\begin{array}{c}\Pi_{\Ecal_\perp}(\hat{V}_{k}-V^{\pi_{\theta_{k}},\,\hat{\mu}_{k}})\\
\hat{J}_k-J(\pi_{\theta_{k}},\hat{\mu}_{k})
\end{array}\right]+\beta_k \left[\begin{array}{c}
\Pi_{\Ecal_\perp}\bar{G}^V(\theta_k,\hat{V}_k,\hat{J}_k,\hat{\mu}_k)\\
\bar{G}^J(\theta_k,\hat{J}_k,\hat{\mu}_k)
\end{array}\right]\right\|^2\notag\\
&\hspace{20pt}+\frac{8}{\gamma\beta_k}\left\|\left[\begin{array}{c}
\Pi_{\Ecal_\perp}(V^{\pi_{\theta_{k+1}},\,\hat{\mu}_{k+1}}-V^{\pi_{\theta_{k}},\,\hat{\mu}_{k}})\\
J(\pi_{\theta_{k+1}},\hat{\mu}_{k+1})-J(\pi_{\theta_{k}},\hat{\mu}_{k})
\end{array}\right]\right\|^2\notag\\
&\leq \frac{\gamma\beta_k}{8}(\varepsilon_k^{V}+\varepsilon_k^{J})+\frac{8}{\gamma\beta_k}\|V^{\pi_{\theta_{k+1}},\,\hat{\mu}_{k+1}}-V^{\pi_{\theta_{k}},\,\hat{\mu}_{k}}\|^2\notag\\
&\hspace{20pt}+\frac{8}{\gamma\beta_k}(J(\pi_{\theta_{k+1}},\hat{\mu}_{k+1})-J(\pi_{\theta_{k}},\,\hat{\mu}_{k}))^2\notag\\
&\leq \frac{\gamma\beta_k}{8}(\varepsilon_k^{V}+\varepsilon_k^{J})+\frac{16L_V^2}{\gamma\beta_k}\left(\|\pi_{\theta_{k+1}}-\pi_{\theta_{k}}\|^2+\|\hat{\mu}_{k+1}-\hat{\mu}_{k}\|^2\right)\notag\\
&\hspace{20pt}+\frac{16L_V^2}{\gamma\beta_k}\left(\|\pi_{\theta_{k+1}}-\pi_{\theta_{k}}\|^2+\|\hat{\mu}_{k+1}-\hat{\mu}_{k}\|^2\right)\notag\\
&\leq \frac{\gamma\beta_k}{8}(\varepsilon_k^{V}+\varepsilon_k^{J})+\frac{32L_V^2}{\gamma\beta_k}(\alpha_k^2\|f_k\|^2+\xi_k^2\|h_k\|^2)\notag\\
&\leq \frac{\gamma\beta_k}{8}(\varepsilon_k^{V}+\varepsilon_k^{J})\notag\\
&\hspace{20pt}+\frac{32L_V^2}{\gamma\beta_k} \left(4\alpha_k^2(\|\Delta f_k\|^2+L_F^2\varepsilon_k^{V}+L_F^2(L_V+1)^2\varepsilon_k^{\mu}+\varepsilon_k^{\pi})+2\xi_k^2(\|\Delta h_k\|^2+L_H^2\varepsilon_k^{\mu})\right)\notag\\
&\leq \frac{\gamma\beta_k}{8}(\varepsilon_k^{V}+\varepsilon_k^{J})+\frac{128L_V^2\alpha_k^2}{\gamma\beta_k} \left(\|\Delta f_k\|^2+L_F^2\varepsilon_k^{V}+4L_F^2 L_V^2\varepsilon_k^{\mu}+\varepsilon_k^{\pi}\right)\notag\\
&\hspace{20pt}+\frac{64L_V^2\xi_k^2}{\gamma\beta_k}\left(\|\Delta h_k\|^2+L_H^2\varepsilon_k^{\mu}\right),
\end{align}
where the third inequality applies Lemma~\ref{lem:Lipschitz_V} and the fifth inequality applies Lemma~\ref{lem:bound_f_g}.

The final term of \eqref{prop:V_conv:proof_eq1} can be bounded simply with the Cauchy-Schwarz inequality
\begin{align}
&2\beta_k\left\langle\Delta g_k, \left[\begin{array}{c}
\Pi_{\Ecal_\perp}(V^{\pi_{\theta_{k+1}},\,\hat{\mu}_{k+1}}-V^{\pi_{\theta_{k}},\,\hat{\mu}_{k}})\\
J(\pi_{\theta_{k+1}},\hat{\mu}_{k+1})-J(\pi_{\theta_{k}},\hat{\mu}_{k})
\end{array}\right] \right\rangle\notag\\
&\leq 2\beta_k\left\|\Delta g_k\right\| \left\|\left[\begin{array}{c}
\Pi_{\Ecal_\perp}(V^{\pi_{\theta_{k+1}},\,\hat{\mu}_{k+1}}-V^{\pi_{\theta_{k}},\,\hat{\mu}_{k}})\\
J(\pi_{\theta_{k+1}},\hat{\mu}_{k+1})-J(\pi_{\theta_{k}},\hat{\mu}_{k})
\end{array}\right]\right\|\notag\\
&\leq 2\beta_k\left\|\Delta g_k\right\|\left(\|\Pi_{\Ecal_\perp}(V^{\pi_{\theta_{k+1}},\,\hat{\mu}_{k+1}}-V^{\pi_{\theta_{k}},\,\hat{\mu}_{k}})\|+|J(\pi_{\theta_{k+1}},\hat{\mu}_{k+1})-J(\pi_{\theta_{k}},\,\hat{\mu}_{k})|\right)\notag\\
&\leq 4B_G\beta_k\cdot \left(L_V(\|\pi_{\theta_{k+1}}-\pi_{\theta_{k}}\|+\|\hat{\mu}_{k+1}-\hat{\mu}_{k}\|)+L_V(\|\pi_{\theta_{k+1}}-\pi_{\theta_{k}}\|+\|\hat{\mu}_{k+1}-\hat{\mu}_{k}\|)\right)\notag\\
&\leq 8L_V B_G \beta_k(B_F\alpha_k+B_H\xi_k)\notag\\
&\leq 16L_V B_F B_G B_H \beta_k\xi_k.
\label{prop:V_conv:proof_eq3}
\end{align}

Plugging \eqref{prop:V_conv:proof_eq2}-\eqref{prop:V_conv:proof_eq3} into \eqref{prop:V_conv:proof_eq1}, we get
\begin{align*}
&\varepsilon_{k+1}^{V}+\varepsilon_{k+1}^{J}\notag\\
&\leq(1-\frac{\gamma \beta_k}{2}) (\varepsilon_k^{V}+\varepsilon_k^{J})+ 4L_V^2(B_F^2+B_G^2+B_H^2)\beta_k^2+\frac{\gamma\beta_k}{8}(\varepsilon_k^{V}+\varepsilon_k^{J})+\frac{8\beta_k}{\gamma}\|\Delta g_k\|^2\notag\\
&\hspace{20pt}+\frac{\gamma\beta_k}{8}(\varepsilon_k^{V}+\varepsilon_k^{J})+\frac{128L_V^2\alpha_k^2}{\gamma\beta_k} \left(\|\Delta f_k\|^2+L_F^2\varepsilon_k^{V}+4L_F^2 L_V^2\varepsilon_k^{\mu}+\varepsilon_k^{\pi}\right)\notag\\
&\hspace{20pt}+\frac{64L_V^2\xi_k^2}{\gamma\beta_k}\left(\|\Delta h_k\|^2+L_H^2\varepsilon_k^{\mu}\right)+16L_V B_F B_G B_H \beta_k\xi_k\notag\\
&\leq (1-\frac{\gamma \beta_k}{4}) (\varepsilon_k^{V}+\varepsilon_k^{J})+\frac{128L_V^2\alpha_k^2}{\gamma\beta_k} \|\Delta f_k\|^2+\frac{8\beta_k}{\gamma}\|\Delta g_k\|^2+\frac{64L_V^2\xi_k^2}{\gamma\beta_k} \|\Delta h_k\|^2\notag\\
&\hspace{20pt}+\frac{128L_V^2\alpha_k^2}{\gamma\beta_k} (L_F^2 \varepsilon_k^{V}+\varepsilon_k^{\pi})+\frac{192L_V^2\xi_k^2}{\gamma\beta_k} \varepsilon_k^{\mu}+28L_V^2 B_F^2 B_G^2 B_H^2 \beta_k^2,
\end{align*}
where we use the conditions $\xi_k\leq\beta_k$ and $\alpha_k\leq\frac{L_H}{2L_F L_V}\xi_k$ in the last inequality to simplify and combine terms.

\qed

\subsection{Proof of Proposition~\ref{prop:g_conv}}\label{sec:g_conv:proof}

The proof of Proposition~\ref{prop:g_conv} relies on the following lemma, the proof of which is presented in Sec.\ref{sec:g_markovian:proof}.
\begin{lem}\label{lem:g_markovian}
We have for all $k\geq\tau_k$
\begin{align*}
&\mathbb{E}[\langle\Delta g_k, G(\theta_{k},\hat{V}_{k},\hat{J}_{k},\hat{\mu}_{k},s_k,a_k,s_{k+1})-\bar{G}(\theta_{k},\hat{V}_{k},\hat{J}_{k},\hat{\mu}_{k})\rangle]\notag\\
&\hspace{150pt}\leq (22L+2|\Acal|) L_F L_{TV} B_F B_G^2 B_H \tau_k^2\lambda_{k-\tau_k}.
\end{align*}
\end{lem}

By the update rule of $f_k$, 
\begin{align*}
\Delta g_{k+1}
&=g_{k+1}-\bar{G}(\theta_{k+1},\hat{V}_{k+1},\hat{J}_{k+1},\hat{\mu}_{k+1})\notag\\
&=(1-\lambda_k)g_k+\lambda_k \bar{G}(\theta_{k},\hat{V}_{k},\hat{J}_{k},\hat{\mu}_{k},s_k,a_k,s_{k+1})-\bar{G}(\theta_{k+1},\hat{V}_{k+1},\hat{J}_{k+1},\hat{\mu}_{k+1})\notag\\
&=(1-\lambda_k)g_k+\lambda_k \bar{G}(\theta_{k},\hat{V}_{k},\hat{J}_{k},\hat{\mu}_{k})-\bar{G}(\theta_{k+1},\hat{V}_{k+1},\hat{J}_{k+1},\hat{\mu}_{k+1})\notag\\
&\hspace{20pt}+\lambda_k\Big( G(\theta_{k},\hat{V}_{k},\hat{J}_{k},\hat{\mu}_{k},s_k,a_k,s_{k+1})-\bar{G}(\theta_{k},\hat{V}_{k},\hat{J}_{k},\hat{\mu}_{k})\Big)\notag\\
&=(1-\lambda_k)\Delta g_k+ \Big(\bar{G}(\theta_{k},\hat{V}_{k},\hat{J}_{k},\hat{\mu}_{k})-\bar{G}(\theta_{k+1},\hat{V}_{k+1},\hat{J}_{k+1},\hat{\mu}_{k+1})\Big)\notag\\
&\hspace{20pt}+\lambda_k\Big( G(\theta_{k},\hat{V}_{k},\hat{J}_{k},\hat{\mu}_{k},s_k,a_k,s_{k+1})-\bar{G}(\theta_{k},\hat{V}_{k},\hat{J}_{k},\hat{\mu}_{k})\Big).
\end{align*}

Taking the norm, we have
\begin{align}
&\|\Delta g_{k+1}\|^2\notag\\
&=(1-\lambda_k)^2\|\Delta G_{k}\|^2+\|\bar{G}(\theta_{k},\hat{V}_{k},\hat{J}_{k},\hat{\mu}_{k})-\bar{G}(\theta_{k+1},\hat{V}_{k+1},\hat{J}_{k+1},\hat{\mu}_{k+1})\|^2\notag\\
&\hspace{20pt}+\lambda_k^2\|G(\theta_{k},\hat{V}_{k},\hat{J}_{k},\hat{\mu}_{k},s_k,a_k,s_{k+1})-\bar{G}(\theta_{k},\hat{V}_{k},\hat{J}_{k},\hat{\mu}_{k})\|^2\notag\\
&\hspace{20pt}+(1-\lambda_k)\langle\Delta g_k,\bar{G}(\theta_{k},\hat{V}_{k},\hat{J}_{k},\hat{\mu}_{k})-\bar{G}(\theta_{k+1},\hat{V}_{k+1},\hat{J}_{k+1},\hat{\mu}_{k+1})\rangle\notag\\
&\hspace{20pt}+(1-\lambda_k)\lambda_k\langle\Delta g_k,G(\theta_{k},\hat{V}_{k},\hat{J}_{k},\hat{\mu}_{k},s_k,a_k,s_{k+1})-\bar{G}(\theta_{k},\hat{V}_{k},\hat{J}_{k},\hat{\mu}_{k})\rangle\notag\\
&\hspace{20pt}+\lambda_k\langle \bar{G}(\theta_{k},\hat{V}_{k},\hat{J}_{k},\hat{\mu}_{k})-\bar{G}(\theta_{k+1},\hat{V}_{k+1},\hat{J}_{k+1},\hat{\mu}_{k+1}), G(\theta_{k},\hat{V}_{k},\hat{J}_{k},\hat{\mu}_{k},s_k,a_k,s_{k+1})-\bar{G}(\theta_{k},\hat{V}_{k},\hat{J}_{k},\hat{\mu}_{k})\rangle\notag\\
&\leq (1-\lambda_k)^2\|\Delta g_{k}\|^2+2\|\bar{G}(\theta_{k},\hat{V}_{k},\hat{J}_{k},\hat{\mu}_{k})-\bar{G}(\theta_{k+1},\hat{V}_{k+1},\hat{J}_{k+1},\hat{\mu}_{k+1})\|^2\notag\\
&\hspace{20pt}+2\lambda_k^2\|G(\theta_{k},\hat{V}_{k},\hat{J}_{k},\hat{\mu}_{k},s_k,a_k,s_{k+1})-\bar{G}(\theta_{k},\hat{V}_{k},\hat{J}_{k},\hat{\mu}_{k})\|^2\notag\\
&\hspace{20pt}+\frac{\lambda_k}{2}\|\Delta g_k\|^2+\frac{2}{\lambda_k}\|\bar{G}(\theta_{k},\hat{V}_{k},\hat{J}_{k},\hat{\mu}_{k})-\bar{G}(\theta_{k+1},\hat{V}_{k+1},\hat{J}_{k+1},\hat{\mu}_{k+1})\|^2\notag\\
&\hspace{20pt}+(1-\lambda_k)\lambda_k\langle\Delta g_k,G(\theta_{k},\hat{V}_{k},\hat{J}_{k},\hat{\mu}_{k},s_k,a_k,s_{k+1})-\bar{G}(\theta_{k},\hat{V}_{k},\hat{J}_{k},\hat{\mu}_{k})\rangle\notag\\
&\leq (1-\lambda_k)\|\Delta g_{k}\|^2+(-\frac{\lambda_k}{2}+\lambda_k^2)\|\Delta g_{k}\|^2\notag\\
&\hspace{20pt}+\frac{4}{\lambda_k}\|\bar{G}(\theta_{k},\hat{V}_{k},\hat{J}_{k},\hat{\mu}_{k})-\bar{G}(\theta_{k+1},\hat{V}_{k+1},\hat{J}_{k+1},\hat{\mu}_{k+1})\|^2\notag\\
&\hspace{20pt}+8B_G^2\lambda_k^2+(1-\lambda_k)\lambda_k\langle\Delta g_k,G(\theta_{k},\hat{V}_{k},\hat{J}_{k},\hat{\mu}_{k},s_k,a_k,s_{k+1})-\bar{G}(\theta_{k},\hat{V}_{k},\hat{J}_{k},\hat{\mu}_{k})\rangle,\label{prop:g_conv:proof_eq1}
\end{align}
where the final inequality follows from the step size condition $\lambda_k\leq 1$ and the boundedness of operator $F$.

Taking expectation and plugging in the result of Lemma~\ref{lem:f_markovian}, we can simplify \eqref{prop:g_conv:proof_eq1} as
\begin{align}
&\mathbb{E}[\|\Delta g_{k+1}\|^2]\notag\\
&\leq \mathbb{E}\Big[(1-\lambda_k)\|\Delta g_{k}\|^2+(-\frac{\lambda_k}{2}+\lambda_k^2)\|\Delta g_{k}\|^2\notag\\
&\hspace{20pt}+\frac{4}{\lambda_k}\|\bar{G}(\theta_{k},\hat{V}_{k},\hat{J}_{k},\hat{\mu}_{k})-\bar{G}(\theta_{k+1},\hat{V}_{k+1},\hat{J}_{k+1},\hat{\mu}_{k+1})\|^2\notag\\
&\hspace{20pt}+8B_G^2\lambda_k^2+(1-\lambda_k)\lambda_k\langle\Delta g_k, G(\theta_{k},\hat{V}_{k},\hat{J}_{k},\hat{\mu}_{k},s_k,a_k,s_{k+1})-\bar{G}(\theta_{k},\hat{V}_{k},\hat{J}_{k},\hat{\mu}_{k})\rangle\Big]\notag\\
&\leq (1-\lambda_k)\mathbb{E}[\|\Delta g_{k}\|^2]+(-\frac{\lambda_k}{2}+\lambda_k^2)\mathbb{E}[\|\Delta g_{k}\|^2]+8B_G^2\lambda_k^2\notag\\
&\hspace{20pt}+\frac{4L_G^2}{\lambda_k}\mathbb{E}\left[\left(\|\theta_{k}-\theta_{k+1}\|+\|\hat{V}_{k}-\hat{V}_{k+1}\|+|\hat{J}_{k}-\hat{J}_{k+1}|+\|\hat{\mu}_{k}-\hat{\mu}_{k+1}\|\right)^2\right]\notag\\
&\hspace{20pt}+(1-\lambda_k)\lambda_k\cdot(22L+2|\Acal|) L_F L_{TV} B_F B_G^2 B_H \tau_k^2\lambda_{k-\tau_k}\notag\\
&\leq (1-\lambda_k)\mathbb{E}[\|\Delta g_{k}\|^2]+(-\frac{\lambda_k}{2}+\lambda_k^2)\mathbb{E}[\|\Delta g_{k}\|^2]+(30L+2|\Acal|) L_F L_{TV} B_F B_G^2 B_H \tau_k^2\lambda_k\lambda_{k-\tau_k}\notag\\
&\hspace{20pt}+\frac{4L_G^2}{\lambda_k}\mathbb{E}[\left(\alpha_k\|f_k\|+\beta_k\|g_k^V\|+\beta_k|g_k^J|+\xi_k\|h_k\|\right)^2]\notag\\
&\leq (1-\lambda_k)\mathbb{E}[\|\Delta g_{k}\|^2]+(-\frac{\lambda_k}{2}+\lambda_k^2)\mathbb{E}[\|\Delta g_{k}\|^2]+(30L+2|\Acal|) L_F L_{TV} B_F B_G^2 B_H \tau_k^2\lambda_k\lambda_{k-\tau_k}\notag\\
&\hspace{20pt}+\frac{4L_G^2}{\lambda_k}\mathbb{E}[\left(\alpha_k\|f_k\|+\sqrt{|\Scal|+1}\beta_k\|g_k\|+\xi_k\|h_k\|\right)^2]\notag\\
&\leq (1-\lambda_k)\mathbb{E}[\|\Delta g_{k}\|^2]+(-\frac{\lambda_k}{2}+\lambda_k^2)\mathbb{E}[\|\Delta g_{k}\|^2]+(30L+2|\Acal|) L_F L_{TV} B_F B_G^2 B_H \tau_k^2\lambda_k\lambda_{k-\tau_k}\notag\\
&\hspace{20pt}+\frac{12L_G^2 \alpha_k^2}{\lambda_k}\mathbb{E}[\left(\|\Delta f_k\|+L_F\sqrt{\varepsilon_k^{V}}+L_F(L_V+1)\sqrt{\varepsilon_k^{\mu}}+\sqrt{\varepsilon_k^{\pi}}\right)^2]\notag\\
&\hspace{20pt}+\frac{24|\Scal|L_G^2 \beta_k^2}{\lambda_k}\mathbb{E}[\left(\|\Delta g_k\|+L_G\sqrt{\varepsilon_k^{V}}+L_G\sqrt{\varepsilon_k^{J}}\right)^2]+\frac{12L_G^2 \xi_k^2}{\lambda_k}\mathbb{E}[\left(\|\Delta h_k\|+L_H\sqrt{\varepsilon_k^{\mu}}\right)^2],\label{prop:g_conv:proof_eq2}
\end{align}
where the fourth inequality follows from $\|g_{k}^V\|+|g_{k}^J|\leq \|g_{k}^V\|_1+|g_{k}^J|=\|g_k\|_1\leq\sqrt{|\Scal|+1}\|g_k\|$.

We can simplify the sum of the last three terms as follows
\begin{align}
&\frac{12L_G^2 \alpha_k^2}{\lambda_k}\mathbb{E}\left[\left(\|\Delta f_k\|+L_F\sqrt{\varepsilon_k^{V}}+L_F(L_V+1)\sqrt{\varepsilon_k^{\mu}}+\sqrt{\varepsilon_k^{\pi}}\right)^2\right]\notag\\
&+\frac{24|\Scal|L_G^2 \beta_k^2}{\lambda_k}\mathbb{E}[\left(\|\Delta g_k\|+L_G\sqrt{\varepsilon_k^{V}}+L_G\sqrt{\varepsilon_k^{J}}\right)^2]+\frac{12L_G^2 \xi_k^2}{\lambda_k}\mathbb{E}[\left(\|\Delta h_k\|+L_H\sqrt{\varepsilon_k^{\mu}}\right)^2]\notag\\
&\leq \frac{48L_G^2 \alpha_k^2}{\lambda_k}\mathbb{E}[\|\Delta f_k\|^2+L_F^2\varepsilon_k^{V}+4L_F^2 L_V^2 \varepsilon_k^{\mu}+\varepsilon_k^{\pi}] + \frac{72|\Scal|L_G^2 \beta_k^2}{\lambda_k}\mathbb{E}[\|\Delta g_k\|^2+L_G^2\varepsilon_k^{V}+L_G^2\varepsilon_k^{J}]\notag\\
&\hspace{20pt}+\frac{24L_G^2 \xi_k^2}{\lambda_k}\mathbb{E}[\|\Delta h_k\|^2+L_H^2\varepsilon_k^{\mu}]\notag\\
&\leq \mathbb{E}\Big[\frac{48L_G^2 \alpha_k^2}{\lambda_k}\|\Delta f_k\|^2+\frac{72|\Scal|L_G^2 \beta_k^2}{\lambda_k}\|\Delta g_k\|^2+\frac{24L_G^2 \xi_k^2}{\lambda_k}\|\Delta h_k\|^2+\frac{48L_G^2 \alpha_k^2}{\lambda_k}\varepsilon_{k}^{\pi}\notag\\
&\hspace{20pt}+\frac{216L_F^2L_G^2L_H^2L_V^2\xi_k^2}{\lambda_k}\varepsilon_{k}^{\mu}+\frac{120|\Scal|L_F^2 L_G^4\beta_k^2}{\lambda_k}(\varepsilon_k^{V}+\varepsilon_k^{J})\Big].
\label{prop:g_conv:proof_eq3}
\end{align}

Combining \eqref{prop:g_conv:proof_eq2} and \eqref{prop:g_conv:proof_eq3}, we have
\begin{align*}
&\mathbb{E}[\|\Delta g_{k+1}\|^2]\notag\\
&\leq (1-\lambda_k)\mathbb{E}[\|\Delta g_{k}\|^2]+(-\frac{\lambda_k}{2}+\lambda_k^2)\mathbb{E}[\|\Delta g_{k}\|^2]+(30L+2|\Acal|) L_F L_{TV} B_F B_G^2 B_H \tau_k^2\lambda_k\lambda_{k-\tau_k}\notag\\
&\hspace{20pt}+\frac{12L_G^2 \alpha_k^2}{\lambda_k}\left(\|\Delta f_k\|+L_F\sqrt{\varepsilon_k^{V}}+L_F(L_V+1)\sqrt{\varepsilon_k^{\mu}}+\sqrt{\varepsilon_k^{\pi}}\right)^2\notag\\
&\hspace{20pt}+\frac{24|\Scal|L_G^2 \beta_k^2}{\lambda_k}\left(\|\Delta g_k\|+L_G\sqrt{\varepsilon_k^{V}}+L_G\sqrt{\varepsilon_k^{J}}\right)^2+\frac{12L_G^2 \xi_k^2}{\lambda_k}\left(\|\Delta h_k\|+L_H\sqrt{\varepsilon_k^{\mu}}\right)^2\notag\\
&\leq(1-\lambda_k)\mathbb{E}[\|\Delta g_{k}\|^2]+(-\frac{\lambda_k}{2}+\lambda_k^2)\mathbb{E}[\|\Delta g_{k}\|^2]+(30L+2|\Acal|) L_F L_{TV} B_F B_G^2 B_H \tau_k^2\lambda_k\lambda_{k-\tau_k}\notag\\
&\hspace{20pt}+\frac{48L_G^2 \alpha_k^2}{\lambda_k}\mathbb{E}[\|\Delta f_k\|^2]+\frac{72|\Scal|L_G^2 \beta_k^2}{\lambda_k}\mathbb{E}[\|\Delta g_k\|^2]+\frac{24L_G^2 \xi_k^2}{\lambda_k}\mathbb{E}[\|\Delta h_k\|^2]+\frac{48L_G^2 \alpha_k^2}{\lambda_k}\mathbb{E}[\varepsilon_{k}^{\pi}]\notag\\
&\hspace{20pt}+\frac{216L_F^2L_G^2L_H^2L_V^2\xi_k^2}{\lambda_k}\mathbb{E}[\varepsilon_{k}^{\mu}]+\frac{120|\Scal|L_F^2 L_G^4\beta_k^2}{\lambda_k}\mathbb{E}[\varepsilon_k^{V}+\varepsilon_k^{J}]\notag\\
&\leq(1-\lambda_k)\mathbb{E}[\|\Delta g_{k}\|^2]+(-\frac{\lambda_k}{2}+\lambda_k^2+\frac{72|\Scal|L_G^2 \beta_k^2}{\lambda_k})\mathbb{E}[\|\Delta g_{k}\|^2]+\frac{48L_G^2 \alpha_k^2}{\lambda_k}\mathbb{E}[\|\Delta f_k\|^2]\notag\\
&\hspace{20pt}+\frac{24L_G^2 \xi_k^2}{\lambda_k}\mathbb{E}[\|\Delta h_k\|^2]+\frac{48L_G^2 \alpha_k^2}{\lambda_k}\mathbb{E}[\varepsilon_{k}^{\pi}]+\frac{216L_F^2L_G^2L_H^2L_V^2\xi_k^2}{\lambda_k}\mathbb{E}[\varepsilon_{k}^{\mu}]\notag\\
&\hspace{20pt}+\frac{120|\Scal|L_F^2 L_G^4\beta_k^2}{\lambda_k}\mathbb{E}[\varepsilon_k^{V}+\varepsilon_k^{J}]+(30L+2|\Acal|) L_F L_{TV} B_F B_G^2 B_H \tau_k^2\lambda_k\lambda_{k-\tau_k}.
\end{align*}

\qed

\section{Proof of Lemmas}

\subsection{Proof of Lemma~\ref{lem:Lipschitz_V}}

The Lipschitz continuity conditions of the value function and $J$ function in the policy are proved in Lemma 3 and Lemma 2 of \citet{kumar2024global}, respectively. The Lipschitz continuity in the mean field can be proved using the same line of argument under Assumption~\ref{assump:Lipschitz}.

The Lipschitz gradient condition of $J$ in $\theta$ is proved in Lemma 4 of \citet{kumar2024global} and can be extended to the gradient of $J$ in $\mu$ by a similar argument.

\qed

\subsection{Proof of Lemma~\ref{lem:boundedness}}

First, by definition in \eqref{eq:def_FGH},
\begin{align*}
\|F(\theta,V,\mu,s,a,s')\|&= \|(r(s,a,\mu) + V(s'))\nabla_{\theta}\log\pi_{\theta}(a\mid s)\|\notag\\
&\leq (|r(s,a,\mu)|+|V(s')|)\|\nabla_{\theta}\log\pi_{\theta}(a\mid s)\|\notag\\
&\leq (1+B_V)\cdot 1\notag\\
&\leq B_V+1,
\end{align*}
where the second inequality is due to the softmax function being Lipschitz with constant $1$.

Similarly, we have
\begin{align*}
\|G^V(V,J,\mu,s,a,s')\|&=\|(r(s,a,\mu) -J + V(s') - V(s))e_{s}\|\notag\\
&\leq (|r(s,a,\mu)|+|J|+|V(s')|-|V(s)|)\|e_{s}\|\notag\\
&\leq (1+1+B_V+B_V)\cdot 1\notag]\\
&\leq 2B_V+2,
\end{align*}
and
\begin{align*}
|G^J(J,\mu,s,a)|=|c_J (r(s,a,\mu) - J)|\leq 2c_J,
\end{align*}
which implies
\begin{align*}
\|G(V,J,\mu,s,a,s')\|\leq\|G^V(V,J,\mu,s,a,s')\|+|G^J(J,\mu,s,a)|\leq 2(B_V+c_J+2).
\end{align*}

Finally, we have
\begin{align*}
\|H(\mu,s)\|&=\|e_s-\mu\|\leq\|e_s\|+\|\mu\|\leq2.
\end{align*}

\qed

\subsection{Proof of Lemma~\ref{lem:Lipschitz_operators}}
By the definition of $\bar{F}(\theta,V,\mu)$ in \eqref{eq:def_FGH_aggregate},
\begin{align}
&\|\bar{F}(\theta_1,V_1,\mu_1)-\bar{F}(\theta_2,V_2,\mu_2)\|\notag\\
&=\|\mathbb{E}_{s\sim \nu^{\pi_{\theta_1},\,\mu_1},a\sim\pi_{\theta_1}(\cdot\mid s),s'\sim\Pcal^{\mu_1}(\cdot\mid s,a)}[F(\theta_1,V_1,\mu_1,s,a,s')]\notag\\
&\hspace{20pt}-\mathbb{E}_{s\sim \nu^{\pi_{\theta_2},\,\mu_2},a\sim\pi_{\theta_2}(\cdot\mid s),s'\sim\Pcal^{\mu_2}(\cdot\mid s,a)}[F(\theta_2,V_2,\mu_2,s,a,s')]\|\notag\\
&=\|\mathbb{E}_{s\sim \nu^{\pi_{\theta_1},\,\mu_1},a\sim\pi_{\theta_1}(\cdot\mid s),s'\sim\Pcal^{\mu_1}(\cdot\mid s,a)}[F(\theta_1,\Pi_{\Ecal_\perp}V_1,\mu_1,s,a,s')]\notag\\
&\hspace{20pt}-\mathbb{E}_{s\sim \nu^{\pi_{\theta_2},\,\mu_2},a\sim\pi_{\theta_2}(\cdot\mid s),s'\sim\Pcal^{\mu_2}(\cdot\mid s,a)}[F(\theta_2,\Pi_{\Ecal_\perp}V_2,\mu_2,s,a,s')]\|\notag\\
&= \|\sum_{s,a,s'}(\nu^{\pi_{\theta_1},\,\mu_1}(s)\pi_{\theta_1}(a\mid s)\Pcal^{\mu_1}(\cdot\mid s,a)-\nu^{\pi_{\theta_2},\,\mu_2}(s)\pi_{\theta_2}(a\mid s)\Pcal^{\mu_2}(\cdot\mid s,a))F(\theta_2,\Pi_{\Ecal_\perp}V_2,\mu_2,s,a,s')\notag\\
&\hspace{20pt}+\mathbb{E}_{s\sim \nu^{\pi_{\theta_1},\,\mu_1},a\sim\pi_{\theta_1}(\cdot\mid s),s'\sim\Pcal^{\mu_1}(\cdot\mid s,a)}[F(\theta_1,\Pi_{\Ecal_\perp}V_1,\mu_1,s,a,s')-F(\theta_2,\Pi_{\Ecal_\perp}V_2,\mu_2,s,a,s')]\|\notag\\
&\leq \|\mathbb{E}_{s\sim \nu^{\pi_{\theta_1},\,\mu_1},a\sim\pi_{\theta_1}(\cdot\mid s),s'\sim\Pcal^{\mu_1}(\cdot\mid s,a,\mu_1)}[F(\theta_1,\Pi_{\Ecal_\perp}V_1,\mu_1,s,a,s')-F(\theta_2,\Pi_{\Ecal_\perp}V_2,\mu_2,s,a,s')]\|\notag\\
&\hspace{20pt}+2B_F d_{TV}(\nu^{\pi_{\theta_1},\,\mu_1}\otimes\pi_{\theta_1}\otimes\Pcal^{\mu_1},\nu^{\pi_{\theta_2},\,\mu_2}\otimes\pi_{\theta_2}\otimes\Pcal^{\mu_2}),\label{lem:Lipschitz_operators:eq1}
\end{align}
where the inequality comes from the definition of TV distance in \eqref{eq:TV_def} and the second equation is a result of the fact that for any constant $c$
\begin{align*}
&\mathbb{E}_{s\sim \nu^{\pi_{\theta},\,\mu},a\sim\pi_{\theta}(\cdot\mid s),s'\sim\Pcal^{\mu}(\cdot\mid s,a)}[F(\theta,V+c\1_{|\Scal|},\mu,s,a,s')] \notag\\
&= \mathbb{E}_{s\sim \nu^{\pi_{\theta},\,\mu},a\sim\pi_{\theta}(\cdot\mid s),s'\sim\Pcal^{\mu}(\cdot\mid s,a)}[(r(s,a,\mu) + (V(s')+c)-(V(s)+c))\nabla_{\theta}\log\pi_{\theta}(a\mid s)]\notag\\
&=\mathbb{E}_{s\sim \nu^{\pi_{\theta},\,\mu},a\sim\pi_{\theta}(\cdot\mid s),s'\sim\Pcal^{\mu}(\cdot\mid s,a)}[(r(s,a,\mu) + V(s')-V(s))\nabla_{\theta}\log\pi_{\theta}(a\mid s)]\notag\\
&=\mathbb{E}_{s\sim \nu^{\pi_{\theta},\,\mu},a\sim\pi_{\theta}(\cdot\mid s),s'\sim\Pcal^{\mu}(\cdot\mid s,a)}[F(\theta,V,\mu,s,a,s')].
\end{align*}

For any $s,a,s'$ we have from \eqref{eq:def_FGH}
\begin{align}
&\|F(\theta_1,\Pi_{\Ecal_\perp}V_1,\mu_1,s,a,s')-F(\theta_2,\Pi_{\Ecal_\perp}V_2,\mu_2,s,a,s')\|\notag\\
&=\|(r(s,a,\mu_1) + \Pi_{\Ecal_\perp}V_1(s')-\Pi_{\Ecal_\perp}V_1(s))\nabla_{\theta}\log\pi_{\theta_1}(a\mid s)\notag\\
&\hspace{20pt}-(r(s,a,\mu_2) + \Pi_{\Ecal_\perp}V_2(s')-\Pi_{\Ecal_\perp}V_2(s))\nabla_{\theta}\log\pi_{\theta_2}(a\mid s)\|\notag\\
&\leq|r(s,a,\mu_1)-r(s,a,\mu_2)|\|\nabla_{\theta}\log\pi_{\theta_1}(a\mid s)\|\notag\\
&\hspace{20pt}+|r(s,a,\mu_2)|\|\nabla_{\theta}\log\pi_{\theta_1}(a\mid s)-\nabla_{\theta}\log\pi_{\theta_2}(a\mid s)\|\notag\\
&\hspace{20pt}+|\Pi_{\Ecal_\perp}V_1(s')-\Pi_{\Ecal_\perp}V_1(s)-\Pi_{\Ecal_\perp}V_2(s')+\Pi_{\Ecal_\perp}V_2(s)|\|\nabla_{\theta}\log\pi_{\theta_1}(a\mid s)\|\notag\\
&\hspace{20pt}+|\Pi_{\Ecal_\perp}V_2(s')-\Pi_{\Ecal_\perp}V_2(s)|\|\nabla_{\theta}\log\pi_{\theta_1}(a\mid s)-\nabla_{\theta}\log\pi_{\theta_2}(a\mid s)\|\notag\\
&\leq |r(s,a,\mu_1)-r(s,a,\mu_2)|+(1+2\|V\|)\|\nabla_{\theta}\log\pi_{\theta_1}(a\mid s)-\nabla_{\theta}\log\pi_{\theta_2}(a\mid s)\|\notag\\
&\leq L\|\mu_1-\mu_2\|+\|\nabla_{\theta}\log\pi_{\theta_1}(a\mid s)-\nabla_{\theta}\log\pi_{\theta_2}(a\mid s)\|\notag\\
&\hspace{20pt}+2\|\Pi_{\Ecal_\perp}V_1-\Pi_{\Ecal_\perp}V_2\|+2\|\Pi_{\Ecal_\perp}V_2\|\|\nabla_{\theta}\log\pi_{\theta_1}(a\mid s)-\nabla_{\theta}\log\pi_{\theta_2}(a\mid s)\|\notag\\
&\leq 5(2B_V+1)\|\theta_1-\theta_2\|+L\|\mu_1-\mu_2\|+2\|\Pi_{\Ecal_\perp}V_1-\Pi_{\Ecal_\perp}V_2\|,\label{lem:Lipschitz_operators:eq2}
\end{align}
where the second inequality bounds $\|\log\pi_{\theta_1}(a\mid s)\|$ by 1 due to the softmax function being Lipschitz with constant 1, the third inequality follows from Assumption~\ref{assump:Lipschitz}, and the final inequality is a result of the fact that the softmax function is smooth with constant 5 (see \citet{agarwal2021theory}[Lemma 52]).

Applying \eqref{lem:Lipschitz_operators:eq2} and the relationship in \eqref{eq:TV_stationary_Lipschitz} to \eqref{lem:Lipschitz_operators:eq1}, we have
\begin{align*}
&\|\bar{F}(\theta_1,V_1,\mu_1)-\bar{F}(\theta_2,V_2,\mu_2)\|\notag\\
&\leq \|\mathbb{E}_{s\sim \nu^{\pi_{\theta_1},\,\mu_1},a\sim\pi_{\theta_1}(\cdot\mid s),s'\sim\Pcal^{\mu_1}(\cdot\mid s,a,\mu_1)}[F(\theta_1,\Pi_{\Ecal_\perp}V_1,\mu_1,s,a,s')-F(\theta_2,\Pi_{\Ecal_\perp}V_2,\mu_2,s,a,s')]\|\notag\\
&\hspace{20pt}+2B_F d_{TV}(\nu^{\pi_{\theta_1},\,\mu_1}\otimes\pi_{\theta_1}\otimes\Pcal^{\mu_1},\nu^{\pi_{\theta_2},\,\mu_2}\otimes\pi_{\theta_2}\otimes\Pcal^{\mu_2})\notag\\
&\leq 5(2B_V+1)\|\theta_1-\theta_2\|+L\|\mu_1-\mu_2\|+2\|\Pi_{\Ecal_\perp}V_1-\Pi_{\Ecal_\perp}V_2\|\notag\\
&\hspace{20pt}+2B_F L_{TV}(\|\theta_1-\theta_2\|+\|\mu_1-\mu_2\|)\notag\\
&\leq (10B_V+L+2B_F L_{TV}+5)\left(\|\theta_1-\theta_2\|+\|\mu_1-\mu_2\|+\|\Pi_{\Ecal_\perp}V_1-\Pi_{\Ecal_\perp}V_2\|\right).
\end{align*}

Following a line of argument similar to \eqref{lem:Lipschitz_operators:eq1},
\begin{align}
&\|\bar{G}(\theta_1,V_1,J_1,\mu_1)-\bar{G}(\theta_2,V_2,J_2,\mu_2)\|\notag\\
&\leq \|\mathbb{E}_{s\sim \nu^{\pi_{\theta_1},\,\mu_1},a\sim\pi_{\theta_1}(\cdot\mid s),s'\sim\Pcal^{\mu_1}(\cdot\mid s,a,\mu_1)}[G(\Pi_{\Ecal_\perp}V_1,J_1,\mu_1,s,a,s')-G(\Pi_{\Ecal_\perp}V_2,J_2,\mu_2,s,a,s')]\|\notag\\
&\hspace{20pt}+2B_G d_{TV}(\nu^{\pi_{\theta_1},\,\mu_1}\otimes\pi_{\theta_1}\otimes\Pcal^{\mu_1},\nu^{\pi_{\theta_2},\,\mu_2}\otimes\pi_{\theta_2}\otimes\Pcal^{\mu_2}).\label{lem:Lipschitz_operators:eq3}
\end{align}

The first term of \eqref{lem:Lipschitz_operators:eq3} can be bounded in a manner similar to \eqref{lem:Lipschitz_operators:eq2}. For any $s,a,s'$, we have
\begin{align}
&\|G(\Pi_{\Ecal_\perp}V_1,J_1,\mu_1,s,a,s')-G(\Pi_{\Ecal_\perp}V_2,J_2,\mu_2,s,a,s')\|\notag\\
&\leq \|(r(s,a,\mu_1)-J_1+\Pi_{\Ecal_\perp}V_1(s')-\Pi_{\Ecal_\perp}V_1(s))e_s\notag\\
&\hspace{50pt}-(r(s,a,\mu_2)-J_2+\Pi_{\Ecal_\perp}V_2(s')-\Pi_{\Ecal_\perp}V_2(s))e_s\|\notag\\
&\hspace{20pt}+ c_J|r(s,a,\mu_1)-J_1-r(s,a,\mu_2)+J_2|\notag\\
&\leq |r(s,a,\mu_1)-r(s,a,\mu_2)|\|e_s\|+|J_1-J_2|\|e_s\|+2\|\Pi_{\Ecal_\perp}V_1-\Pi_{\Ecal_\perp}V_2\|\|e_s\|\notag\\
&\hspace{20pt}+ c_J |r(s,a,\mu_1)-r(s,a,\mu_2)|+c_J|J_1-J_2|\notag\\
&\leq (c_J+1)|r(s,a,\mu_1)-r(s,a,\mu_2)|+(c_J+1)|J_1-J_2|+2\|\Pi_{\Ecal_\perp}V_1-\Pi_{\Ecal_\perp}V_2\|\notag\\
&\leq (c_J+1)L|\mu_1-\mu_2|+(c_J+1)|J_1-J_2|+2\|\Pi_{\Ecal_\perp}V_1-\Pi_{\Ecal_\perp}V_2\|.\label{lem:Lipschitz_operators:eq4}
\end{align}

Plugging \eqref{lem:Lipschitz_operators:eq4} into \eqref{lem:Lipschitz_operators:eq3}, we get
\begin{align*}
&\|\bar{G}(\theta_1,V_1,J_1,\mu_1)-\bar{G}(\theta_2,V_2,J_2,\mu_2)\|\notag\\
&\leq \|\mathbb{E}_{s\sim \nu^{\pi_{\theta_1},\,\mu_1},a\sim\pi_{\theta_1}(\cdot\mid s),s'\sim\Pcal^{\mu_1}(\cdot\mid s,a,\mu_1)}[G(\Pi_{\Ecal_\perp}V_1,J_1,\mu_1,s,a,s')-G(\Pi_{\Ecal_\perp}V_2,J_2,\mu_2,s,a,s')]\|\notag\\
&\hspace{20pt}+2B_G d_{TV}(\nu^{\pi_{\theta_1},\,\mu_1}\otimes\pi_{\theta_1}\otimes\Pcal^{\mu_1},\nu^{\pi_{\theta_2},\,\mu_2}\otimes\pi_{\theta_2}\otimes\Pcal^{\mu_2})\notag\\
&\leq (c_J+1)L|\mu_1-\mu_2|+(c_J+1)|J_1-J_2|+2\|\Pi_{\Ecal_\perp}V_1-\Pi_{\Ecal_\perp}V_2\|\notag\\
&\hspace{20pt}+2B_G L_{TV}(\|\theta_1-\theta_2\|+\|\mu_1-\mu_2\|)\notag\\
&\leq L_G\left(\|\theta_1-\theta_2\|+\|\mu_1-\mu_2\|+2\|\Pi_{\Ecal_\perp}V_1-\Pi_{\Ecal_\perp}V_2\|+|J_1-J_2|\right),
\end{align*}
with $L_G=2B_G L_{TV}+(L+1)(c_J+1)+2$.

Finally, again following steps similar to \eqref{lem:Lipschitz_operators:eq1} we can show
\begin{align}
&\|\bar{H}(\theta_1,\mu_1)-\bar{H}(\theta_2,\mu_2)\|\notag\\
&\leq\|\mathbb{E}_{s\sim \nu^{\pi_{\theta_1},\,\mu_1}}[H(\mu_1,s)-H(\mu_2,s)]\|+2B_H d_{TV}(\nu^{\pi_{\theta_1},\,\mu_1},\nu^{\pi_{\theta_2},\,\mu_2}).\label{lem:Lipschitz_operators:eq5}
\end{align}

From the definition of $H(\mu,s)$ in \eqref{eq:def_FGH}, we have for any $s$
\begin{align}
\|H(\mu_1,s)-H(\mu_2,s)\|=\|(e_s-\mu_1)-(e_s-\mu_2)\|=\|\mu_1-\mu_2\|.\label{lem:Lipschitz_operators:eq6}
\end{align}

By Assumption~\ref{assump:Lipschitz},
\begin{align}
d_{TV}(\nu^{\pi_{\theta_1},\,\mu_1},\nu^{\pi_{\theta_2},\,\mu_2})&=\frac{1}{2}\|\nu^{\pi_{\theta_1},\,\mu_1}-\nu^{\pi_{\theta_2},\,\mu_2}\|_1\notag\\
&\leq L(\|\pi_{\theta_1}-\pi_{\theta_2}\|+\|\mu_1-\mu_2\|)\notag\\
&\leq L(\|\theta_1-\theta_2\|+\|\mu_1-\mu_2\|),
\label{lem:Lipschitz_operators:eq7}
\end{align}
where the final inequality is a result of the $1$-Lipschitz continuity of the softmax function.

Plugging \eqref{lem:Lipschitz_operators:eq6} and \eqref{lem:Lipschitz_operators:eq7} into \eqref{lem:Lipschitz_operators:eq5}, we have
\begin{align}
\|\bar{H}(\theta_1,\mu_1)-\bar{H}(\theta_2,\mu_2)\|
&\leq\|\mathbb{E}_{s\sim \nu^{\pi_{\theta_1},\,\mu_1}}[H(\mu_1,s)-H(\mu_2,s)]\|+2B_H d_{TV}(\nu^{\pi_{\theta_1},\,\mu_1},\nu^{\pi_{\theta_2},\,\mu_2})\notag\\
&\leq\|\mu_1-\mu_2\|+L(\|\theta_1-\theta_2\|+\|\mu_1-\mu_2\|)\notag\\
&\leq (L+1)(\|\theta_1-\theta_2\|+\|\mu_1-\mu_2\|).
\end{align}

\qed

\subsection{Proof of Lemma~\ref{lem:bound_f_g}}

By the definition $\Delta f_k$,
\begin{align*}
&\|f_k\|\notag\\
&=\|\Delta f_k+\bar{F}(\theta_k,\hat{V}_k,\hat{\mu}_k)-\bar{F}(\theta_k,V^{\pi_{\theta_k},\mu^{\star}(\pi_{\theta_k})},\mu^{\star}(\pi_{\theta_k}))+\bar{F}(\theta_k,V^{\pi_{\theta_k},\mu^{\star}(\pi_{\theta_k})},\mu^{\star}(\pi_{\theta_k}))\|\notag\\
&\leq \|\Delta f_k\|+\|\bar{F}(\theta_k,\hat{V}_k,\hat{\mu}_k)-\bar{F}(\theta_k,V^{\pi_{\theta_k},\mu^{\star}(\pi_{\theta_k})},\mu^{\star}(\pi_{\theta_k}))\|+\|\bar{F}(\theta_k,V^{\pi_{\theta_k},\mu^{\star}(\pi_{\theta_k})},\mu^{\star}(\pi_{\theta_k}))\|\notag\\
&\leq \|\Delta f_k\|+L_F\|\Pi_{\Ecal_\perp}(V^{\pi_{\theta_k},\mu^{\star}(\pi_{\theta_k})}-\hat{V}_k)\|+L_F\|\hat{\mu}_k-\mu^{\star}(\pi_{\theta_k})\|+\|\nabla_{\theta} J(\pi_{\theta_k},\mu)\mid_{\mu=\mu^{\star}(\pi_{\theta_k})}\|\notag\\
&\leq \|\Delta f_k\|+L_F\|\Pi_{\Ecal_\perp}(V^{\pi_{\theta_k},\mu^{\star}(\pi_{\theta_k})}-V^{\pi_{\theta_k},\hat{\mu}_k})\|+L_F\|\Pi_{\Ecal_\perp}(V^{\pi_{\theta_k},\hat{\mu}_k} -\hat{V}_k)\|\notag\\
&\hspace{20pt}+L_F\|\hat{\mu}_k-\mu^{\star}(\pi_{\theta_k})\|+\sqrt{\varepsilon_k^{\pi}}\notag\\
&\leq \|\Delta f_k\|+L_F\sqrt{\varepsilon_k^{V}}+L_F (L_V+1)\sqrt{\varepsilon_k^{\mu}}+\sqrt{\varepsilon_k^{\pi}},
\end{align*}
where the last inequality follows from the Lipschitz continuity of the value function in the mean field and the fact that linear projection is non-expansive, and the second inequality follows from the Lipschitz continuity of operator $F$ and the relationship 
\[\nabla_{\theta} J(\pi_{\theta_k},\mu)\mid_{\mu=\mu^{\star}(\pi_{\theta_k})}=\bar{F}(\theta_k,V^{\pi_{\theta_k},\mu^{\star}(\pi_{\theta_k})},\mu^{\star}(\pi_{\theta_k})).\]

Similarly, by the definition of $\Delta g_k$, we have
\begin{align*}
\|g_k\|&=\|\Delta g_k+\bar{G}(\theta_k,\hat{V}_k,\hat{J}_k,\hat{\mu}_k)-\bar{G}(\theta_k,V^{\pi_{\theta_k},\hat{\mu}_k},J(\pi_{\theta_k},\hat{\mu}_k),\hat{\mu}_k)\|\notag\\
&\leq \|\Delta g_k\|+\|\bar{G}(\theta_k,\hat{V}_k,\hat{J}_k,\hat{\mu}_k)-\bar{G}(\theta_k,V^{\pi_{\theta_k},\hat{\mu}_k},J(\pi_{\theta_k},\hat{\mu}_k),\hat{\mu}_k)\|\notag\\
&\leq \|\Delta g_k\|+L_G\|\Pi_{\Ecal_\perp}(V^{\pi_{\theta_k},\hat{\mu}_k}-\hat{V}_k)\|+L_G|J(\pi_{\theta_k},\hat{\mu}_k)-\hat{J}_k|\notag\\
&=\|\Delta g_k\|+L_G\sqrt{\varepsilon_k^V}+L_G\sqrt{\varepsilon_k^J},
\end{align*}
where the first equation follows from the fact that $G(\theta_k,V^{\pi_{\theta_k},\hat{\mu}_k},J(\pi_{\theta_k},\hat{\mu}_k),\hat{\mu}_k)=0$.

Finally, by the definition of $\Delta h_k$, we have
\begin{align*}
\|h_k\|&=\|\Delta h_k+\bar{H}(\theta_k,\hat{\mu}_k)\|\notag\\
&=\|\Delta h_k+\bar{H}(\theta_k,\hat{\mu}_k)-\bar{H}(\theta_k,\mu^{\star}(\pi_{\theta_k}))\|\notag\\
&\leq \|\Delta h_k\|+\|\bar{H}(\theta_k,\hat{\mu}_k)-\bar{H}(\theta_k,\mu^{\star}(\pi_{\theta_k}))\|\notag\\
&\leq \|\Delta h_k\|+L_H\|\hat{\mu}_k-\mu^{\star}(\pi_{\theta_k})\|\notag\\
&= \|\Delta h_k\|+L_H\sqrt{\epsilon_k^{\mu}},
\end{align*}
where the second equation follows from the fact that $H(\theta_k,\mu^{\star}(\pi_{\theta_k}))=0$.

\qed

\subsection{Proof of Lemma~\ref{lem:negative_drift}}

See \citet{zhang2021finite}[Lemma 2] or \citet{tsitsiklis1999average}[Lemma 7].

\subsection{Proof of Lemma~\ref{lem:gradient_domination}}

Adapted from Lemma 19 of \citet{ganesh2024variance}.

\subsection{Proof of Lemma~\ref{lem:f_markovian}}\label{sec:f_markovian:proof}

The proof of this lemma proceeds in a manner similar to that of Lemma~\ref{lem:h_markovian}. We note that the samples generated in the algorithm follow the time-varying Markov chain
\begin{align}
s_{k-\tau_k} \stackrel{\theta_{k-\tau_k}}{\longrightarrow}  a_{k-\tau_k} \stackrel{\hat{\mu}_{k-\tau_k}}{\longrightarrow}  s_{k-\tau_k+1} \stackrel{\theta_{k-\tau_k+1}}{\longrightarrow} a_{k-\tau_k+1} \stackrel{\hat{\mu}_{k-\tau_k+1}}{\longrightarrow}  \cdots s_{k-1} \stackrel{\theta_{k-1}}{\longrightarrow} a_{k-1} \stackrel{\hat{\mu}_{k-1}}{\longrightarrow}  s_{k}.
\end{align}
We construct an auxiliary Markov chain generated under a constant control
\begin{align}
s_{k-\tau_k} \stackrel{\theta_{k-\tau_k}}{\longrightarrow}  a_{k-\tau_k} \stackrel{\hat{\mu}_{k-\tau_k}}{\longrightarrow}  \widetilde{s}_{k-\tau_k+1} \stackrel{\theta_{k-\tau_k}}{\longrightarrow} \widetilde{a}_{k-\tau_k+1} \stackrel{\hat{\mu}_{k-\tau_k}}{\longrightarrow}  \cdots \widetilde{s}_{k-1} \stackrel{\theta_{k-\tau_k}}{\longrightarrow} \widetilde{a}_{k-1} \stackrel{\hat{\mu}_{k-\tau_k}}{\longrightarrow}  \widetilde{s}_{k}
\label{lem:f_markovian:distribution_auxiliary}
\end{align}

Let $\widetilde{\mu}$ denote the stationary distribution of state, action, and next state under \eqref{lem:f_markovian:distribution_auxiliary}.
We denote $p_k(s,a,s')=\mathbb{P}(s_k=s,a_k=a,s_{k+1}=s')$ and $\widetilde{p}_k(s,a,s')=\mathbb{P}(\widetilde{s}_k=s,\widetilde{a}_k=a,\widetilde{s}_{k+1}=s')$ and define
\begin{align*}
T_1&\triangleq\mathbb{E}[\langle\Delta f_k-\Delta f_{k-\tau_k}, F(\theta_{k},\hat{V}_{k},\hat{\mu}_{k},s_k,a_k,s_{k+1})-\bar{F}(\theta_{k},\hat{V}_{k},\hat{\mu}_{k})\rangle],\\
T_2&\triangleq\mathbb{E}[\langle\Delta f_{k-\tau_k}, F(\theta_{k},\hat{V}_{k},\hat{\mu}_{k},s_k,a_k,s_{k+1})-F(\theta_{k},\hat{V}_{k},\hat{\mu}_{k},\widetilde{s}_k,\widetilde{a}_k,\widetilde{s}_{k+1})\rangle],\\
T_3&\triangleq\mathbb{E}[\langle\Delta f_{k-\tau_k},F(\theta_{k},\hat{V}_{k},\hat{\mu}_{k},\widetilde{s}_k,\widetilde{a}_k,\widetilde{s}_{k+1})-\mathbb{E}_{(s,a,s')\sim\widetilde{\mu}}[F(\theta_{k},\hat{V}_{k},\hat{\mu}_{k},s,a,s')]\rangle]\\
T_4&\triangleq\mathbb{E}[\langle\Delta f_{k-\tau_k},\mathbb{E}_{(s,a,s')\sim\widetilde{\mu}}[F(\theta_{k},\hat{V}_{k},\hat{\mu}_{k},s,a,s')]-\bar{F}(\theta_{k},\hat{V}_{k},\hat{\mu}_{k})\rangle].
\end{align*}

It is obvious to see
\begin{align}
\mathbb{E}[\langle\Delta f_k, F(\theta_{k},\hat{V}_{k},\hat{\mu}_{k},s_k,a_k,s_{k+1})-\bar{F}(\theta_{k},\hat{V}_{k},\hat{\mu}_{k})\rangle]=T_1+T_2+T_3+T_4.\label{lem:f_markovian:proof_eq1}
\end{align}

We bound the terms individually. First, we treat $T_1$
\begin{align*}
T_1 &= \mathbb{E}[\langle\Delta f_k-\Delta f_{k-\tau_k}, F(\theta_{k},\hat{V}_{k},\hat{\mu}_{k},s_k,a_k,s_{k+1})-\bar{F}(\theta_{k},\hat{V}_{k},\hat{\mu}_{k})\rangle]\notag\\
&\leq \mathbb{E}[\|f_k- f_{k-\tau_k}\|\|F(\theta_{k},\hat{V}_{k},\hat{\mu}_{k},s_k,a_k,s_{k+1})-\bar{F}(\theta_{k},\hat{V}_{k},\hat{\mu}_{k})\|]\notag\\
&\hspace{20pt}+\mathbb{E}\Big[\|\bar{F}(\theta_k,\hat{V}_k,\hat{\mu}_k)-\bar{F}(\theta_{k-\tau_k},\hat{V}_{k-\tau_k},\hat{\mu}_{k-\tau_k})\|\notag\\
&\hspace{100pt}\cdot\|F(\theta_{k},\hat{V}_{k},\hat{\mu}_{k},s_k,a_k,s_{k+1})-\bar{F}(\theta_{k},\hat{V}_{k},\hat{\mu}_{k})\|\Big]\notag\\
&\leq 2B_F\sum_{t=0}^{\tau_k-1}\mathbb{E}[\|f_{k-t}-f_{k-t-1}\|]\notag\\
&\hspace{20pt}+2L_F B_F\sum_{t=0}^{\tau_k-1}\mathbb{E}[\|\theta_{k-t}-\theta_{k-t-1}\|+\|\hat{V}_{k-t}-\hat{V}_{k-t-1}\|+\|\hat{\mu}_{k-t}-\hat{\mu}_{k-t-1}\|]\notag\\
&\leq 4B_F^2\tau_k\lambda_{k-\tau_k}+2L_F B_F\tau_k(B_F\alpha_{k-\tau_k}+B_G\beta_{k-\tau_k}+B_H\xi_{k-\tau_k})\notag\\
&\leq 10 L_F B_F^2 B_G B_H\tau_k\lambda_{k-\tau_k},
\end{align*}
where the second inequality bounds $\|F(\theta_{k},\hat{V}_{k},\hat{\mu}_{k},s_k,a_k,s_{k+1})-\bar{F}(\theta_{k},\hat{V}_{k},\hat{\mu}_{k})\|$ by $2B_F$ and $\|\bar{F}(\theta_k,\hat{V}_k,\hat{\mu}_k)-\bar{F}(\theta_{k-\tau_k},\hat{V}_{k-\tau_k},\hat{\mu}_{k-\tau_k})\|$ using the Lipschitz continuity established in Lemma~\ref{lem:Lipschitz_operators}.
The last inequality follows from the step size condition $\alpha_k\leq\xi_k\leq\beta_k\leq\lambda_k$ for all $k$. The third inequality follows from the fact that $\|f_{k+1}-f_k\|=\lambda_k\|f_k-F(\theta_{k},\hat{V}_{k},\hat{\mu}_{k},s_k,a_k,s_{k+1})\|\leq 2B_F\lambda_k$ for all $k$ and that the per-iteration drift of $\theta_k$, $\hat{V}_k$, and $\hat{\mu}_k$ can be similarly bounded
\begin{align*}
\|\theta_{k+1}-\theta_k\|\leq B_F\alpha_k,\quad \|\hat{V}_{k+1}-\hat{V}_k\|\leq B_G\beta_k,\quad \|\hat{\mu}_{k+1}-\hat{\mu}_k\|\leq B_H\xi_k.
\end{align*}

We next bound $T_2$
\begin{align*}
T_2&=\mathbb{E}[\langle\Delta f_{k-\tau_k}, F(\theta_{k},\hat{V}_{k},\hat{\mu}_{k},s_k,a_k,s_{k+1})-F(\theta_{k},\hat{V}_{k},\hat{\mu}_{k},\widetilde{s}_k,\widetilde{a}_k,\widetilde{s}_{k+1})\rangle]\notag\\
&\leq 2B_F \mathbb{E}_{\Fcal_{k-\tau_k}}[\mathbb{E}[\|F(\theta_{k},\hat{V}_{k},\hat{\mu}_{k},s_k,a_k,s_{k+1})-F(\theta_{k},\hat{V}_{k},\hat{\mu}_{k},\widetilde{s}_k,\widetilde{a}_k,\widetilde{s}_{k+1})\|\mid \Fcal_{k-\tau_k}]]\notag\\
&\leq 2B_F \mathbb{E}[\int_{\Scal}\int_{\Acal}\int_{\Scal}F(\theta_{k},\hat{V}_{k},\hat{\mu}_{k},s,a,s')\left(p_k(s,a,s')-\widetilde{p}_k(s,a,s')\right)ds\,da\,ds']\notag\\
&\leq 2B_F^2\mathbb{E}[d_{TV}(p_k,\widetilde{p}_k)].
\end{align*}
where the last inequality follows from the definition of TV distance in \eqref{eq:TV_def}.

Applying Lemma B.2 from \citet{wu2020finite}, we then have
\begin{align*}
T_2&\leq 2B_F^2\mathbb{E}[d_{TV}(p_k,\widetilde{p}_k)]\notag\\
&\leq 2B_F^2\mathbb{E}[d_{TV}(\mathbb{P}(s_k=\cdot),\mathbb{P}(\widetilde{s}_k=\cdot))+\frac{|\Acal|}{2}\|\theta_{k-1}-\theta_{k-\tau_k}\|]\notag\\
&\leq 2B_F^2\mathbb{E}\Big[d_{TV}(\mathbb{P}(s_{k-1}=\cdot),\mathbb{P}(\widetilde{s}_{k-1}=\cdot))+L\|\theta_{k-1}-\theta_{k-\tau_k}\|+L\|\hat{\mu}_{k-1}-\hat{\mu}_{k-\tau_k}\|\notag\\
&\hspace{20pt}+\frac{|\Acal|}{2}\|\theta_{k-1}-\theta_{k-\tau_k}\|\Big]\notag\\
&\leq |\Acal|B_F^2\mathbb{E}[\|\theta_{k-1}-\theta_{k-\tau_k}\|]+2LB_F^2 \sum_{t=k-\tau_k}^{k-1}\mathbb{E}[\|\theta_t-\theta_{k-\tau_k}\|+\|\hat{\mu}_t-\hat{\mu}_{k-\tau_k}\|]\notag\\
&\leq (2L+|\Acal|) B_F^2 \tau_k^2(B_F\alpha_{k-\tau_k}+B_H\xi_{k-\tau_k})\notag\\
&\leq (4L+2|\Acal|) B_F^3 B_H\tau_k^2 \lambda_{k-\tau_k},
\end{align*}
where the third inequality is a result of \eqref{assump:Lipschitz:eq1}, and the fourth inequality recursively applies the inequality above it.

The term $T_3$ is proportional to the distance between the distribution of the auxiliary Markov chain \eqref{lem:f_markovian:distribution_auxiliary} at time $k$ and its stationary distribution. To bound $T_3$, 
\begin{align*}
T_3&=\mathbb{E}[\langle\Delta f_{k-\tau_k},F(\theta_{k},\hat{V}_{k},\hat{\mu}_{k},\widetilde{s}_k,\widetilde{a}_k,\widetilde{s}_{k+1})-\mathbb{E}_{(s,a,s')\sim\widetilde{\mu}}[F(\theta_{k},\hat{V}_{k},\hat{\mu}_{k},s,a,s')]\rangle]\notag\\
&\leq 2B_F \mathbb{E}_{\Fcal_{k-\tau_k}}[\mathbb{E}[\|F(\theta_{k},\hat{V}_{k},\hat{\mu}_{k},\widetilde{s}_k,\widetilde{a}_k,\widetilde{s}_{k+1})-\mathbb{E}_{(s,a,s')\sim\widetilde{\mu}}[F(\theta_{k},\hat{V}_{k},\hat{\mu}_{k},s,a,s')]\|\mid \Fcal_{k-\tau_k}]]\notag\\
&\leq 2B_F\mathbb{E}[\int_{\Scal}\int_{\Acal}\int_{\Scal}F(\theta_{k},\hat{V}_{k},\hat{\mu}_{k},s,a,s')\left(\widetilde{p}_k(s)-\widetilde{\mu}(s)\right)ds\, da\, ds']\notag\\
&\leq 2B_F^2\mathbb{E}[d_{TV}(\widetilde{p}_k,\widetilde{\mu})]\notag\\
&\leq 2B_F^2\alpha_k,
\end{align*}
where the final inequality follows from the definition of the mixing time $\tau_k$ as the number of iterations for the TV distance between $\widetilde{p}_k$ and $\widetilde{\mu}$ to drop below $\alpha_k$.

Finally, we bound the term $T_4$
\begin{align*}
T_4 &= \mathbb{E}[\langle\Delta f_{k-\tau_k},\mathbb{E}_{(s,a,s')\sim\widetilde{\mu}}[F(\theta_{k},\hat{V}_{k},\hat{\mu}_{k},s,a,s')]-\bar{F}(\theta_{k},\hat{V}_{k},\hat{\mu}_{k})\rangle]\notag\\
&\leq 2B_F\mathbb{E}[\|\mathbb{E}_{(s,a,s')\sim\widetilde{\mu}}[F(\theta_{k},\hat{V}_{k},\hat{\mu}_{k},s,a,s')]-\bar{F}(\theta_{k},\hat{V}_{k},\hat{\mu}_{k})\|]\notag\\
&\leq 2B_F^2 \mathbb{E}[d_{TV}(\widetilde{\mu},\nu^{\pi_{\theta_k},\,\hat{\mu}_{k}}\otimes\pi_{\theta_k}\otimes\Pcal^{\hat{\mu}_{k}})]\notag\\
&\leq 2 L_{TV} B_F^2\mathbb{E}[\|\pi_{\theta_k}-\pi_{\theta_{k-\tau_k}}\|+\|\hat{\mu}_k-\hat{\mu}_{k-\tau_k}\|]\notag\\
&\leq 2 L_{TV} B_F^2\tau_k(B_F\alpha_{k-\tau_k}+B_H\xi_{k-\tau_k})\notag\\
&\leq 4 L_{TV} B_F^3 B_H\xi_{k-\tau_k},
\end{align*}
where the third inequality applies the result in \eqref{eq:TV_stationary_Lipschitz}.

Collecting the bounds on $T_1$-$T_4$ and plugging them into \eqref{lem:f_markovian:proof_eq1}, we get
\begin{align*}
&\mathbb{E}[\langle\Delta f_k, F(\theta_{k},\hat{V}_{k},\hat{\mu}_{k},s_k,a_k,s_{k+1})-\bar{F}(\theta_{k},\hat{V}_{k},\hat{\mu}_{k})\rangle]\notag\\
&=T_1+T_2+T_3+T_4\notag\\
&\leq 10 L_F B_F^2 B_G B_H\tau_k\lambda_{k-\tau_k}
+ (4L+2|\Acal|) B_F^3 B_H\tau_k^2 \lambda_{k-\tau_k} + 2B_F^2\alpha_k + 4 L_{TV} B_F^3 B_H\xi_{k-\tau_k}\notag\\
&\leq (20L+2|\Acal|) L_F L_{TV} B_F^3 B_G B_H^2\tau_k^2\lambda_{k-\tau_k}.
\end{align*}

\qed

\subsection{Proof of Lemma~\ref{lem:H_strongmonotone}}\label{sec:H_strongmonotone:proof}

By the definition of $\bar{H}$, we have for any $\mu\in\Delta_{\Scal}$
\begin{align*}
&\langle \mu-\mu^{\star}(\pi_{\theta}),\bar{H}(\theta,\mu)-\bar{H}(\theta,\mu^{\star}(\pi_{\theta}))\rangle\notag\\
&=\langle \mu-\mu^{\star}(\pi_{\theta}),\mu^{\star}(\pi_{\theta})-\mu\rangle + \langle \mu-\mu^{\star}(\pi_{\theta}),\nu^{\pi_{\theta},\,\mu}-\nu^{\pi_{\theta},\,\mu^{\star}(\pi_{\theta})}\rangle\notag\\
&\leq-\|\mu-\mu^{\star}(\pi_{\theta})\|^2+\|\mu-\mu^{\star}(\pi_{\theta})\|\|\nu^{\pi_{\theta},\,\mu}-\nu^{\pi_{\theta},\,\mu^{\star}(\pi_{\theta})}\|\notag\\
&\leq-(1-\delta)\|\mu-\mu^{\star}(\pi_{\theta})\|^2,
\end{align*}
where the second inequality follows from Assumption~\ref{assump:nu}.

\qed

\subsection{Proof of Lemma~\ref{lem:h_markovian}}\label{sec:h_markovian:proof}

The cause of the gap between $\mathbb{E}[e_{s_k}]$ and $\mathbb{E}_{s\sim\nu^{\pi_{\theta_k},\,\hat{\mu}_k}}[e_s]$ is a time-varying Markovian noise. To elaborate, we first show how the sample $s_k$ is generated below
\begin{align}
s_{k-\tau_k} \stackrel{\theta_{k-\tau_k},\,\hat{\mu}_{k-\tau_k}}{\longrightarrow}  s_{k-\tau_k+1} \stackrel{\theta_{k-\tau_k+1},\,\hat{\mu}_{k-\tau_k+1}}{\longrightarrow}  \cdots s_{k-1} \stackrel{\theta_{k-1},\,\hat{\mu}_{k-1}}{\longrightarrow}  s_{k}.
\end{align}
This Markov chain is ``time-varying'' as its stationary distribution changes over iterations as the control changes.
We introduce an auxiliary Markov chain, which is ``time-invariant'' in the sense that it is generated under a constant control, starting from state $s_{k-\tau_k}$.
\begin{align}
s_{k-\tau_k} \stackrel{\theta_{k-\tau_k},\,\hat{\mu}_{k-\tau_k}}{\longrightarrow}  \widetilde{s}_{k-\tau_k+1} \stackrel{\theta_{k-\tau_k},\,\hat{\mu}_{k-\tau_k}}{\longrightarrow}  \cdots \widetilde{s}_{k-1} \stackrel{\theta_{k-\tau_k},\,\hat{\mu}_{k-\tau_k}}{\longrightarrow}  \widetilde{s}_{k}.
\label{lem:h_markovian:distribution_auxiliary}
\end{align}

Defining
\begin{align*}
T_1 &\triangleq \mathbb{E}[\langle\Delta h_k-\Delta h_{k-\tau_k}, e_{s_k}-\mathbb{E}_{s\sim\nu^{\pi_{\theta_k},\,\hat{\mu}_k}}[e_s]\rangle]\\
T_2 &\triangleq \mathbb{E}[\langle\Delta h_{k-\tau_k}, e_{s_k}-e_{\widetilde{s}_k}\rangle]\\
T_3 &\triangleq \mathbb{E}[\langle\Delta h_{k-\tau_k}, e_{\widetilde{s}_k}-\mathbb{E}_{s\sim\nu^{\pi_{\theta_{k-\tau_k}},\,\hat{\mu}_{k-\tau_k}}}[e_{s}]\rangle]\\
T_4 &\triangleq \mathbb{E}[\langle\Delta h_{k-\tau_k}, \mathbb{E}_{s\sim\nu^{\pi_{\theta_{k-\tau_k}},\,\hat{\mu}_{k-\tau_k}}}[e_{s}]-\mathbb{E}_{s\sim\nu^{\pi_{\theta_k},\,\hat{\mu}_k}}[e_s]\rangle],
\end{align*}
we see that 
\begin{align}
    \mathbb{E}[\langle\Delta h_k, e_{s_k}-\mathbb{E}_{s\sim\nu^{\pi_{\theta_k},\,\hat{\mu}_k}}[e_s]\rangle]=T_1+T_2+T_3+T_4.\label{lem:h_markovian:proof_eq1}
\end{align}

We bound the terms individually. First, we treat $T_1$
\begin{align*}
T_1 &= \mathbb{E}[\langle h_k- h_{k-\tau_k}+\bar{H}(\theta_{k-\tau_k},\hat{\mu}_{k-\tau_k})-\bar{H}(\theta_k,\hat{\mu}_k), e_{s_k}-\mathbb{E}_{s\sim\nu^{\pi_{\theta_k},\,\hat{\mu}_k}}[e_s]\rangle]\notag\\
&\leq \mathbb{E}[\|h_k- h_{k-\tau_k}\|\|e_{s_k}-\mathbb{E}_{s\sim\nu^{\pi_{\theta_k},\,\hat{\mu}_k}}[e_s]\|]\notag\\
&\hspace{20pt}+\mathbb{E}[\|\bar{H}(\theta_k,\hat{\mu}_k)-\bar{H}(\theta_{k-\tau_k},\hat{\mu}_{k-\tau_k})\|\|e_{s_k}-\mathbb{E}_{s\sim\nu^{\pi_{\theta_k},\,\hat{\mu}_k}}[e_s]\|]\notag\\
&\leq 2\sum_{t=0}^{\tau_k-1}\mathbb{E}[\|h_{k-t}-h_{k-t-1}\|]+2L_H\sum_{t=0}^{\tau_k-1}\mathbb{E}[\|\theta_{k-t}-\theta_{k-t-1}\|+\|\hat{\mu}_{k-t}-\hat{\mu}_{k-t-1}\|]\notag\\
&\leq 4B_H\tau_k\lambda_{k-\tau_k}+2B_F\tau_k\alpha_{k-\tau_k}+2B_H\tau_k\xi_{k-\tau_k}\notag\\
&\leq 8 B_F B_H\tau_k\lambda_{k-\tau_k},
\end{align*}
where the last inequality follows from the step size condition $\alpha_k\leq\xi_k\leq\lambda_k$ for all $k$, and the third inequality follows from the fact that $\|h_{k+1}-h_{k}\| \leq \lambda_k\|h_k+\hat{\mu}_k-e_{s_k}\|\leq 2B_H\lambda_k$ for all $k$ and that the per-iteration drift of $\theta_k$ and $\hat{\mu}_k$ can be similarly bounded
\begin{align*}
\|\theta_{k+1}-\theta_k\|\leq B_F\alpha_k,\quad \|\hat{\mu}_{k+1}-\hat{\mu}_k\|\leq B_H\xi_k.
\end{align*}

We next bound $T_2$. We denote $p_k(s)=\mathbb{P}(s_k=s)$ and $\widetilde{p}_k(s)=\mathbb{P}(\widetilde{s}_k=s)$.
\begin{align*}
T_2&=\mathbb{E}_{\Fcal_{k-\tau_k}}[\mathbb{E}[\langle h_{k-\tau_k}-\bar{H}(\theta_{k-\tau_k},\hat{\mu}_{k-\tau_k}), e_{s_k}-e_{\widetilde{s}_k}\rangle\mid \Fcal_{k-\tau_k}]]\notag\\
&\leq 2B_H \mathbb{E}_{\Fcal_{k-\tau_k}}[\mathbb{E}[\|e_{s_k}-e_{\widetilde{s}_k}\|\mid \Fcal_{k-\tau_k}]]\notag\\
&\leq 2B_H\mathbb{E}[\int_{\Scal}e_s\left(p_k(s)-\widetilde{p}_k(s)\right)ds]\notag\\
&\leq 2B_H\mathbb{E}[d_{TV}(p_k,\widetilde{p}_k)]\notag\\
&\leq 2B_H\mathbb{E}[d_{TV}(p_{k-1},\widetilde{p}_{k-1})+L\|\theta_{k-1}-\theta_{k-\tau_k}\|+L\|\hat{\mu}_{k-1}-\hat{\mu}_{k-\tau_k}\|]\notag\\
&\leq 2LB_H \sum_{t=k-\tau_k}^{k-1}\mathbb{E}[\|\theta_t-\theta_{k-\tau_k}\|+\|\hat{\mu}_t-\hat{\mu}_{k-\tau_k}\|]\notag\\
&\leq 2L B_H\tau_k^2(B_F\alpha_{k-\tau_k}+B_H\xi_{k-\tau_k})\notag\\
&\leq 4L B_F B_H^2\tau_k^2 \lambda_{k-\tau_k},
\end{align*}
where the third inequality follows from the definition of TV distance in \eqref{eq:TV_def}, and the fourth and fifth inequalities are a result of \eqref{assump:Lipschitz:eq1}.

The term $T_3$ is proportional to the distance between the distribution of the auxiliary Markov chain \eqref{lem:h_markovian:distribution_auxiliary} at time $k$ and its stationary distribution. Let $\widetilde{\mu}$ denote the stationary distribution of \eqref{lem:h_markovian:distribution_auxiliary}. We can bound this term as follows under Assumption~\ref{assump:ergodic}
\begin{align*}
T_3&=\mathbb{E}[\langle\Delta h_{k-\tau_k}, e_{\widetilde{s}_k}-\mathbb{E}_{s\sim\nu^{\pi_{\theta_{k-\tau_k}},\,\hat{\mu}_{k-\tau_k}}}[e_{s}]\rangle]\notag\\
&\leq 2B_H \mathbb{E}_{\Fcal_{k-\tau_k}}[\mathbb{E}[\|e_{\widetilde{s}_k}-\mathbb{E}_{s\sim\nu^{\pi_{\theta_{k-\tau_k}},\,\hat{\mu}_{k-\tau_k}}}[e_{s}]\|\mid \Fcal_{k-\tau_k}]]\notag\\
&\leq 2B_H\mathbb{E}[\int_{\Scal}e_s\left(\widetilde{p}_k(s)-\widetilde{\mu}(s)\right)ds]\notag\\
&\leq 2B_H\mathbb{E}[d_{TV}(\widetilde{p}_k,\widetilde{\mu})]\notag\\
&\leq 2B_H\alpha_k,
\end{align*}
where the final inequality follows from the definition of the mixing time $\tau_k$ as the number of iterations for the TV distance between $\widetilde{p}_k$ and $\widetilde{\mu}$ to drop below $\alpha_k$.

The term $T_4$ can be treated by the Lipschitz continuity of $\nu$
\begin{align*}
T_4 &= \mathbb{E}[\langle\Delta h_{k-\tau_k}, \mathbb{E}_{s\sim\nu^{\pi_{\theta_{k-\tau_k}},\,\hat{\mu}_{k-\tau_k}}}[e_{s}]-\mathbb{E}_{s\sim\nu^{\pi_{\theta_k},\,\hat{\mu}_k}}[e_s]\rangle]\notag\\
&\leq 2B_H\mathbb{E}[\|\nu^{\pi_{\theta_{k-\tau_k}},\,\hat{\mu}_{k-\tau_k}}-\nu^{\pi_{\theta_k},\,\hat{\mu}_k}\|]\notag\\
&\leq 2B_H L\mathbb{E}[\|\pi_{\theta_{k}}-\pi_{\theta_{k-\tau_k}}\|]+2B_H \delta \mathbb{E}[\|\hat{\mu}_{k-\tau_k}-\hat{\mu}_k\|]\notag\\
&\leq 2B_H L\sum_{t=k-\tau_k}^{k}\mathbb{E}[\|\alpha_t f_t\|] + 2B_H L\sum_{t=k-\tau_k}^{k} \mathbb{E}[\|\xi_t h_t\|]\notag\\
&\leq 2B_H L\tau_k\left(B_F\alpha_{k-\tau_k} + B_H\xi_{k-\tau_k} \right)\notag\\
&\leq 2 L B_F B_H^2\xi_{k-\tau_k}
\end{align*}
where the last inequality follows from the step size condition $\alpha_k\leq\xi_k$ for all $k$.

Collecting the bounds on $T_1$-$T_4$ and plugging them into \eqref{lem:h_markovian:proof_eq1}, we get
\begin{align*}
&\mathbb{E}[\langle\Delta h_k, e_{s_k}-\mathbb{E}_{s\sim\nu^{\pi_{\theta_k},\,\hat{\mu}_k}}[e_s]\rangle]\notag\\
&=T_1+T_2+T_3+T_4\notag\\
&\leq 8 B_F B_H\tau_k\lambda_{k-\tau_k}
+4L B_F B_H^2\tau_k^2 \lambda_{k-\tau_k}+2B_H\alpha_k+2 L B_F B_H^2\xi_{k-\tau_k}\notag\\
&\leq 16L B_F B_H^2\tau_k^2\lambda_{k-\tau_k}.
\end{align*}

\qed

\subsection{Proof of Lemma~\ref{lem:G_strongmonotone}}\label{sec:G_strongmonotone:proof}

By the definition of operators $G^V$ and $G^J$ in \eqref{eq:def_FGH}, for any $V\in\mathbb{R}^{|\Scal|}$ and $J\in\mathbb{R}$
\begin{align*}
&\left\langle \left[\begin{array}{c}
\Pi_{\Ecal_\perp}(V - V^{\pi_{\theta},\,\mu})\\
J - J(\pi_{\theta},\mu)
\end{array}\right],\left[\begin{array}{c}
\Pi_{\Ecal_\perp}\bar{G}^V(\theta,V,J,\mu)\\
\bar{G}^J(\theta,J,\mu)
\end{array}\right]\right\rangle \notag\\
&\leq \langle \Pi_{\Ecal_\perp}(V - V^{\pi_{\theta},\,\mu}), \Pi_{\Ecal_\perp}\mathbb{E}_{s\sim \nu^{\pi_{\theta},\,\mu},a\sim\pi_{\theta}(\cdot\mid s),s'\sim\Pcal^{\mu}(\cdot\mid s,a)}[r(s,a,\mu)-J+e_s(e_{s'}-e_s)^{\top}V]\rangle\notag\\
&\hspace{20pt}+c_J \langle J - J(\pi_{\theta},\mu), \mathbb{E}_{s\sim \nu^{\pi_{\theta},\,\mu},a\sim\pi_{\theta}(\cdot\mid s)}[r(s,a,\mu)-J]\rangle\notag\\
&= \langle \Pi_{\Ecal_\perp}(V - V^{\pi_{\theta},\,\mu}), \Pi_{\Ecal_\perp} \mathbb{E}_{s\sim \nu^{\pi_{\theta},\,\mu},a\sim\pi_{\theta}(\cdot\mid s),s'\sim\Pcal^{\mu}(\cdot\mid s,a)}\left[\left(r(s,a,\mu)-J(\pi_{\theta},\mu)+(e_{s'}-e_s)^{\top}\Pi_{\Ecal_\perp}V\right)e_s\right]\rangle\notag\\
&\hspace{20pt}+\langle \Pi_{\Ecal_\perp}(V - V^{\pi_{\theta},\,\mu}), \Pi_{\Ecal_\perp}\mathbb{E}_{s\sim \nu^{\pi_{\theta},\,\mu}}[(J(\pi_{\theta},\mu)-J)e_s]\rangle\notag\\
&\hspace{20pt}+c_J\langle J - J(\pi_{\theta},\mu), \mathbb{E}_{s\sim \nu^{\pi_{\theta},\,\mu},a\sim\pi_{\theta}(\cdot\mid s)}[r(s,a,\mu)-J]\rangle\notag\\
&=\langle \Pi_{\Ecal_\perp}(V - V^{\pi_{\theta},\,\mu}), \Pi_{\Ecal_\perp}\mathbb{E}_{s\sim \nu^{\pi_{\theta},\,\mu},a\sim\pi_{\theta}(\cdot\mid s),s'\sim\Pcal^{\mu}(\cdot\mid s,a)}\left[e_s(e_{s'}-e_s)^{\top}\right]\Pi_{\Ecal_\perp} \left(V-V^{\pi_{\theta},\,\mu}\right)\rangle\notag\\
&\hspace{20pt}+\langle \Pi_{\Ecal_\perp}(V - V^{\pi_{\theta},\,\mu}), \mathbb{E}_{s\sim \nu^{\pi_{\theta},\,\mu}}[(J(\pi_{\theta},\mu)-J)e_s]\rangle- c_J(J - J(\pi_{\theta},\mu))^2\notag\\
&\leq \left(\Pi_{\Ecal_\perp}(V - V^{\pi_{\theta},\,\mu})\right)^{\top} \Pi_{\Ecal_\perp}\mathbb{E}_{s\sim \nu^{\pi_{\theta},\,\mu},a\sim\pi_{\theta}(\cdot\mid s),s'\sim\Pcal^{\mu}(\cdot\mid s,a)}\left[e_s(e_{s'}-e_s)^{\top}\right]\Pi_{\Ecal_\perp} \left(V-V^{\pi_{\theta},\,\mu}\right)\notag\\
&\hspace{20pt}+\frac{\gamma}{2}\|\Pi_{\Ecal_\perp}(V - V^{\pi_{\theta},\,\mu})\|^2+\frac{1}{2\gamma}\|\mathbb{E}_{s\sim \nu^{\pi_{\theta},\,\mu}}[(J(\pi_{\theta},\mu)-J)e_s]\|^2-c_J(J - J(\pi_{\theta},\mu))^2\notag\\
&= \left(\Pi_{\Ecal_\perp}(V - V^{\pi_{\theta},\,\mu})\right)^{\top} \mathbb{E}_{s\sim \nu^{\pi_{\theta},\,\mu},a\sim\pi_{\theta}(\cdot\mid s),s'\sim\Pcal^{\mu}(\cdot\mid s,a)}\left[e_s(e_{s'}-e_s)^{\top}\right]\Pi_{\Ecal_\perp} \left(V-V^{\pi_{\theta},\,\mu}\right)\notag\\
&\hspace{20pt}+\frac{\gamma}{2}\|\Pi_{\Ecal_\perp}(V - V^{\pi_{\theta},\,\mu})\|^2+\frac{1}{2\gamma}\|\mathbb{E}_{s\sim \nu^{\pi_{\theta},\,\mu}}[(J(\pi_{\theta},\mu)-J)e_s]\|^2-c_J(J - J(\pi_{\theta},\mu))^2\notag\\
&\leq -\frac{\gamma}{2}\|\Pi_{\Ecal_\perp}(V - V^{\pi_{\theta},\,\mu})\|^2-\frac{1}{2\gamma}(J - J(\pi_{\theta},\mu))^2,
\end{align*}
where the second inequality follows from the fact that $\langle\vec{a},\vec{b}\rangle\leq \frac{c}{2}\|\vec{a}\|^2+\frac{1}{2c}\|\vec{b}\|^2$ for any vectors $\vec{a},\vec{b}$ and scalar $c>0$, the third inequality applies Lemma~\ref{lem:negative_drift} and the condition $c_J\geq1/\gamma$, the third equation uses the property of the projection matrix $\Pi_{\Ecal_\perp}^2=\Pi_{\Ecal_\perp}=\Pi_{\Ecal_\perp}^{\top}$, and the second equation is a result of the equation below
\[\mathbb{E}_{s\sim \nu^{\pi_{\theta},\,\mu},a\sim\pi_{\theta}(\cdot\mid s),s'\sim\Pcal^{\mu}(\cdot\mid s,a)}\left[\left(r(s,a,\mu)-J(\pi_{\theta},\mu)+(e_{s'}-e_s)^{\top}\Pi_{\Ecal_\perp}V^{\pi_{\theta},\,\mu}\right)e_s\right]=0.\]

Since $\gamma\in(0,1)$, we have $\frac{1}{2\gamma}\geq\frac{\gamma}{2}$. This leads to the claimed result.

\qed

\subsection{Proof of Lemma~\ref{lem:g_markovian}}\label{sec:g_markovian:proof}

The proof of this lemma proceeds in a manner similar to that of Lemma~\ref{lem:f_markovian}. We note that the samples generated in the algorithm follow the time-varying Markov chain
\begin{align}
s_{k-\tau_k} \stackrel{\theta_{k-\tau_k}}{\longrightarrow}  a_{k-\tau_k} \stackrel{\hat{\mu}_{k-\tau_k}}{\longrightarrow}  s_{k-\tau_k+1} \stackrel{\theta_{k-\tau_k+1}}{\longrightarrow} a_{k-\tau_k+1} \stackrel{\hat{\mu}_{k-\tau_k+1}}{\longrightarrow}  \cdots s_{k-1} \stackrel{\theta_{k-1}}{\longrightarrow} a_{k-1} \stackrel{\hat{\mu}_{k-1}}{\longrightarrow}  s_{k}.
\end{align}
We construct an auxiliary Markov chain generated under a constant control
\begin{align}
s_{k-\tau_k} \stackrel{\theta_{k-\tau_k}}{\longrightarrow}  a_{k-\tau_k} \stackrel{\hat{\mu}_{k-\tau_k}}{\longrightarrow}  \widetilde{s}_{k-\tau_k+1} \stackrel{\theta_{k-\tau_k}}{\longrightarrow} \widetilde{a}_{k-\tau_k+1} \stackrel{\hat{\mu}_{k-\tau_k}}{\longrightarrow}  \cdots \widetilde{s}_{k-1} \stackrel{\theta_{k-\tau_k}}{\longrightarrow} \widetilde{a}_{k-1} \stackrel{\hat{\mu}_{k-\tau_k}}{\longrightarrow}  \widetilde{s}_{k}
\label{lem:g_markovian:distribution_auxiliary}
\end{align}

Let $\widetilde{\mu}$ denote the stationary distribution of state, action, and next state under \eqref{lem:g_markovian:distribution_auxiliary}.
We denote $p_k(s,a,s')=\mathbb{P}(s_k=s,a_k=a,s_{k+1}=s')$ and $\widetilde{p}_k(s,a,s')=\mathbb{P}(\widetilde{s}_k=s,\widetilde{a}_k=a,\widetilde{s}_{k+1}=s')$ and define
\begin{align*}
T_1&\triangleq\mathbb{E}[\langle\Delta g_k-\Delta g_{k-\tau_k}, G(\theta_{k},\hat{V}_{k},\hat{J}_{k},\hat{\mu}_{k},s_k,a_k,s_{k+1})-\bar{G}(\theta_{k},\hat{V}_{k},\hat{J}_{k},\hat{\mu}_{k})\rangle],\\
T_2&\triangleq\mathbb{E}[\langle\Delta g_{k-\tau_k}, G(\theta_{k},\hat{V}_{k},\hat{J}_{k},\hat{\mu}_{k},s_k,a_k,s_{k+1})-G(\theta_{k},\hat{V}_{k},\hat{J}_{k},\hat{\mu}_{k},\widetilde{s}_k,\widetilde{a}_k,\widetilde{s}_{k+1})\rangle],\\
T_3&\triangleq\mathbb{E}[\langle\Delta g_{k-\tau_k},G(\theta_{k},\hat{V}_{k},\hat{J}_{k},\hat{\mu}_{k},\widetilde{s}_k,\widetilde{a}_k,\widetilde{s}_{k+1})-\mathbb{E}_{(s,a,s')\sim\widetilde{\mu}}[G(\theta_{k},\hat{V}_{k},\hat{J}_{k},\hat{\mu}_{k},s,a,s')]\rangle]\\
T_4&\triangleq\mathbb{E}[\langle\Delta g_{k-\tau_k},\mathbb{E}_{(s,a,s')\sim\widetilde{\mu}}[G(\theta_{k},\hat{V}_{k},\hat{J}_{k},\hat{\mu}_{k},s,a,s')]-\bar{G}(\theta_{k},\hat{V}_{k},\hat{J}_{k},\hat{\mu}_{k})\rangle].
\end{align*}

It is obvious to see
\begin{align}
\mathbb{E}[\langle\Delta g_k,G(\theta_{k},\hat{V}_{k},\hat{J}_{k},\hat{\mu}_{k},s_k,a_k,s_{k+1})-\bar{G}(\theta_{k},\hat{V}_{k},\hat{J}_{k},\hat{\mu}_{k})\rangle\rangle]=T_1+T_2+T_3+T_4.\label{lem:g_markovian:proof_eq1}
\end{align}

We bound the terms individually. First, we treat $T_1$
\begin{align*}
T_1 &= \mathbb{E}[\langle\Delta g_k-\Delta g_{k-\tau_k}, G(\theta_{k},\hat{V}_{k},\hat{J}_{k},\hat{\mu}_{k},s_k,a_k,s_{k+1})-\bar{G}(\theta_{k},\hat{V}_{k},\hat{J}_{k},\hat{\mu}_{k})\rangle]\notag\\
&\leq \mathbb{E}[\|g_k- g_{k-\tau_k}\|\|G(\theta_{k},\hat{V}_{k},\hat{J}_{k},\hat{\mu}_{k},s_k,a_k,s_{k+1})-\bar{G}(\theta_{k},\hat{V}_{k},\hat{J}_{k},\hat{\mu}_{k})\|]\notag\\
&\hspace{20pt}+\mathbb{E}\Big[\|\bar{G}(\theta_k,\hat{V}_k,\hat{J}_{k},\hat{\mu}_k)-\bar{G}(\theta_{k-\tau_k},\hat{V}_{k-\tau_k},\hat{J}_{k-\tau_k},\hat{\mu}_{k-\tau_k})\|\notag\\
&\hspace{100pt}\cdot\|G(\theta_{k},\hat{V}_{k},\hat{J}_{k},\hat{\mu}_{k},s_k,a_k,s_{k+1})-\bar{G}(\theta_{k},\hat{V}_{k},\hat{J}_{k},\hat{\mu}_{k})\|\Big]\notag\\
&\leq 2B_G\sum_{t=0}^{\tau_k-1}\|g_{k-t}-g_{k-t-1}\|\notag\\
&\hspace{20pt}+2L_G B_G\sum_{t=0}^{\tau_k-1}\left(\|\theta_{k-t}-\theta_{k-t-1}\|+\|\hat{V}_{k-t}-\hat{V}_{k-t-1}\|+|\hat{J}_{k-t}-\hat{J}_{k-t-1}|+\|\hat{\mu}_{k-t}-\hat{\mu}_{k-t-1}\|\right)\notag\\
&\leq 4B_G^2\tau_k\lambda_{k-\tau_k}+2L_G B_G\tau_k(B_F\alpha_{k-\tau_k}+B_G\beta_{k-\tau_k}+B_G\beta_{k-\tau_k}+B_H\xi_{k-\tau_k})\notag\\
&\leq 12 L_G B_F B_G^2 B_H \tau_k\lambda_{k-\tau_k},
\end{align*}
where the second inequality bounds $\|G(\theta_{k},\hat{V}_{k},\hat{J}_{k},\hat{\mu}_{k},s_k,a_k,s_{k+1})-\bar{G}(\theta_{k},\hat{V}_{k},\hat{J}_{k},\hat{\mu}_{k})\|$ by $2B_G$ and $\|\bar{G}(\theta_k,\hat{V}_k,\hat{J}_{k},\hat{\mu}_k)-\bar{G}(\theta_{k-\tau_k},\hat{V}_{k-\tau_k},\hat{J}_{k-\tau_k},\hat{\mu}_{k-\tau_k})\|$ using the Lipschitz continuity established in Lemma~\ref{lem:Lipschitz_operators}. The last inequality follows from the step size condition $\alpha_k\leq\xi_k\leq\beta_k\leq\lambda_k$ for all $k$. The third inequality follows from the fact that $\|g_{k+1}-g_k\|=\lambda_k\|g_k-G(\theta_{k},\hat{V}_{k},\hat{\mu}_{k},s_k,a_k,s_{k+1})\|\leq 2B_G\lambda_k$ for all $k$ and that the per-iteration drift of $\theta_k$, $\hat{V}_k$, and $\hat{\mu}_k$ can be similarly bounded due to Lemma~\ref{lem:boundedness}
\begin{align*}
\|\theta_{k+1}-\theta_k\|\leq B_F\alpha_k,\; \|\hat{V}_{k+1}-\hat{V}_k\|\leq B_G\beta_k,\; |\hat{J}_{k+1}-\hat{J}_{k}|\leq B_G\beta_k,\;\|\hat{\mu}_{k+1}-\hat{\mu}_k\|\leq B_H\xi_k.
\end{align*}

We next bound $T_2$
\begin{align*}
T_2&=\mathbb{E}[\langle\Delta g_{k-\tau_k}, G(\theta_{k},\hat{V}_{k},\hat{J}_{k},\hat{\mu}_{k},s_k,a_k,s_{k+1})-G(\theta_{k},\hat{V}_{k},\hat{J}_{k},\hat{\mu}_{k},\widetilde{s}_k,\widetilde{a}_k,\widetilde{s}_{k+1})\rangle]\notag\\
&\leq 2B_G \mathbb{E}_{\Fcal_{k-\tau_k}}[\mathbb{E}[\|G(\theta_{k},\hat{V}_{k},\hat{J}_{k},\hat{\mu}_{k},s_k,a_k,s_{k+1})-G(\theta_{k},\hat{V}_{k},\hat{J}_{k},\hat{\mu}_{k},\widetilde{s}_k,\widetilde{a}_k,\widetilde{s}_{k+1})\|\mid \Fcal_{k-\tau_k}]]\notag\\
&\leq 2B_G \mathbb{E}[\int_{\Scal}\int_{\Acal}\int_{\Scal}G(\theta_{k},\hat{V}_{k},\hat{J}_{k},\hat{\mu}_{k},s,a,s')\left(p_k(s,a,s')-\widetilde{p}_k(s,a,s')\right)ds\,da\,ds']\notag\\
&\leq 2B_G^2\mathbb{E}[d_{TV}(p_k,\widetilde{p}_k)].
\end{align*}
where the last inequality follows from the definition of TV distance in \eqref{eq:TV_def}.

Applying Lemma B.2 from \citet{wu2020finite}, we then have
\begin{align*}
&T_2\notag\\
&\leq 2B_G^2\mathbb{E}[d_{TV}(p_k,\widetilde{p}_k)]\notag\\
&\leq 2B_G^2\mathbb{E}[d_{TV}(\mathbb{P}(s_k=\cdot),\mathbb{P}(\widetilde{s}_k=\cdot))+\frac{|\Acal|}{2}\|\theta_{k-1}-\theta_{k-\tau_k}\|]\notag\\
&\leq 2B_G^2\mathbb{E}[d_{TV}(\mathbb{P}(s_{k-1}=\cdot),\mathbb{P}(\widetilde{s}_{k-1}=\cdot))+L\|\theta_{k-1}-\theta_{k-\tau_k}\|+L\|\hat{\mu}_{k-1}-\hat{\mu}_{k-\tau_k}\|+\frac{|\Acal|}{2}\|\theta_{k-1}-\theta_{k-\tau_k}\|]\notag\\
&\leq |\Acal|B_G^2\mathbb{E}[\|\theta_{k-1}-\theta_{k-\tau_k}\|]+2L B_G^2 \sum_{t=k-\tau_k}^{k-1}\mathbb{E}[\|\theta_t-\theta_{k-\tau_k}\|+\|\hat{\mu}_t-\hat{\mu}_{k-\tau_k}\|]\notag\\
&\leq (2L+|\Acal|) B_G^2 \tau_k^2(B_F\alpha_{k-\tau_k}+B_H\xi_{k-\tau_k})\notag\\
&\leq (4L+2|\Acal|) B_F B_G^2 B_H\tau_k^2 \lambda_{k-\tau_k},
\end{align*}
where the third inequality is a result of Assumption~\ref{assump:Lipschitz}, and the fourth inequality recursively applies the inequality above it.

The term $T_3$ is proportional to the distance between the distribution of the auxiliary Markov chain \eqref{lem:g_markovian:distribution_auxiliary} at time $k$ and its stationary distribution. To bound $T_3$, 
\begin{align*}
T_3&=\mathbb{E}[\langle\Delta g_{k-\tau_k},G(\theta_{k},\hat{V}_{k},\hat{J}_{k},\hat{\mu}_{k},\widetilde{s}_k,\widetilde{a}_k,\widetilde{s}_{k+1})-\mathbb{E}_{(s,a,s')\sim\widetilde{\mu}}[G(\theta_{k},\hat{V}_{k},\hat{J}_{k},\hat{\mu}_{k},s,a,s')]\rangle]\notag\\
&\leq 2B_G \mathbb{E}_{\Fcal_{k-\tau_k}}[\mathbb{E}[\|G(\theta_{k},\hat{V}_{k},\hat{J}_{k},\hat{\mu}_{k},\widetilde{s}_k,\widetilde{a}_k,\widetilde{s}_{k+1})-\mathbb{E}_{(s,a,s')\sim\widetilde{\mu}}[G(\theta_{k},\hat{V}_{k},\hat{J}_{k},\hat{\mu}_{k},s,a,s')]\|\mid \Fcal_{k-\tau_k}]]\notag\\
&\leq 2B_G\mathbb{E}[\int_{\Scal}\int_{\Acal}\int_{\Scal} G(\theta_{k},\hat{V}_{k},\hat{J}_{k},\hat{\mu}_{k},s,a,s')\left(\widetilde{p}_k(s)-\widetilde{\mu}(s)\right)ds\, da\, ds']\notag\\
&\leq 2B_G^2\mathbb{E}[d_{TV}(\widetilde{p}_k,\widetilde{\mu})]\notag\\
&\leq 2B_G^2\alpha_k,
\end{align*}
where the final inequality follows from the definition of the mixing time $\tau_k$ as the number of iterations for the TV distance between $\widetilde{p}_k$ and $\widetilde{\mu}$ to drop below $\alpha_k$.

Finally, we bound the term $T_4$
\begin{align*}
T_4 &= \mathbb{E}[\langle\Delta g_{k-\tau_k},\mathbb{E}_{(s,a,s')\sim\widetilde{\mu}}[G(\theta_{k},\hat{V}_{k},\hat{J}_{k},\hat{\mu}_{k},s,a,s')]-\bar{G}(\theta_{k},\hat{V}_{k},\hat{J}_{k},\hat{\mu}_{k})\rangle]\notag\\
&\leq 2B_G\mathbb{E}[\|\mathbb{E}_{(s,a,s')\sim\widetilde{\mu}}[G(\theta_{k},\hat{V}_{k},\hat{J}_{k},\hat{\mu}_{k},s,a,s')]-\bar{G}(\theta_{k},\hat{V}_{k},\hat{J}_{k},\hat{\mu}_{k})\|]\notag\\
&\leq 2B_G^2 \mathbb{E}[d_{TV}(\widetilde{\mu},\nu^{\pi_{\theta_k},\,\hat{\mu}_{k}}\otimes\pi_{\theta_k}\otimes\Pcal^{\hat{\mu}_{k}})]\notag\\
&\leq 2 L_{TV} B_G^2\left(\|\pi_{\theta_k}-\pi_{\theta_{k-\tau_k}}\|+\|\hat{\mu}_k-\hat{\mu}_{k-\tau_k}\|\right)\notag\\
&\leq 2 L_{TV} B_G^2\tau_k(B_F\alpha_{k-\tau_k}+B_H\xi_{k-\tau_k})\notag\\
&\leq 4 L_{TV} B_F B_G^2 B_H\xi_{k-\tau_k},
\end{align*}
where the third inequality applies the result in \eqref{eq:TV_stationary_Lipschitz}.

Collecting the bounds on $T_1$-$T_4$ and plugging them into \eqref{lem:g_markovian:proof_eq1}, we get
\begin{align*}
&\mathbb{E}[\langle\Delta g_k, G(\theta_{k},\hat{V}_{k},\hat{\mu}_{k},s_k,a_k,s_{k+1})-\bar{G}(\theta_{k},\hat{V}_{k},\hat{\mu}_{k})\rangle]\notag\\
&=T_1+T_2+T_3+T_4\notag\\
&\leq 12 L_F B_F B_G^2 B_H \tau_k\lambda_{k-\tau_k}
+ (4L+2|\Acal|) B_F B_G^2 B_H\tau_k^2 \lambda_{k-\tau_k} \hspace{-2pt}+\hspace{-2pt} 2B_G^2\alpha_k \hspace{-2pt}+\hspace{-2pt} 4 L_{TV} B_F B_G^2 B_H\xi_{k-\tau_k}\notag\\
&\leq (22L+2|\Acal|) L_F L_{TV} B_F B_G^2 B_H \tau_k^2\lambda_{k-\tau_k}.
\end{align*}

\qed

\section{Details for Example~\ref{example:assumption_Delta}}\label{sec:details_example}

We prove that the mentioned class of MFGs satisfies \eqref{eq:def_Delta} with $\rho=2$ and $\kappa=0$ when $p=1$ and note that a similar line of argument can be made for other $p$.

Equivalent to \eqref{eq:def_Delta} with $\rho=2,\kappa=0,p=0$ is
\begin{align}
J(\pi',\mu^{\star}(\pi))-J(\pi',\mu^{\star}(\pi')) \leq J(\pi,\mu^{\star}(\pi))-J(\pi',\mu^{\star}(\pi)).\label{example:assumption_Delta:eq1}
\end{align}

As the transition kernel does not depend on $\mu$ here, we use $\nu^{\pi}$ to denote the stationary distribution of states under policy $\pi$. Note in this case that $\mu^{\star}(\pi)=\nu^{\pi}$.

We first compute $J(\pi',\mu^{\star}(\pi))$
\begin{align}
J(\pi',\mu^{\star}(\pi)) = \langle \nu^{\pi'},q\sum_{a}\pi'(a\mid \cdot)r(\cdot,a,\nu^{\pi})\rangle=q\sum_{s}\nu^{\pi'}(s)\nu^{\pi}(s).
\end{align}
Similarly, we have
\begin{align*}
J(\pi,\mu^{\star}(\pi))=q\sum_{s}\Big(\nu^{\pi}(s)\Big)^2,\quad J(\pi',\mu^{\star}(\pi'))=q\sum_{s}(\nu^{\pi'}(s)\Big)^2
\end{align*}
As a result,
\begin{align*}
J(\pi',\mu^{\star}(\pi)-J(\pi',\mu^{\star}(\pi'))=q\sum_{s}\nu^{\pi'}(s)\Big(\nu^{\pi}(s)-\nu^{\pi'}(s)\Big),\\
J(\pi,\mu^{\star}(\pi)-J(\pi',\mu^{\star}(\pi))=q\sum_{s}\nu^{\pi}(s)\Big(\nu^{\pi}(s)-\nu^{\pi'}(s)\Big).
\end{align*}
This obvious leads to \eqref{example:assumption_Delta:eq1} as
\[\Big(J(\pi,\mu^{\star}(\pi)-J(\pi',\mu^{\star}(\pi))\Big)-\Big(J(\pi',\mu^{\star}(\pi)-J(\pi',\mu^{\star}(\pi'))\Big)=q\sum_{s}\Big(\nu^{\pi}(s)-\nu^{\pi'}(s)\Big)^2\geq0.\]

\begin{figure*}[!ht]
  \centering
  \includegraphics[width=.5\linewidth]{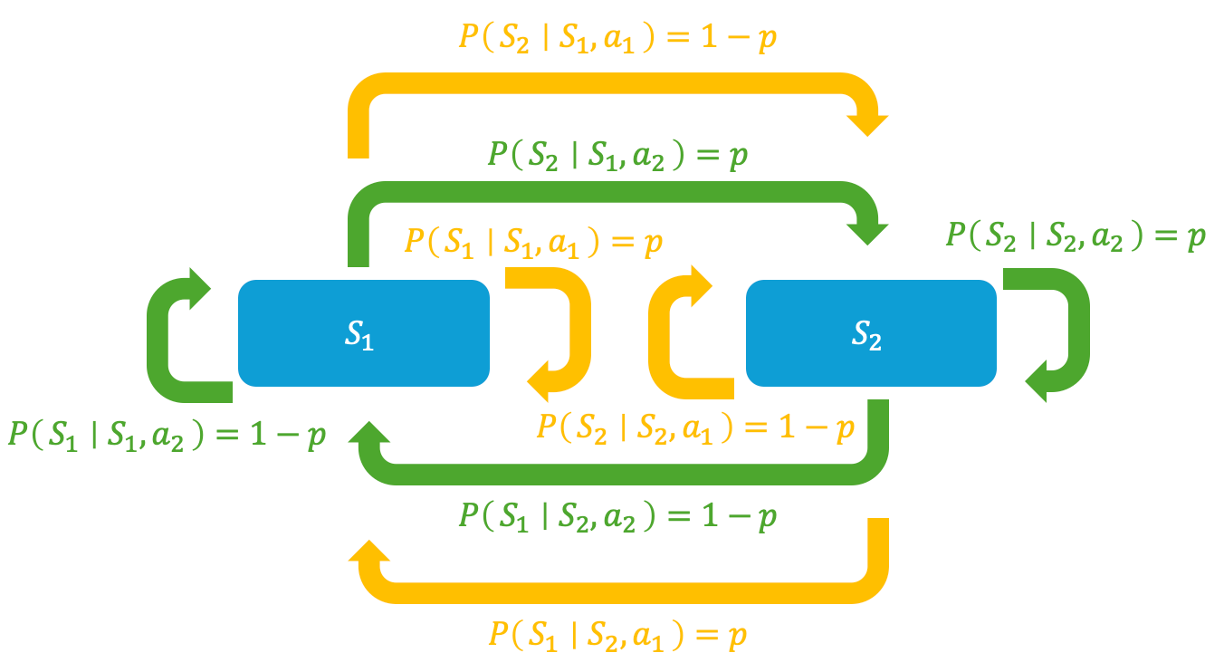}
  
  \caption{Example Mean Field Game Transition}
  \label{fig:example}
\end{figure*}

Next, we provide the detailed derivation on the equilibrium of the MFG in the special case $|\Scal|=|\Acal|=2$ under the transition kernel such that in either state $s\in\{s_1,s_2\}$, the action $a_1$ (resp. $a_2$) leads the next state to $s_1$ (resp. $s_2$) with probability $p=3/4$. A visualization of the transition kernel can be found in Figure.~\ref{fig:example}.

Under any policy $\pi$, the transition matrix is
\begin{align*}
P^{\pi}=\left[\begin{array}{ll}
p\pi(a_1\mid s_1)+(1-p)\pi(a_2\mid s_1) & p\pi(a_1\mid s_2)+(1-p)\pi(a_2\mid s_2) \\
(1-p)\pi(a_1\mid s_1)+p\pi(a_2\mid s_1) & (1-p)\pi(a_1\mid s_2)+p\pi(a_2\mid s_2)
\end{array}\right],
\end{align*}
under which the stationary distribution (induced mean field) is
\[\nu^{\pi}\propto\left[\frac{\pi(a_2\mid s_2)+p-2p\pi(a_2\mid s_2)}{\pi(a_1\mid s_1)+p-2p\pi(a_1\mid s_1)},1\right]^{\top}.\]
In the case $p=3/4$ we have
\[\mu^{\star}(\pi)=\nu^{\pi}=\frac{1}{1+\frac{3/4-\pi(a_2\mid s_2)/2}{3/4-\pi(a_1\mid s_1)/2}}\left[\frac{3/4-\pi(a_2\mid s_2)/2}{3/4-\pi(a_1\mid s_1)/2},1\right]^{\top}.\]

The fact that $\bar{\pi}_1$, $\bar{\pi}_2$, and any policy inducing $[1/2,1/2]^{\top}$ as the mean field can be easily verified at this point.

\section{Average-Reward MDP -- Detailed Formulation and Algorithm}\label{sec:appendix_MDP}

Consider a standard average-reward MDP characterized by state space $\Scal$, action space $\Acal$, transition kernel $\Pcal:\Scal\times\Acal\rightarrow\Delta_{\Scal}$, and reward function $r:\Scal\times\Acal\rightarrow[0,1]$. 
The cumulative reward collected by a policy $\pi:\Scal\rightarrow\Delta_{\Acal}$ is denoted by $J_{\text{MDP}}(\pi)$
\begin{align}
\textstyle J_{\text{MDP}}(\pi)\triangleq\lim_{T\rightarrow\infty}\frac{1}{T}\mathbb{E}_{a_t\sim\pi(\cdot\mid s_t),s_{t+1}\sim\Pcal(\cdot\mid s_t,a_t)}[ \sum_{t=0}^{T-1}r(s_t, a_t) \mid s_0].\label{eq:obj_MDP}
\end{align}

The policy optimization objective under softmax parameterization is
\begin{align}
\max_{\theta}\quad J_{\text{MDP}}(\pi_{\theta}).\label{eq:MDP_obj}
\end{align}

The differential value function under policy $\pi_{\theta}$ is
\begin{align*}
    V_{\text{MDP}}^{\pi_{\theta}}(s) = \mathbb{E}_{a_t\sim\pi_{\theta}(\cdot\mid s_t),s_{t+1}\sim\Pcal(\cdot\mid s_t,a_t)} \Big[ \sum_{t=0}^{\infty}\big(r(s_t, a_t)-J_{\text{MDP}}(\pi) \big) \mid s_0=s\Big].
\end{align*}

We use $P^{\pi}$ and $\nu^{\pi}$ to denote the transition probability matrix and the stationary distribution of states under the control of $\pi$. The policy gradient is 
\begin{align}
\nabla_{\theta} J_{\text{MDP}}(\pi_{\theta}) 
=\mathbb{E}_{s\sim \nu^{\pi_{\theta}},a\sim\pi_{\theta}(\cdot\mid s),s'\sim\Pcal(\cdot\mid s,a)}\Big[(r(s,a)+V_{\text{MDP}}^{\pi_{\theta}}(s')-V_{\text{MDP}}^{\pi_{\theta}}(s))\nabla_{\theta}\log\pi_{\theta}(a\mid s)\Big],\label{eq:policy_grad_MDP}
\end{align}
and $V_{\text{MDP}}^{\pi}$ satisfies the Bellman equation
\begin{align}
V_{\text{MDP}}^{\pi_{\theta}}= \sum_{a}\pi_{\theta}(a\mid \cdot) r(\cdot,a)+J_{\text{MDP}}(\pi_{\theta})\1_{|\Scal|}+(P^{\pi_{\theta}})^{\top}V_{\text{MDP}}^{\pi_{\theta}}.\label{eq:Bellman_MDP}
\end{align}

The algorithm for optimizing $J_{\text{MDP}}$ in an average-reward MDP, simplified from Algorithm~\ref{alg:main}, is presented in Algorithm~\ref{alg:MDP}. We have three main iterates in the algorithm, namely, policy parameter $\theta_k$ and value function estimates $\hat{V}_k$ and $\hat{V}_k$ which are used to track $V_{\text{MDP}}^{\pi_{\theta_k}}$ and $J_{\text{MDP}}(\pi_{\theta_k})$. The policy parameter is updated along the direction of an approximated policy gradient, while the value functions are updated to solve \eqref{eq:Bellman_MDP} and \eqref{eq:MDP_obj} using stochastic approximation.

\begin{algorithm}[!ht]
\caption{Online Actor Critic Algorithm for Average-Reward MDP}
\label{alg:MDP}
\begin{algorithmic}[1]
\STATE{\textbf{Initialize:} policy parameter $\theta_0\in\mathbb{R}^{|\Scal|\times|\Acal|}$, value function estimate $\hat{V}_0\in\mathbb{R}^{|\Scal|},\hat{J}_0\in\mathbb{R}$, gradient/operator estimates $f_0=0\in\mathbb{R}^{|\Scal||\Acal|},g_0^V=0\in\mathbb{R}^{|\Scal|},g_0^J=0\in\mathbb{R}$}
\STATE{\textbf{Sample:} initial state $s_0\in\Scal$ randomly}
\FOR{iteration $k=0,1,2,...$}
\STATE{Take action $a_k \sim \pi_{\theta_k}(\cdot\mid s_k)$. Observe reward $r(s_k,a_k)$ and next state $s_{k+1}\sim \Pcal(\cdot\mid s_k,a_k)$}
\STATE{Policy (actor) update:
\begin{align}
\theta_{k+1} = \theta_k + \alpha_k f_k.\notag
\end{align}
}
\STATE{Value function (critic) update:
\begin{align}
\hat{V}_{k+1} = \Pi_{B_V}(\hat{V}_{k} + \beta_k g_k^V),\quad\hat{J}_{k+1}=\Pi_{[0,1]}(\hat{J}_k+\beta_k g_k^J).\notag
\end{align}
}
\STATE{Gradient/Operator estimate update:
\begin{align}
f_{k+1}&=(1-\lambda_k)f_{k}+\lambda_k(r(s_k,a_k)+\hat{V}_k(s_{k+1}))\nabla\log\pi_{\theta_k}(a_k\mid s_k),\notag\\
g_{k+1}^V&=(1-\lambda_k)g_{k}^V+\lambda_k(r(s_k,a_k)-\hat{J}_k+\hat{V}_k(s_{k+1})-\hat{V}_k(s_k))e_{s_k}\notag\\
g_{k+1}^J&=(1-\lambda_k) g_{k}^J+\lambda_k c_J( r(s_k,a_k)-\hat{J}_k).\notag
\end{align}
}
\ENDFOR
\end{algorithmic}
\end{algorithm}

\section{Simulation Details}\label{sec:simulations_appendix}

We choose the reward function to be
\begin{align*}
r(s,a,\mu)=\mu(s)+\omega_r(s,a)*0.01,\quad\forall s,a,
\end{align*}
where $\omega_r(s,a)\in\mathbb{R}$ is sampled from the standard normal distribution.

For Environments 1 and 2, the transition kernel $\Pcal$ is randomly generated such that for all $s,a$
\begin{align*}
\Pcal^{\mu}(\cdot\mid s,a)\propto\omega_P(s,a),
\end{align*}
where $\omega_P(s,a)\in\mathbb{R}^{|\Scal|}$ is drawn element-wise i.i.d. from the standard uniform distribution.

For Environment 3, the transition kernel $\Pcal$ is also randomly generated such that for all $s,a$
\begin{align*}
\Pcal^{\mu}(\cdot\mid s,a)\propto\omega_P(s,a)+\mu,
\end{align*}
where $\omega_P(s,a)\in\mathbb{R}^{|\Scal|}$ is drawn element-wise i.i.d. from the standard uniform distribution.

For the proposed algorithm, we select the initial step size parameters to be $\alpha_0=10$, $\beta_0=0.1$, $\xi_0=0.02$, and $\lambda_0=1$. The step size parameters for the algorithm in \citet{zaman2023oracle} are taken from the paper in the Numerical Results section. We tried to adjust the parameters of their algorithm in an attempt to see whether we can get it to converge faster, and found out that the parameters prescribed in the paper are good enough and hard to improve at least locally.

\end{document}